\newif\iffinal
\finalfalse	
\finaltrue	

\documentclass[letterpaper,12pt,reqno]{amsart}
\RequirePackage[utf8]{inputenc}
\usepackage[portrait,margin=3cm]{geometry}
\usepackage{blindtext}

\iffinal\else\usepackage[notref,notcite]{showkeys}\fi
\usepackage{mathrsfs,soul}
\usepackage{hyperref}
\usepackage{lipsum}

\usepackage[foot]{amsaddr}
\usepackage{amssymb,amsthm,amsfonts,amsbsy,latexsym,dsfont,amsmath}
\usepackage[textsize=small]{todonotes}       
\usepackage{graphicx}
\usepackage[numeric,initials,nobysame]{amsrefs}
\usepackage{upref,setspace}
\usepackage{xcolor,colortbl}

\usepackage{enumerate}
\newenvironment{enumeratei}{\begin{enumerate}[\upshape (i)]}{\end{enumerate}}
\newenvironment{enumeratea}{\begin{enumerate}[\upshape (a)]}{\end{enumerate}}

\makeatletter
\newcommand{\pushright}[1]{\ifmeasuring@#1\else\omit\hfill$\displaystyle#1$\fi\ignorespaces}
\newcommand{\pushleft}[1]{\ifmeasuring@#1\else\omit$\displaystyle#1$\hfill\fi\ignorespaces}
\makeatother

\usepackage{tikz}
\usetikzlibrary{calc,intersections,through,backgrounds,shapes.geometric}
\usetikzlibrary{graphs}
 	\usepackage{pgfplots}

\usepackage{paralist}

\newenvironment{inparaenumi}{\begin{inparaenum}[\upshape \bfseries (i) ]}{\end{inparaenum}}
\newenvironment{inparaenuma}{\begin{inparaenum}[\upshape \bfseries (a) ]}{\end{inparaenum}}

\newenvironment{inparaenumaa}{\begin{inparaenum}[\upshape(a)]}{\end{inparaenum}}
\newenvironment{inparaenumii}{\begin{inparaenum}[\upshape  (i) ]}{\end{inparaenum}}

\numberwithin{equation}{section}
\numberwithin{figure}{section}
\numberwithin{table}{section}

\usepackage{caption}
\usepackage{subcaption}

\newtheorem{thm}{Theorem}[section]
\newtheorem{lem}[thm]{Lemma}
\newtheorem{theorem}[thm]{Theorem}

\newtheorem{prop}[thm]{Proposition}
\newtheorem{defn}[thm]{Definition}

\newtheorem{ass}[thm]{Assumption}

\newtheorem{conj}[thm]{Conjecture}

\newtheorem{constr}[thm]{Construction}

\theoremstyle{definition}
\newtheorem{rem}{Remark}

\renewcommand{\leq}{\le}
\renewcommand{\geq}{\ge}

\newcommand{\ind}{\mathds{1}}
\newcommand{\eps}{\varepsilon}

\newcommand{\set}[1]{\left\{#1\right\}}

\newcommand{\equald}{\stackrel{\mathrm{d}}{=}}
\newcommand{\probc}{\stackrel{\mathrm{P}}{\longrightarrow}}

\newcommand{\weakc}{\stackrel{\mathrm{d}}{\longrightarrow}}
\newcommand{\convas}{\stackrel{\mathrm{a.s.}}{\longrightarrow}}

\newcommand{\vac}{\mathrm{vac}}
\newcommand{\hght}{\mathrm{ht}}
\newcommand{\anc}{\mathrm{Anc}}
\newcommand{\DF}{\mathrm{DF}}
\newcommand{\lab}{\mathrm{lab}}
\def\qed{ \hfill $\blacksquare$}

\newcommand{\cC}{\mathcal{C}}
\newcommand{\cE}{\mathcal{E}}\newcommand{\cF}{\mathcal{F}}
\newcommand{\cG}{\mathcal{G}}\newcommand{\cI}{\mathcal{I}}
\newcommand{\cL}{\mathcal{L}}

\newcommand{\cP}{\mathcal{P}}\newcommand{\cQ}{\mathcal{Q}}\newcommand{\cR}{\mathcal{R}}
\newcommand{\cS}{\mathcal{S}}\newcommand{\cT}{\mathcal{T}}
\newcommand{\cV}{\mathcal{V}}

\newcommand{\vD}{\mathbf{D}}

\newcommand{\vP}{\mathbf{P}}
\newcommand{\vT}{\mathbf{T}}

\newcommand{\vd}{\mathbf{d}}\newcommand{\ve}{\mathbf{e}}

\newcommand{\vl}{\mathbf{l}}

\newcommand{\vp}{\mathbf{p}}

\newcommand{\vx}{\mathbf{x}}

\newcommand{\mvA}{\boldsymbol{A}}\newcommand{\mvB}{\boldsymbol{B}}

\newcommand{\mvM}{\boldsymbol{M}}

\newcommand{\mvU}{\boldsymbol{U}}
\newcommand{\mvV}{\boldsymbol{V}}\newcommand{\mvX}{\boldsymbol{X}}
\newcommand{\mvZ}{\boldsymbol{Z}}

\newcommand{\mvf}{\boldsymbol{f}}

\newcommand{\mvk}{\boldsymbol{k}}

\newcommand{\mvq}{\boldsymbol{q}}\newcommand{\mvs}{\boldsymbol{s}}
\newcommand{\mvt}{\boldsymbol{t}}\newcommand{\mvu}{\boldsymbol{u}}\newcommand{\mvv}{\boldsymbol{v}}
\newcommand{\mvx}{\boldsymbol{x}}

\newcommand{\mvmu}{\boldsymbol{\mu}}

\newcommand{\mvPi}{\boldsymbol{\Pi}}

\newcommand{\mvxi}{\boldsymbol{\xi}}

\newcommand{\bE}{\mathbb{E}}\newcommand{\bF}{\mathbb{F}}
\newcommand{\bG}{\mathbb{G}}

\newcommand{\bN}{\mathbb{N}}
\newcommand{\bP}{\mathbb{P}}\newcommand{\bR}{\mathbb{R}}
\newcommand{\bT}{\mathbb{T}}

\newcommand{\bZ}{\mathbb{Z}}

\newcommand{\dU}{\mathds{U}}

\newcommand{\sS}{\mathfrak{S}}

\DeclareMathOperator{\shape}{shape}
\DeclareMathOperator{\point}{pt}

\DeclareMathOperator{\E}{\mathbb{E}}
\DeclareMathOperator{\pr}{\mathbb{P}}

\DeclareMathOperator{\var}{Var}
\DeclareMathOperator{\cov}{Cov}
\DeclareMathOperator{\dis}{dis}

\DeclareMathOperator{\GH}{GH}
\DeclareMathOperator{\GHP}{GHP}
 
 \DeclareMathOperator{\height}{ht}
 \DeclareMathOperator{\ord}{ord}
 
 \DeclareMathOperator{\con}{con}

 \DeclareMathOperator{\exec}{exc}
 \DeclareMathOperator{\ERRG}{ERRG}
 
 \DeclareMathOperator{\CM}{CM}
  \DeclareMathOperator{\conf}{Conf}
 \DeclareMathOperator{\Unif}{Unif}

\DeclareMathOperator{\er}{er}
\DeclareMathOperator{\cm}{cm}

\newcommand{\sss}{\scriptscriptstyle}

\newcommand{\erdos}{Erd\H{o}s-R\'enyi }
\newcommand{\ldown}{l^2_{\downarrow}}

\definecolor{aqua}{rgb}{0.0, 1.0, 1.0}
\definecolor{webbrown}{rgb}{.6,0,0}
\definecolor{pinegreen}{rgb}{0.0, 0.47, 0.44}
\definecolor{ultramarineblue}{rgb}{0.25, 0.4, 0.96}
\definecolor{jrnl}{rgb}{0.0, 0.5, 0.0}
\definecolor{lincolngreen}{rgb}{0.11, 0.35, 0.02}
\definecolor{green(html/cssgreen)}{rgb}{0.0, 0.5, 0.0}
\definecolor{airforceblue}{rgb}{0.36, 0.54, 0.66}
\definecolor{azure}{rgb}{0.0, 0.5, 1.0}
\definecolor{bleudefrance}{rgb}{0.19, 0.55, 0.91}
\definecolor{cobalt}{rgb}{0.0, 0.28, 0.67}

\hypersetup{colorlinks=true, linktocpage=true, pdfstartpage=1, pdfstartview=FitV,
	breaklinks=true, pdfpagemode=UseNone, pageanchor=true, pdfpagemode=UseOutlines,
	plainpages=false, bookmarksnumbered, bookmarksopen=true, bookmarksopenlevel=1,
	hypertexnames=true, pdfhighlight=/O,
	urlcolor=webbrown, linkcolor=cobalt, citecolor=jrnl
}

\usepackage{fourier}

\definecolor{aqua}{rgb}{0.0, 1.0, 1.0}
\definecolor{webbrown}{rgb}{.6,0,0}
\definecolor{pinegreen}{rgb}{0.0, 0.47, 0.44}
\definecolor{ultramarineblue}{rgb}{0.25, 0.4, 0.96}
\definecolor{jrnl}{rgb}{0.0, 0.5, 0.0}
\definecolor{lincolngreen}{rgb}{0.11, 0.35, 0.02}
\definecolor{green(html/cssgreen)}{rgb}{0.0, 0.5, 0.0}
\definecolor{airforceblue}{rgb}{0.36, 0.54, 0.66}
\definecolor{azure}{rgb}{0.0, 0.5, 1.0}
\definecolor{bleudefrance}{rgb}{0.19, 0.55, 0.91}
\definecolor{cobalt}{rgb}{0.0, 0.28, 0.67}
\definecolor{darkolivegreen}{rgb}{0.33, 0.42, 0.18}
\definecolor{darkspringgreen}{rgb}{0.09, 0.45, 0.27}
\definecolor{dartmouthgreen}{rgb}{0.05, 0.5, 0.06}
\definecolor{ferngreen}{rgb}{0.31, 0.47, 0.26}
\definecolor{huntergreen}{rgb}{0.21, 0.37, 0.23}

\newcommand{\chnr}[1]{\textcolor{black}{{#1}}}
\newcommand{\chr}[1]{\textcolor{black}{{#1}}}

\begin{document}

\title[Vacant set left by random walks, Wright's constants, and critical random graphs]{Geometry of the vacant set left by random walk on random graphs, Wright's constants, and critical random graphs with prescribed degrees}

\date{}

\author[Bhamidi]{Shankar Bhamidi$^1$}
\address{\hskip-15pt $^1$Department of Statistics and Operations Research, University of North Carolina, Chapel Hill, USA}
\author[Sen]{Sanchayan Sen$^2$}
\address{\hskip-15pt $^2$Department of Mathematics, Indian Institute of Science, Bangalore, India}
\email{bhamidi@email.unc.edu, sanchayan.sen1@gmail.com}

\maketitle
\begin{abstract}
We provide an explicit algorithm for sampling a uniform simple connected random graph with a given degree sequence. 
By products of this central result include:

\vskip1pt

\begin{inparaenumii}
\noindent\item continuum scaling limits of uniform simple connected graphs with given degree sequence and asymptotics for the number of simple connected graphs with given degree sequence under some regularity conditions, and

\vskip1pt

\noindent\item scaling limits for the metric space structure of the maximal components in the critical regime of both the configuration model and the uniform simple random graph model with prescribed degree sequence under finite third moment assumption on the degree sequence.
\end{inparaenumii}

\vskip1pt

As a substantive application we answer a question raised by \v{C}ern{\'y} and Teixeira \cite{cerny-teixeira} by obtaining the metric space scaling limit of maximal components in the  vacant set left by random walks on random regular graphs.

\end{abstract}

\section{Introduction}
\label{sec:intro}
Motivated by applications in a wide array of fields ranging from sociology to systems biology and most closely related to this work, in probabilistic combinatorics and statistical physics, the last few years have witnessed an explosion in both network models as well as interacting particle systems on these models. In this context, the two major themes of this work are as follows:

\medskip

\begin{inparaenumaa}
\noindent\item {\bf Connectivity, percolation and critical random graphs:} A fundamental question in this general area is understanding connectivity properties of the network model, including the time and nature of emergence of the giant component.  Writing $[n] := \set{1,2,\ldots, n}$ for the vertex set, most of these models have a parameter $t$ (related to the edge density) and a model dependent critical time $t_c$ such that for $t< t_c$ (subcritical regime), there exists no giant component (size of the largest component $|\cC_1(t)| = o_P(n)$), while for $t> t_c$ (supercritical regime), the size of the largest component scales like $f(t)n$ with $f(t) > 0$ and is model dependent. Behavior in the so-called ``critical regime'' (i.e., when $t = t_c$) is the main content of this paper.
 To prime the reader let us informally describe the types of results closest in spirit to this work. We defer precise definitions of technical aspects (e.g., definitions of the limiting objects, the underlying topology etc.) to Section \ref{sec:notation} and precise statements of related results to Section \ref{sec:disc}.
The prototypical example of the ``critical'' phenomenon is the \erdos random graph at criticality which is constructed as follows: Fix a parameter $\lambda\in \bR$ and vertex set $[n]$ and let $\ERRG(n,\lambda)$ be the random graph obtained by placing each of the ${n\choose 2}$ possible edges independently with probability $n^{-1} + \lambda n^{-4/3}$.
Maximal component sizes in $\ERRG(n,\lambda)$ were studied extensively in \cite{bollobas1984evolution,luczak1990component,janson1993birth,aldous-crit,nachmias-peres-diameter}. The scaling limit of the maximal components of $\ERRG(n,\lambda)$ when viewed as metric spaces was identified in \cite{BBG-12}.  
It is believed that a large class of random discrete structures, in the critical regime, belong to the ``\erdos universality class.''
Soon after the work  \cite{BBG-12}, an abstract universality principle was developed in \cite{SBSSXW14,bhamidi-broutin-sen-wang} which was used to establish \erdos type scaling limits for a wide array of critical random graph models including the configuration model and various models of inhomogeneous random graphs.
It is strongly believed that the components of critical percolation on high-dimensional tori \cite{hofstad-sapozhnikov,heydenreich-hofstad-1,heydenreich-hofstad-2}, and the hypercube \cite{hofstad-nachmias} also share the \erdos scaling limit, but these problems are open at this point.

\medskip

\noindent\item {\bf Vacant set left by random walk (VSRW) on graphs:} The second main theme is the area of random interlacements and percolative properties of the vacant set of random walks on finite graphs, see e.g. \cite{sznitman2010vacant}. See \cite{vcerny2012random} for a recent survey most closely related to this paper, and \cite{drewitz-rath-sapozhnikov} for an introduction to random interlacements. This question was initially posed by Hilhorst who wanted to understand the geometry of crystals affected by corrosion. The precise mathematical model is as follows: consider a finite  graph on $[n]$ vertices (and to fix ideas assumed connected) which represents the crystalline structure of the object of interest. Now suppose a ``corrosive particle'' wanders through the structure via a simple random walk $\set{X_t:t\geq 0}$ (started from say a uniformly chosen vertex),  marking each vertex it visits as ``corroded'' (this marking does not affect the dynamics of the walk). For a fixed parameter $u\geq 0$, define the vacant set as the set of all vertices that have not been ``corroded'' (i.e., not visited by the walk) by time $u n$,
\begin{equation}
\label{eqn:vac-set-def}
	\cV^u = [n]\setminus \big\{X_t: 0 \leq t \leq n u\big\}.
\end{equation}
When $u$ is ``small'' one expects that only a small fraction of the vertices have been visited by the corrosive particle and thus the maximal connected component $\cC_1(u)$ of the non-corroded set $\cV^u$ has a large connected component of size $\cC_1(u) = \Theta_P(n)$, while if $u$ increases beyond a ``critical point'' $u_{\star}$ then the corrosion in the crystal has spread far enough that the maximal connected component in $\cC_1(u) = o_P(n)$. The ``critical'' $u=u_{\star}$ regime and in particular the fractal properties of connected components in this regime are of great interest.
\end{inparaenumaa}

\subsection{Organization of the paper}
\label{sec:org}
In the remaining subsections of the introduction, we describe the random graph models considered in this paper and give an informal description of our results. Section \ref{sec:res} contains precise statements of our main results regarding scaling limits. Section \ref{sec:sampling-conn-rg} contains the explicit algorithm for generating a uniform connected graph with prescribed degrees (Theorem \ref{lem:alternate-construction}). 
We have deferred major definitions to Section \ref{sec:not}. We discuss the relevance of this work and connections to existing results in Section \ref{sec:disc}.
In Section \ref{sec:proof-lem-alternate-construction}, we present the proof of Theorem \ref{lem:alternate-construction}.
We start Section \ref{sec:proof-plane-tree} with the statement of Lemma \ref{lem:plane-trees} that contains the main technical estimates related to the uniform measure on plane trees with a prescribed child sequence. This lemma forms the crucial work horse in the rest of the proofs. The proof of this lemma occupies the rest of Section \ref{sec:proof-plane-tree}.
Section \ref{sec:proof-conn} and Section \ref{sec:proof-main-deg-vac-thm} contain the proofs of our main results.

%
%

 \subsection{Random graph models}
 \label{sec:rg-models}
  Fix a collection of $n$ vertices labeled by $[n] = \set{1,2,\ldots, n}$ and an associated degree sequence $\vd=\vd^{\sss(n)} = (d_v^{\sss(n)}, v\in [n])$ where $\ell_n := \sum_{v\in [n]}d_v^{\sss(n)}$ is assumed even. There are two natural constructions resulting in generating a random graph on the above vertex set with the prescribed degree sequence.

 \begin{enumeratea}
 	\item {\bf Uniform distributed simple graph:} Let $\bG_{n,\vd}$ denote the space of all simple graphs on $[n]$-labeled vertices with degree sequence $\vd$. Let $\pr_{n,\vd}$ denote the uniform distribution on this space and write $\cG_{n,\vd}$ for the random graph with distribution $\pr_{n, \vd}$.
	\item {\bf Configuration model \cite{Boll-book,molloy1995critical,bender1978asymptotic}:} Recall that a multigraph is a graph where we allow multiple edges and self-loops. Write $\overline\bG_{n,\vd}$ for the space of all multigraphs on vertex set $[n]$ with prescribed degree sequence $\vd$.
Write $\CM_{n}(\vd)$ for the random multigraph 
constructed sequentially as follows: Equip each vertex $v\in [n]$ with $d_v^{\sss(n)}$ half-edges or stubs. Pick two half-edges uniformly from the set of half-edges that have not yet been paired, and pair them to form a full edge. Repeat till all half-edges have been paired. Write $\overline\pr_{n, \vd}$ for the law of $\CM_n(\vd)$.
 \end{enumeratea}

\subsection{Informal description of our contribution}
This work has five major contributions which we now informally describe:

\vskip3pt

\begin{inparaenuma}
\noindent\item  We provide an explicit algorithm (Theorem \ref{lem:alternate-construction}) for sampling a uniform {\bf connected} random graph with given degree sequence by first sampling a planar tree with a modified degree sequence via an appropriately defined tilt with respect to the uniform distribution on the space of trees with this modified degree sequence. This allows us to derive scaling limits for the uniform distribution on the space of simple connected graphs with degree sequence satisfying regularity conditions including a finite number of surplus edges (Theorem \ref{prop:condition-on-connectivity}).

\vskip3pt
	
\noindent\item We then use this result to derive scaling limits for the critical regime of both the configuration model as well as the uniform random graph model with prescribed degree sequence (Theorem \ref{thm:graphs-given-degree-scaling}). 
\chr{Theorem \ref{thm:graphs-given-degree-scaling} is an improvement over the result concerning the scaling limit of the configuration model under critical percolation proved in \cite[Theorem 4.7]{bhamidi-broutin-sen-wang}.
Indeed, the result in \cite[Theorem 4.7]{bhamidi-broutin-sen-wang} follows from Theorem \ref{thm:graphs-given-degree-scaling}.
However, Theorem \ref{thm:graphs-given-degree-scaling} is stronger than \cite[Theorem 4.7]{bhamidi-broutin-sen-wang} as we explain below.}

\chr{
The proof of \cite[Theorem 4.7]{bhamidi-broutin-sen-wang} proceeds in the following steps:
One constructs a process $(\conf_n(t), t\geq 0)$ that is a dynamical version of the configuration model.
(In \cite{bhamidi-broutin-sen-wang}, this process is actually denoted by $\CM_n(\cdot)$, but we use $\conf_n(\cdot)$ to avoid confusing this process with $\CM_n(\vd)$ defined above.)
The scaling limit of this process inside the critical window can be obtained by using a general universality principle developed in \cite{bhamidi-broutin-sen-wang}.
Now, for a supercritical configuration model with degree sequence $(d_1,\ldots,d_n)$, the graph obtained under percolation with edge retention probability $p\in(0,1)$ can be generated as follows: 
(i)
First sample $M\sim\mathrm{Binomial}\big(\sum d_i/2, p\big)$.
Then $M$ has the same distribution as the number of edges remaining after percolation. 
(ii)
Conditional on $M$, uniformly sample $2M$ many half-edges from the set of $\sum d_i$ many half-edges, and construct a configuration model with these half-edges.
This construction enables one to couple the graph obtained under critical percolation on a supercritical configuration model with $\conf_n(t)$ in the critical window. 
From there, one obtains the scaling limit of the percolated graph by comparing it with the scaling limit of $\conf_n(t)$ inside the critical window.
}

\chr{However, we do not see an obvious coupling with the process $\conf_n(\cdot)$ when working with a critical configuration model with a given {\bf deterministic} degree sequence. 
Thus, the above approach fails in this case.
Further, \cite[Theorem 4.7]{bhamidi-broutin-sen-wang} was proved under exponential moment condition on the degree sequence, whereas Theorem \ref{thm:graphs-given-degree-scaling} requires only finite third moment.}
The technique used in this paper is also completely disjoint from \cite{bhamidi-broutin-sen-wang}.

\vskip3pt
	
\noindent\item Write $C(n, n+k)$ for the number of {\bf connected graphs} with $n$ labeled vertices and $n+k$ edges. Deriving asymptotics for $C(n, n+k)$ for fixed $k$ as $n\to\infty$ has inspired a large body work both in the combinatorics community \cite{wright1977number,bender1990asymptotic,spencer-count,flajolet-salvy-schaeffer} as well as in the probability community \cite{aldous-crit,janson2007brownian,takacs1991bernoulli}. Extending such results to count {\bf connected graphs} with {\bf prescribed degree} sequence seems beyond the ken of existing techniques.   As a consequence of our proof technique, we derive asymptotics for such counts (Theorem \ref{thm:number-of-connected-graphs}).  

\vskip3pt

\noindent\item As a substantive application, we answer a question raised by \v{C}ern{\'y} and Teixeira (\cite[Remark 6.1 (2)]{cerny-teixeira} (also see the work of Sznitman, e.g. \cite[Remark 4.5]{sznitman2011lower}) by obtaining the metric space scaling limit of the VSRW on random regular graphs. This is the first result about scaling limit of maximal components in the critical regime for this model. 
\chr{We remark here that we do not see a clear way of applying the techniques of \cite{bhamidi-broutin-sen-wang} to obtain the scaling limit of the VSRW on random regular graphs. 
Indeed, it is easy to check that conditional on the number of edges remaining in the VSRW (run up to some time) on a regular configuration model being $M$, the degree sequence of the vacant set cannot be generated by uniformly sampling $2M$ half-edges from the set of all half-edges.
So we do not see an obvious way of coupling the VSRW with the process $\conf_n(\cdot)$ as explained in {\bf (b)} above.}

The eventual hope, albeit not addressed in this paper is as follows: Consider spatial systems such as the $d$-dimensional lattice or perhaps more relevant to this paper, asymptotics for the vacant set left by random walks on  the $d$-dimensional torus $(\bZ/n\bZ)^d$ in the large $n\to\infty$ network limit.  As described in \cite{sznitman2011lower}, the basic intuition is that for high enough dimensions $d$, the corresponding objects should behave similar to what one sees in the context of $2d$-regular random graphs.

\vskip3pt

\noindent\item \chr{The results of this paper are stepping stones in the study of the minimal spanning tree (MST) of graphs with given degree sequence. In \cite{addarioberry-sen}, using Theorem \ref{thm:graphs-given-degree-scaling} and Theorem \ref{lem:alternate-construction} of this paper, the scaling limit of the MST constructed by assigning exchangeable pairwise distinct weights to the edges of a random (simple) $3$-regular graph and a $3$-regular configuration model is obtained.}

The only other model for which the scaling limit of the MST has been established is the complete graph \cite{AddBroGolMie13}. 
This scaling limit should be universal in the sense that the MST of many standard models exhibiting mean-field behavior should coincide with it (possibly up to some multiplicative factor).
We would eventually like to address the question of universality of the MST of a wide array of models including graphs with given degree sequences. We expect the results of this paper to play a key role in this program.
\end{inparaenuma}

\section{Main results}
\label{sec:res}
We will now describe our main results.  We \chr{describe the general results in Section \ref{sec:abs-res}}. \chr{We then discuss the application regarding vacant sets left by random walks in Section \ref{sec:res-vac}}.  We first fix a convention that we will follow throughout this paper.

\medskip

\noindent{\bf Convention.} For any metric measure space $\mvX=(X,d,\mu)$ and $\alpha>0$, $\alpha\mvX$ will denote the metric measure space $(X,\alpha d,\mu)$, i.e, the space where the metric has been multiplied by $\alpha$ and the measure $\mu$ has remained unchanged. Precise definitions of metric space convergence including the Gromov-Hausdorff-Prokhorov (GHP) topology are deferred to Section \ref{sec:not}.

\subsection{Scaling limits of random graphs with prescribed degrees}
\label{sec:abs-res}
This section describes our main results on graphs with prescribed degree sequence. The first result describes maximal component structure for critical random graphs under appropriate assumptions. For each $n\geq 1$, let $\vd=\vd^{\sss(n)}=\chr{(d_v^{\sss (n)}, v\in[n])}$ be a degree sequence with vertex set $[n]$.
\chr{For simplicity, we will omit the superscript and write $d_v$, $v\in[n]$.}
We will work with degree sequences that satisfy the following assumption.

\begin{ass}
	\label{ass:cm-deg}
	Let $D_n$ be a random variable with distribution given by
\[\bP\big(D_n=i\big)=\frac{1}{n}\#\big\{j\ :\ d_j=i\big\},\]
i.e., $D_n$ has the law of the degree of a vertex selected uniformly at random from $[n]$. Assume the following hold as $n\to\infty$:
\begin{enumeratei}
\item There exists a limiting random variable $D$ with $\pr(D = 1)> 0$ such that $D_n \weakc D$.
\item Convergence of third moments (and hence all lower moments):
\[\E\big[D_n^3\big]:= \frac{1}{n} \sum_{v\in [n]} d_v^3 \to \E\big[D^3\big]<\infty. \]
\item We are in the critical scaling window, i.e., there exists $\lambda \in \bR$ such that
\[\nu_n:= \frac{\sum_{v\in [n]} d_v(d_v-1)}{\sum_{v\in [n]} d_v} = 1+\frac{\lambda}{n^{1/3}} + o(n^{-1/3}). \]
In particular, $\E[D^2]=2\E[D]$.
\end{enumeratei}
\end{ass}

Recall the definitions of the random graphs $\cG_{n,\vd}$ and $\CM_n(\vd)$ from Section \ref{sec:rg-models}.

\begin{thm}[Scaling limit of graphs with given degree sequence]\label{thm:graphs-given-degree-scaling}
Suppose the sequence $\{\vd^{\sss(n)}\}_{n\geq 1}$ satisfies Assumption \ref{ass:cm-deg} with limiting random variable $D$.
\begin{enumeratei}
\item Let $\cC_{(i)}^{\sss n}$ be the $i$-th largest component of $\cG_{n,\vd}$. Endow $\cC_{(i)}^{\sss n}$ with the graph distance and the uniform probability measure on its vertices.  Then there exists a sequence $\mvM^D(\lambda)=(M_1^D(\lambda), M_2^D(\lambda),\ldots)$ of (random) compact metric measure spaces such that
\[\frac{1}{n^{1/3}}\big(\cC_{(1)}^{\sss n}, \cC_{(2)}^{\sss n},\ldots\big)\weakc\mvM^D(\lambda)\]
with respect to product topology induced by GHP distance on each coordinate.
\item The conclusion of part (i) continues to hold with the same limiting sequence $\mvM^D(\lambda)$ if we replace $\cG_{n,\vd}$ by $\CM_n(\vd)$.
\end{enumeratei}
\end{thm}
\begin{rem}
	The limit objects $\mvM^D(\lambda)$ are described explicitly in Construction \ref{constr:M-D}.
The limiting spaces, up to a multiplicative factor, coincide with the scaling limit of the maximal components of the critical \erdos random graph; see \cite[P748]{BBG-limit-prop-11}.
\end{rem}

\begin{rem}
	\chr{Assume the setup of Theorem \ref{thm:graphs-given-degree-scaling}, and let
\[
N_{\sss(i)}^n(\lambda) := 
\big|E\big(\cC_{\sss(i)}^n(\lambda)\big)\big| - \big|\chr{\cC_{\sss(i)}^n}(\lambda)\big| +1
\]
denote the number of surplus edges in $\chr{\cC_{\sss(i)}^n}(\lambda)$. Then it was shown in \cite{dhara-hofstad-leeuwaarden-sen} that there exists a random sequence $\mvZ^D(\lambda)$ such that under the conditions of Theorem \ref{thm:graphs-given-degree-scaling}, as $n \to \infty$,
\[\bigg(\big(n^{-2/3}|\cC_{\sss(i)}^n(\lambda)|, N_{\sss(i)}^n(\lambda)\big),~ i\geq 1\bigg)\weakc \mvZ^D(\lambda).\]
Thus Theorem \ref{thm:graphs-given-degree-scaling} gives further insight into the geometry of these maximal components. The precise theorem and limit object $\mvZ^D(\lambda)$ as constructed in \cite{dhara-hofstad-leeuwaarden-sen} are described in Theorem \ref{thm:cm-component-sizes} below. 
}
\end{rem}

The main ingredient in proving Theorem \ref{thm:graphs-given-degree-scaling} is the following result about the uniform distribution on the space of all {\bf connected} simple graphs with a prescribed degree sequence.
For each fixed $\widetilde m\geq 1$, let $\widetilde\vd^{\sss(\widetilde m)}=(\widetilde d_1^{\sss(\widetilde m)},\hdots,\widetilde d_{\widetilde m}^{\sss(\widetilde m)})$ be a given degree sequence. 
We will often suppress the superscript and write $\widetilde\vd$, $\widetilde d_i$ etc.
Consider the following assumption on the sequence $\{\widetilde\vd^{\sss(\widetilde m)}\}_{\widetilde m\geq 1}$:
\begin{ass}\label{ass:degree}\hfill
\begin{enumeratei}
\item\label{it:1}
$\widetilde d_j\ge 1$ for $1\le j\le\widetilde m$, and $\widetilde d_1=1$.
\item\label{it:2}
There exists a \chr{probability mass function (p.m.f.)} $(\tilde p_1,\tilde p_2,\ldots)$ with
\[\widetilde p_1>0,\quad  
\sum_{i\ge 1}i\widetilde p_i=2,\quad\text{and}\quad
\sum_{i\ge 1}i^2\widetilde p_i<\infty\]
such that
\[\frac{1}{\widetilde m}\#\set{j: \widetilde d_j=i}\to\widetilde p_i\ \text{ for }\ i\ge 1,\ \text{ and }\
\frac{1}{\widetilde m}\sum_{i\ge 1}\widetilde d_i^2\to\sum_{i\ge 1}i^2\widetilde p_i.\]
In particular, $\max_{1\le j\le \widetilde m}\widetilde d_j=o(\sqrt{\widetilde m})$.
\end{enumeratei}
\end{ass}
\begin{rem}\label{rem:d1=1-irrelevant}
We make two observations about the above set of assumptions.
\begin{enumeratei}
\item
The assumption $\widetilde d_1=1$ makes the notation in the proofs simpler. It has no other special relevance. Indeed,
since $\widetilde p_1>0$, a positive proportion of vertices have degree one when $\widetilde m$ is large. Thus, we can always consider the vertex that has the smallest label among all vertices that have degree one.
\item We will work with connected graphs with fixed complexity, i.e., for all $\widetilde m\geq 1$, the degree sequence $\widetilde\vd^{\sss(\widetilde m)}$ will satisfy $\sum_{j\in[\widetilde m]}\widetilde d_j=2(\widetilde m-1)+2k$ for some fixed $k\geq 0$.
    Hence in this case, the assumption $\sum_{i\geq 1}i\widetilde p_i=2$ is redundant as it follows from the other assumptions.
\end{enumeratei}
\end{rem}
Let $\bG_{\widetilde\vd}^{\con}$ be the set of all {\bf connected}, simple, labeled (by $[\widetilde m]$) graphs with degree sequence $\widetilde\vd$ where the vertex labeled $j$ has degree $\tilde d_j$.

\begin{thm}[Scaling limit of {\bf connected} graphs with given degree sequence]\label{prop:condition-on-connectivity}
Consider a sequence of degree sequences $\widetilde\vd^{\sss(\widetilde m)}=(\widetilde d_1,\hdots,\widetilde d_{\widetilde m})$ satisfying Assumption \ref{ass:degree}. In addition, assume that for all $\widetilde m$,
\begin{align}\label{eqn:666}
\sum_{j\in[\widetilde m]}\widetilde d_j=2(\widetilde m-1)+2k
\end{align}
for some (fixed) nonnegative integer $k$.
Sample $\cG_{\widetilde\vd}^{\con}$ uniformly from $\bG_{\widetilde{\vd}}^{\con}$, and endow it with the graph distance and the uniform probability measure on vertices. Then there exists a random compact metric measure space $M^{\sss(k)}$ such that
\[\frac{1}{\sqrt{\widetilde m}}\cG_{\widetilde\vd}^{\con}\weakc \frac{1}{\sigma}M^{\sss(k)}\]
in the GHP sense, where $\sigma^2=\sum_{i\ge 1}i^2\tilde p_i-4$ is the asymptotic variance.
\end{thm}

\begin{rem}
	The limit object $M^{\sss(k)}$ is described explicitly in Construction \ref{constr:M-k}.
\end{rem}

\chr{Part of the conclusion in Theorem \ref{prop:condition-on-connectivity} is that when the degree sequence satisfies \eqref{eqn:666} and Assumption \ref{ass:degree}, $\bG_{\widetilde\vd}^{\con}$ is nonempty for large $\widetilde m$.}
Our next result concerns enumeration of the set $\bG_{\widetilde\vd}^{\con}$.
\begin{thm}[Asymptotic number of connected graphs with given degree sequence when the complexity is fixed]\label{thm:number-of-connected-graphs}
Consider a sequence of degree sequences $\widetilde\vd^{\sss(\widetilde m)}=(\widetilde d_1,\hdots,\widetilde d_{\widetilde m})$ satisfying Assumption \ref{ass:degree}. 
Assume further that for all $\widetilde m$, \eqref{eqn:666} holds for some fixed nonnegative integer $k$. 
Let $\sigma^2=\sum_{i\ge 1}i^2\tilde p_i-4$ be the asymptotic variance. Then
\[\lim_{\widetilde m}\ \frac{\big|\bG_{\widetilde\vd}^{\con}\big|\times
\prod_{i=1}^{\widetilde m}\big(\widetilde d_i-1\big)!\times
\widetilde m^{k/2}}{\big(\widetilde m+2k-2\big)!}=\frac{\sigma^k}{k!}\E\bigg[\bigg(\int_0^1\ve(x)dx\bigg)^k\bigg],\]
where $\big(\ve(x),\ 0\leq x\leq 1\big)$ is a standard Brownian excursion.
\end{thm}
\begin{rem}\label{rem:wright}
Let $C(n,n+k)$ denote the number of connected graphs with $n$ labeled vertices and $n+k$ edges. Wright \cite{wright1977number} showed that for any fixed $k\geq -1$,
\[C(n, n+k)\sim\rho_k n^{n+(3k-1)/2}\qquad \text{as }\ n\to\infty,\]
where the constants $\rho_k$ satisfy a certain recurrence relation. Spencer \cite{spencer-count} proved a connection between this purely combinatorial result and a probabilistic object by showing that
\[\rho_k=\frac{1}{(k+1)!}\E\bigg[\bigg(\int_0^1\ve(x)dx\bigg)^{k+1}\bigg],\ \ k\geq -1,\]
where $(\ve(x),~0\leq x\leq 1)$ is a standard Brownian excursion. Theorem \ref{thm:number-of-connected-graphs} proves the analogue of this result for connected graphs when the degree sequence is fixed.
\end{rem}

\subsection{Geometry of vacant sets left by random walk}
\label{sec:res-vac}
Fix $r\geq 3$ and $n\geq 1$. Here and throughout we assume $nr$ is even.
Recall the definitions of $\bG_{n,\vd}$, $\pr_{n,\vd}$, $\overline\bG_{n,\vd}$, and $\overline\pr_{n,\vd}$ from Section \ref{sec:rg-models}.
Let $\vd_r^{\sss(n)}=(r,r,\ldots,r)$, and
define
\[\bG_{n,r}:=\bG_{n,\vd_r^{\sss(n)}},\ \ \text{ and }\ \ \overline\bG_{n,r}:=\overline\bG_{n,\vd_r^{\sss(n)}}.\]
Analogously write $\pr_{n,r}$ (resp. $\overline\pr_{n,r}$) for $\pr_{n,\vd_r^{\sss(n)}}$ (resp. $\overline\pr_{n,\vd_r^{\sss(n)}}$).
Let $\cG_{n,r}\sim \pr_{n,r}$.
For any (multi)graph $G$, write $P^G$ for the distribution of a simple random walk $(X_t,~ t\geq 0)$ on $G$ with the initial state $X_0$ chosen uniformly at random. Recall the definition of the vacant set from \eqref{eqn:vac-set-def}. Define,
\begin{equation}
\label{eqn:ustar-def}
	u_{\star}=\frac{r(r-1)\ln(r-1)}{(r-2)^2}.
\end{equation}

\begin{figure}[htbp]
\centering
\begin{minipage}{.5\textwidth}
  \centering
  \includegraphics[trim=5.5cm 5cm 3.8cm 4.5cm, clip=true, angle=0, scale=0.9]{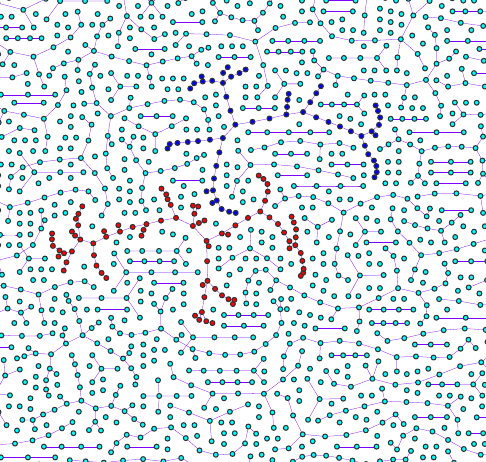}
  \captionof{figure}{Connectivity structure at $u=u_{\star}+0.5$.}
  \label{fig:test1}
\end{minipage}%
\begin{minipage}{.5\textwidth}
  \centering
  \includegraphics[trim=5.7cm 5.5cm 3.6cm 4cm, clip=true, angle=0, scale=0.9]{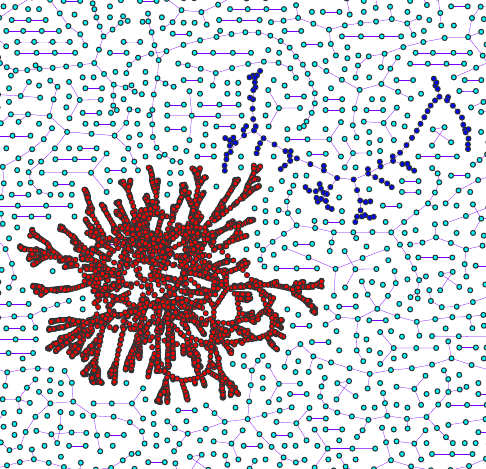}
  \captionof{figure}{Connectivity structure at $u=u_{\star}-0.5$. }
  \label{fig:test2}
\end{minipage}
\end{figure}

Write $\cC_1(u)$ for the maximal connected component in $\cV^u$.  Then the following was shown in \cite{vcerny2011giant}: With $\pr_{n,r}$-probability converging to one as $n\to\infty$,
\begin{enumeratea}
	\item given any $u< u_{\star}$ and $\sigma > 0$, there exist strictly positive constants $\rho, c> 0$ depending only on $u,\sigma,r$ such that
	\[P^{\cG_{n,r}}\big(\big|\cC_{\sss(1)}(u)\big|\geq \rho n\big)\geq 1-cn^{-\sigma};\]
	\item for any fixed $u> u_{\star}$ and $\sigma > 0$, there exists $\rho^\prime >0$ depending only on $u, \sigma, r$ such that
	\[P^{\cG_{n,r}}\big(\big|\cC_{\sss(1)}(u)\big|\geq \rho^\prime \log(n)\big)\leq cn^{-\sigma}.\]
\end{enumeratea}
Figures \ref{fig:test1} and \ref{fig:test2} display the connectivity structure of the vacant set just above and just below $u_{\star}$ respectively where the underlying graph is a $r=4$-regular random graph on $n=50,000$ vertices. The maximal component has been colored red, the second largest component blue, and all other components have been colored cyan.

The main aim of this section is the study of the annealed measures:
\[\vP_{n,r}(\cdot) := \frac{1}{|\bG_{n,r}|} \sum_{G\in \bG_{n,r}} P^G(\cdot)\ \ \text{ and }\ \
\overline\vP_{n,r}(\cdot) := \sum_{G\in\overline\bG_{n,r}}\overline\pr_{n,r}\big(G\big) P^G(\cdot).\]
Building on the work of Cooper and Frieze \cite{cooper-frieze}, \v{C}ern{\'y} and Teixeira in \cite[Theorem 1.1]{cerny-teixeira} showed the following for the above annealed distribution: Let $\{u_n\}_{n\geq 1}$ be any sequence such that there exists fixed $\beta < \infty$ such that, $n^{1/3}|u_n - u_{\star}|\leq \beta$ for all large $n$. Then given any $\eps> 0$, there exists $A:= A(\eps,r,\beta) > 0$ such that for all $n$ large,
\begin{equation}
\label{eqn:cerny-tiex}
	\vP_{n,r}\left(A^{-1} n^{2/3} \leq |\cC_{\sss(1)}\chr{(u_n)}|\leq An^{2/3}\right)\geq 1-\eps.
\end{equation}


They also showed that analogous results hold for $\overline\vP_{n,r}(\cdot)$.  The $n^{2/3}$ scaling of the maximal component size suggests that the critical behavior for this model resembles that of the critical \erdos random graph. Our next theorem confirms this assertion.

\begin{thm}[Scaling limit of the vacant set]\label{thm:vacant-set-scaling}
Let $r\geq 3$.

\begin{inparaenumi}
\noindent\item
Let $\cG_{n,r}\sim\bP_{n,r}$ and $u_\star$ be as in \eqref{eqn:ustar-def}.
Run a simple random walk on $\cG_{n,r}$ up to time $nu_n$ starting from a uniformly chosen vertex, where
\begin{align}\label{eqn:40}
n^{1/3}(u_{\star}-u_n)\to a_0\in\bR.
\end{align}
Let $\cC_{(j)}$ be $j$-th largest component of the subgraph of $\cG_{n,r}$ induced by the vacant set $\cV^{u_n}$. Endow $\cC_{(j)}$ with the graph distance and the uniform probability measure on its vertices.
Then there exists a sequence $\mvM^{\vac}(a_0)=\big(M_1^{\vac}(a_0),M_2^{\vac}(a_0),\ldots\big)$ of random compact metric measure spaces such that under the annealed measure $\vP_{n,r}$,
\[n^{-1/3}\cdot\big(\cC_{(1)}, \cC_{(2)},\ldots\big)\weakc
\mvM^{\vac}(a_0)=\big(M_1^{\vac}(a_0),M_2^{\vac}(a_0),\ldots\big)\]
with respect to product topology induced by GHP distance (see Section \ref{sec:gh-mc} for definition) on each coordinate.

\vskip2pt

\noindent\item
The conclusion in part (i) continues to hold with the same limiting sequence $\mvM^{\vac}(a_0)$ if we replace $\cG_{n,r}$ by $\CM_n(\vd_r^{\sss(n)})$ and $\vP_{n,r}$ by the corresponding annealed measure $\overline\vP_{n,r}$.
\end{inparaenumi}
\end{thm}

\begin{rem}
A complete description of the limiting spaces appearing in Theorem \ref{thm:vacant-set-scaling}
requires certain definitions and is thus deferred to Section \ref{sec:not}. The limiting object $\mvM^{\vac}(a_0)$ is explicitly defined in Construction \ref{constr:M-vac}.   The connection between the scaling limit of the critical \erdos random graph $\ERRG(n,\lambda)$
and the limiting spaces in the results stated in this section is also explained in Section \ref{sec:not}.
\end{rem}

The above result deals with the VSRW on the random $r$-regular graph. We in fact conjecture that for the corresponding problem on random graphs with general prescribed degree sequence, one has analogous results with a {\bf universality} phenomenon under moment conditions on the degree sequence.
\begin{conj}\label{conj:general-degree-vacant}
Let $\vd = \vd^{\sss(n)}=\chr{(d_1^{\sss (n)},\ldots,d_n^{\sss (n)})}$ be a degree sequence, and let $D_n$ denote the degree of a vertex chosen uniformly from $[n]$.
Assume that as $n\to\infty$,  $D_n\weakc D$ with $\E(D^2) < \infty$ and $\E(D_n^2)\to \E(D^2)$. Further assume
\[\nu := \frac{\E[D(D-1)]}{\E[D]} > 1 \mbox{ and } \pr(D \geq 3) > 0.\]
Consider the VSRW on $\cG_{n, \vd}$ or $\CM_n(\vd)$ at level $u$. We conjecture that the following hold:
\begin{enumeratea}
\item There exists a (model dependent) critical point $u_{\star}$ such that for $u< u_{\star}$, size of the maximal component $|\cC_{\sss(1)}(u)| = \Theta_P(n)$ whilst for $u > u_{\star}$,  $|\cC_{\sss(1)}| = \chr{o_P(n)}$.

\item\label{item:erdos} If $~\E(D^3) < \infty$ and $\E(D_n^3)\to E(D^3)$, then for $u_n$ satisfying
	\[\lim_{n\to\infty} n^{1/3} (u_{\star} - u_n)=a_0\]
for some $a_0\in \bR$, the connectivity structure of VSRW at level $u_n$ with edges in the maximal components rescaled by $n^{-1/3}$ satisfy results analogous to Theorem \ref{thm:vacant-set-scaling}.
	
\item {\it (Personal communication from Remco van der Hofstad)} Let $p_k:=\pr(D = k)$, $k=0, 1, \ldots$. Assume that there exists $C>0$ and $\tau\in (3,4)$ such that
\[p_k \sim C k^{-\tau}\ \ \text{ as }\ \ k\to\infty.\]
(In particular, $\E[D^2]<\infty$, but $\E[D^3]=\infty$.)
Then the maximal components in the critical scaling window still belong to the \erdos universality class as in \eqref{item:erdos} above with distances scaling like $n^{-1/3}$. This contrasts drastically with critical percolation on these random graphs where maximal components with distances scaled by $n^{-\frac{\tau-3}{\tau -1}}$ converge to limiting random fractals \cite{bhamidi-hofstad-sen,SB-SD-vdH-SS}.  
 \end{enumeratea}

\end{conj}

\begin{rem}
	\label{rem:clar}
	At this point, we owe the reader two clarifications regarding the conjecture.
First we need to clarify the phrase ``results analogous to Theorem \ref{thm:vacant-set-scaling}'' in (b).
Second we need to explain how the claim in (c) differs from critical bond percolation.
Both of these clarifications are deferred to Section \ref{sec:disc} \eqref{sec:disc-vacant-set}.
\end{rem}

\section{Sampling connected uniform random graphs with prescribed degrees}
\label{sec:sampling-conn-rg}
In this section we describe an explicit algorithm for generating connected random graphs with prescribed degree sequence. This is a core ingredient in the proofs of all the main results.
 We start by setting up notation related to plane trees that will be used both in the statement of the result and throughout the proof sections below.

\subsection{Plane trees: Basic functionals and exploration}
\label{sec:plane-tree-def}
Throughout the sequel we let $\mvt$ denote a plane tree, and use $\rho$ to denote the root. 
Write $\cL(\mvt)$ for the set of leaves of  $\mvt$, i.e., the vertices that have no children.
For each non-root vertex $u\in \mvt$, let $\overleftarrow{u}$ denote the parent of $u$. Let $[\rho, u]$ (resp. $[\rho, u)$) denote the ancestral line of $u$ including (resp. excluding) $u$. Thus $[\rho, u)=[\rho,\overleftarrow{u}]$ \chr{if $u\neq\rho$}. Using the planar embedding, any plane tree $\mvt$ can be explored in a depth-first manner (our convention is to go to the ``leftmost'' child first). Let $\prec_{\DF}$ be the linear order on vertices of a plane tree induced by a depth-first exploration starting from the root, i.e., $x\prec_{\DF} y$ if $x$ is explored \chr{strictly before} $y$ in a depth-first search of the plane tree.

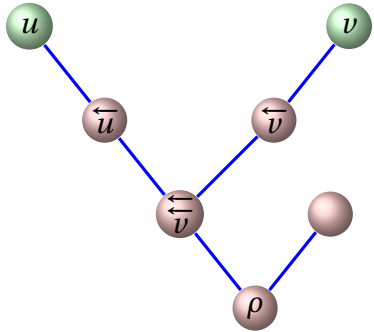
\begin{figure}[htbp]
	\centering
	\pgfmathsetmacro{\nodebasesize}{1} 
	
		\begin{tikzpicture}[scale=.25,  iron/.style={circle, minimum size=6mm, inner sep=0pt, ball color=red!20}, wat/.style= {circle, inner sep=3pt, ball color=green!20}
]

		    \node (1) [iron] at (0,0) {$\rho$};
		    \node (11) [iron] at (-4,5) {$\overleftarrow{\overleftarrow{v}}$};
		    \node (12) [iron] at (4,5) {};
		    \node (111) [iron] at (-8,10) {$\overleftarrow{u}$};
			\node (112) [iron] at (1,10) {$\overleftarrow{v}$};
			\node (1111) [wat] at (-12,15) {${u}$};
			\node (1121) [wat] at (5,15) {${v}$};

	        \draw[blue, very thick] (1)-- (11) -- (111) -- (1111);
	        \draw[blue, very thick] (1) -- (12);
	        \draw[blue, very thick] (11) -- (112) -- (1121);

		\end{tikzpicture}
	\caption{An example of an admissible pair of leaves.  }
	\label{fig:admiss}
\end{figure}

\begin{defn}[Admissible pairs of leaves]
	\label{def:admissible}
	For leaves $u,v\in\cL(\mvt)$, we say that the {\bf ordered} pair $(u,v)$ is {\bf admissible} if
	\chr{$\overleftarrow{v}\neq\rho$, and}
	\[\overleftarrow{\overleftarrow{v}}\in[\rho, \overleftarrow{u}),\ \text{ and }\
	\overleftarrow{u}\prec_{\DF}\overleftarrow{v}.\]
	Let $\mvA(\mvt)$ denote the set of admissible pairs of $\mvt$.
\end{defn}
See Figure \ref{fig:admiss} for an example of an admissible pair.
	We introduce a linear order $\ll$ on $\mvA(\mvt)$ as follows: \chr{For $(u_1,v_1), (u_2,v_2)\in\mvA(\mvt)$, we write} $(u_1,v_1)\ll (u_2,v_2)$ if either $\overleftarrow{u}_1\prec_{\DF}\overleftarrow{u}_2$ or if $\overleftarrow{u}_1=\overleftarrow{u}_2$ and $\overleftarrow{v}_1\prec_{\DF}\overleftarrow{v}_2$. For $u\in\cL(\mvt)$, define
\begin{equation}
\label{eqn:ftu-atu-def}
	\mvA(\mvt, u):=\set{v\in\cL(\mvt)\ :\ (u,v)\in\mvA(\mvt)},\ \text{ and }\
	f_{\mvt}(u):=|\mvA(\mvt, u)|.
\end{equation}
Note that
\begin{equation}
\label{eqn:At=sum-ftu}
|\mvA(\mvt)|=\sum_{u\in\cL(\mvt)} f_{\mvt}(u).
\end{equation}

Now fix $k\ge 1$. Define
\[\mvA_k(\mvt)=\bigg\{\big\{(u_1,v_1),\ldots,(u_k,v_k)\big\}\ \mid\ (u_j,v_j)\in\mvA(\mvt)\text{ and }u_1,v_1,\ldots,u_k,v_k\text{ are }2k\text{ distinct leaves}\bigg\}.\]
Let $\mvA_k^{\ord}(\mvt)$ be the collection of all such ordered $k$-tuples of admissible pairs. 
Clearly,
\[\mvA_1(\mvt)=\mvA(\mvt),\ |\mvA_k^{\ord}(\mvt)|=k!\times|\mvA_k(\mvt)|,\ \text{ and }|\mvA_k^{\ord}(\mvt)|\le |\mvA(\mvt)|^k.\]
For later use, define $\mvA(\mvt)^k = \otimes_{i=1}^k \mvA(\mvt)$ be the $k$-fold Cartesian product of $\mvA(\mvt)$.   

Given a plane tree $\mvt$, write $\mvxi(\mvt) = (\xi_v(\mvt), v\in \mvt)$, where $\xi_v(\mvt)$ is the number of children of $v$ in $\mvt$. Further, write $\mvs(\mvt) = (s_i(\mvt), i\geq 0)$ for the \chr{{\bf child frequency distribution (CFD)}} of $\mvt$. Namely,
 	\[s_i(\mvt):= \#\set{v\in \mvt: \xi_v(\mvt) = i}, \qquad i\geq 0.\]
Given a sequence of integers $\mvs = (s_i,\ i\geq 0)$, we say that the sequence is a \emph{tenable} \chr{CFD} for a tree if there exists a finite plane tree $\mvt$ with $\mvs(\mvt) = \mvs$. It is easy to check $\mvs$ is tenable if and only if $s_i \geq 0$ for all $i$ with $s_0 \geq 1 $, and
\[\sum_{i\ge 0} s_i = 1+ \sum_{i\ge 0} i s_i<\infty.\]
Given a tenable \chr{CFD} $\mvs$, let  $\bT_{\mvs}$ denote the set of all plane trees having \chr{CFD} $\mvs$.

Finally fix $k\geq 1$ and let $\bT_{\mvs}^{\sss(k)}$ denote the set of all pairs $(\mvt, \mvx)$, where $\mvt\in\bT_{\mvs}$ and $\mvx\in\mvA_k(\mvt)$. For a plane tree $\mvt$ and $\mvx=\big\{(u_1,v_1),\ldots,(u_k,v_k)\big\}\in\mvA_k(\mvt)$, let $\cI(\mvt, \mvx)$ be the rooted space obtained by adding an edge between $\overleftarrow{u}_j$ and $\overleftarrow{v}_j$, and deleting $u_j,v_j$ and the two edges incident to them for $j=1,\ldots,k$. (See Figure \ref{fig:admiss-ident} for an illustration.) We endow the space $\cI(\mvt, \mvx)$ with the graph distance and the uniform probability measure on all vertices. Similarly if $\mvt^{\lab}$ is a labeled plane tree (i.e., a plane tree where the vertices are labeled) and $\mvx\in\mvA_k(\mvt^\lab)$, then $\cI(\mvt^\lab, \mvx)$ is the {\bf labeled} graph obtained by following the same construction and retaining the vertex labels.

\begin{figure}[htbp]
	\centering
	\pgfmathsetmacro{\nodebasesize}{1} 
	
		\begin{tikzpicture}[scale=.25,  iron/.style={circle, minimum size=6mm, inner sep=0pt, ball color=red!20}, wat/.style= {circle, inner sep=3pt, ball color=green!20}
]

		    \node (1) [iron] at (0,0) {$\rho$};
		    \node (11) [iron] at (-4,5) {$\overleftarrow{\overleftarrow{v}}$};
		    \node (12) [iron] at (4,5) {};
		    \node (111) [iron] at (-8,10) {$\overleftarrow{u}$};
			\node (112) [iron] at (1,10) {$\overleftarrow{v}$};

	        \draw[blue, very thick] (1)-- (11) -- (111);
	        \draw[blue, very thick] (1) -- (12);
	        \draw[blue, very thick] (11) -- (112);
			\draw[red, very thick] (111) -- (112);

		\end{tikzpicture}
	\caption{An example of the operation $\cI$ applied on the tree $\mvt$ and admissible pair $(u,v)$ in Figure \ref{fig:admiss}.  }
	\label{fig:admiss-ident}
\end{figure}

\subsection{Algorithm for sampling \chr{connected} random graphs with given degree sequence}
\label{sec:algo-rg}

Let $\widetilde{\vd} = \widetilde\vd^{\sss(\widetilde m)}$ be as in Theorem \ref{prop:condition-on-connectivity}, and recall that $\cG_{\widetilde{\vd}}^{\con}$ represents a random connected graph with degree sequence $\widetilde{\vd}$ sampled uniformly from $\bG_{\widetilde{\vd}}^{\con}$. Recall also that under Assumption \ref{ass:degree},  $\widetilde d_1 = 1$. Consider the remaining vertices $\set{2,\ldots, \widetilde m}$, and form the child sequence $\mvxi = \big(\xi_j,\ 2\leq j \leq (\widetilde m+2k)\big)$ via
\begin{equation}\label{eqn:children}
\mvxi:=\big(\widetilde d_2-1,\ldots,\widetilde d_{\widetilde m}-1, 0,\ldots,0\big)\ \ \ (\text{with } 2k\text{ zeros at the end}).
\end{equation}
By the hypothesis of Theorem \ref{prop:condition-on-connectivity},
\begin{align}\label{eqn:18}
\sum_{j=2}^{\widetilde m+2k} \xi_j = (\widetilde m-1+2k) - 1,
\end{align}
and thus $\mvxi$ represents a valid child sequence for a tree on
\begin{align}\label{eqn:def-m}
m:=\widetilde m-1+2k
\end{align}
vertices.
Let $\mvs=\mvs^{\sss(m)}=(s_0, s_1,\ldots)$ be the \chnr{frequency} distribution of $\mvxi$, i.e.,
\begin{align}\label{eqn:35}
s_i=\#\big\{2\leq j\leq \widetilde m+2k :\ \xi_j=i\big\}.
\end{align}

Sample $(\widetilde\cT_{\mvs}, \widetilde\mvX)$ from $\bT_{\mvs}^{\sss(k)}$ uniformly. Assume that $\widetilde\mvX=\{(u_1, v_1),\ldots, (u_k,v_k)\}$, where
\[(u_1, v_1)\ll\ldots\ll (u_k, v_k).\]
Label $u_j$ as $\widetilde m+2j-1$ and $v_j$ as $\widetilde m+2j$, $1\leq j\leq k$. Label the other $\widetilde m-1$ vertices of $\widetilde\cT_{\mvs}$ uniformly using labels $2,\ldots, \widetilde m$ so that in the resulting labeled plane tree $j$ has $\widetilde d_j-1$ many children. (Thus there are $(s_0-2k)!\times\prod_{i\geq 1} s_i!$ many ways of obtaining such a labeling of $\widetilde\cT_{\mvs}$.) 
Call this labeled plane tree $\widetilde\cT_{\mvs}^\lab$. Construct the graph $\cI(\widetilde\cT_{\mvs}^\lab, \widetilde\mvX)$, attach a vertex labeled $1$ to the root, and then forget about the planar order and the root. Let $\cG$ be the resulting graph.
\begin{thm}\label{lem:alternate-construction}
Let $\cG$ be the random graph resulting from the above construction. Then  $\cG\sim\Unif(\bG_{\widetilde\vd}^{\con})$, i.e., $\cG\equald\cG_{\widetilde\vd}^{\con}$.
\end{thm}
\begin{rem}
\chr{Note that
\[
\big|\bT_{\mvs}^{\sss(k)}\big|=\sum_{\mvt\in\bT_{\mvs}}|\mvA_k(\mvt)|
=\sum_{\mvt\in\bT_{\mvs}}|\mvA_k^{\ord}(\mvt)|/k!
=|\bT_{\mvs}|\times\bE\big(|\mvA_k^{\ord}(\cT_{\mvs})|\big)/k!,
\]
where $\cT_{\mvs}\sim\mathrm{Unif}(\bT_{\mvs})$.
It thus follows from Lemma \ref{lem:plane-trees}\eqref{lem:independent-sampling} and the convergence of the second coordinate in Equation \eqref{eqn:nnn} given below that $\bT_{\mvs}^{\sss(k)}$ is nonempty for large $\widetilde m$ when $\widetilde\vd^{\sss(\widetilde m)}$ is as in Theorem \ref{prop:condition-on-connectivity}. 
In particular, it is possible to sample $(\widetilde \cT_{\mvs}, \widetilde X)$ as described above.}
\end{rem}
\begin{rem}
\chr{To implement the above algorithm, we have to sample $(\widetilde \cT_{\mvs}, \widetilde X)$ as described above.
	  Let us briefly discuss how this can be done in polynomial time.
	  We can sample $(\widetilde \cT_{\mvs}, \widetilde X)$ in two steps: 
	  \vskip2pt
	  \noindent(i) Let $\cT_{\mvs}\sim\mathrm{Unif}(\bT_{\mvs})$.
	  Sample $\widetilde \cT_{\mvs}$ according to law
	  \begin{align}\label{eqn:888}
	  \frac{\pr\big(\widetilde \cT_{\mvs}=\mvt\big)}{\pr\big(\cT_{\mvs}=\mvt\big)}
	  =
	  \frac{|\mvA_k(\mvt)|}{\bE\big[|\mvA_k(\cT_{\mvs})|\big]}
	  =: g(\mvt)~, \ \ \mvt\in\bT_{\mvs}.
	  \end{align}
	  \noindent(ii) Conditional on $\widetilde\cT_{\mvs}$, sample $\widetilde X$ uniformly from $\mvA_k(\widetilde\cT_{\mvs})$.
	  \vskip2pt
	  To perform the sampling in (ii), we can list the elements in $\mvA(\widetilde\cT_{\mvs})$ by checking if $(u, v)\in\mvA(\widetilde\cT_{\mvs})$ for every $u,v\in\cL(\widetilde\cT_{\mvs})$. 
	  This requires one to perform $O(m^2)$ many checks. 
	  Turning to (i), we can use rejection sampling to generate $\widetilde\cT_{\mvs}$:
	  First note that $|\mvA_k(\mvt)|\leq m^{2k}$ for any $\mvt\in\bT_{\mvs}$.
	  If we had a lower bound $\bE[|\mvA_k(\cT_{\mvs})|]\geq\eps$ for some $\eps>0$, then it would follow from \eqref{eqn:888} that $g(\mvt)\leq m^{2k}/\eps$.
	  We can generate $\cT_{\mvs}$ using a uniform permutation on $m$ elements as described in Lemma \ref{lem:generating-tree-from-permutation} below.
	  Let $\cT_{\mvs}^{(j)}, j\geq 1$, be i.i.d. copies of $\cT_{\mvs}$, and independent of this sequence sample i.i.d. $\mathrm{Uniform}[0,1]$ random variables $U_j, j\geq 1$.
	  Set $\widetilde\cT_{\mvs}$ equal to $\cT_{\mvs}^{(j_{\star})}$, where 
	  \[
	  j_{\star}=\inf\big\{j\geq 1 : U_j\times m^{2k}\leq \eps\cdot g(\cT_{\mvs}^{(j)})\big\}.	  
	  \]
	  This sampling scheme terminates in $O_P(m^{2k}/\eps)$ many steps.
	  The only caveat here is that one has to find $\eps>0$ that satisfies $\bE[|\mvA_k(\cT_{\mvs})|]\geq\eps$.
	  However, it is often easy to obtain bounds of the form $\bE[|\mvA_k(\cT_{\mvs})|]\geq m^{-r}$ for some integer $r$.
	  For example, if $k=1$, $s_1\geq 2$ and $s_2\geq 1$, then one can show that $\bE[|\mvA_1(\cT_{\mvs})|]\geq m^{-4}$ just by considering admissible pairs of leaves $(u,v)$ for which 
	  \[
	  \xi_{\overleftarrow{u}}(\cT_{\mvs})=1=\xi_{\overleftarrow{v}}(\cT_{\mvs})
	  \qquad \text{ and }\qquad
	  \overleftarrow{\overleftarrow{u}}=\overleftarrow{\overleftarrow{v}}.
	  \]
      (Recall that $\xi_w$ denotes the number of children of $w$.)
      In the case of degree sequences as in Theorem \ref{prop:condition-on-connectivity}, the corresponding child sequence $\mvs$ (as defined in \eqref{eqn:35}) satisfies Assumption \ref{ass:ecd} stated below.
	  In this case it follows from Lemma \ref{lem:plane-trees} \eqref{lem:uniform-integrability}, \eqref{lem:independent-sampling}, and \eqref{lem:conditional-expectation} stated below that $\bE[|\mvA_k(\cT_{\mvs})|]=\Theta(m^{3k/2})$.
	  In particular, the above sampling scheme terminates in $O_P(m^{k/2})$ many steps provided one can find a constant $C$ that satisfies $\bE[|\mvA_k(\cT_{\mvs})|]\geq Cm^{3k/2}$.
	}
  	\end{rem}

\section{Discussion}
\label{sec:disc}

Here we briefly discuss related work, the relevance of the work in this paper and possible extensions and questions raised by this work.

\vskip5pt

\begin{inparaenumaa}
\noindent\item {\bf Graphs with prescribed degree distribution:} Graphs with prescribed degree sequence have played an integral part in probabilistic combinatorics over the last decade and have also been  heavily used in the applied fields including epidemic modeling \cite{britton2007graphs,newman2002spread,keeling2005networks} community detection and clustering \cite{fortunato2010community} and so on. In the context of this paper, the critical point for existence of a giant component was established in \cite{molloy1995critical}. When the degree sequence results in trees, under suitable assumptions on the degree sequence, Broutin and Marckert in \cite{broutin-marckert} showed that these trees appropriately normalized converge to Aldous's continuum random tree; this result will show up in a number of our proofs.

\medskip

\noindent\item {\bf Critical random graphs: } In the context of continuum scaling limits of maximal components in the critical regime, the only other result for the configuration model was derived in \cite{bhamidi-broutin-sen-wang}; here using completely different techniques, critical percolation on the {\bf supercritical} regime of the configuration model where the degree distribution has exponential tails was studied. Associated dynamic versions of this model were constructed and coupled appropriately to Aldous's multiplicative coalescent. A general universality principle also derived in the same paper then resulted in the scaling limits of maximal components at critical percolation. The techniques in that paper, however, do not extend to this work. Here we need start directly with a {\bf critical} prescribed degree sequence; the proof techniques in this paper are completely different and use a combinatorial description of the uniform distribution on the space of {\bf connected} simple graphs with a prescribed degree sequence.

\vskip5pt

\noindent\item\label{sec:disc-vacant-set} {\bf Vacant sets and random interlacements on general random graphs:} With regards to vacant sets, Theorem \ref{thm:vacant-set-scaling} applies to random regular graphs. However as elucidated in Conjecture \ref{conj:general-degree-vacant}, we believe that analogous results hold for the VSRW problem on $\cG_{n, \vd}$ or $\CM_n(\vd)$ constructed using general degree sequence satisfying the hypothesis of Conjecture \ref{conj:general-degree-vacant}. Let us now address the two clarifications described in Remark \ref{rem:clar}. Assuming one can establish the critical point for VSRW for such graphs, we can look at the following two regimes:

\vskip5pt

\noindent{\upshape(i)} Finite third moment: we conjecture that one can construct a (model dependent) random variable $D^*_{\vac}$ (analogous to \eqref{eqn:D-vac-def} for the random regular graph) and $\lambda_{\vac}^*$ a function of both the distribution of $D$ and $a_0$ (analogous to \eqref{eqn:lambda-vac-def}) such that the maximal connected components in the critical regime with edges rescaled by $n^{-1/3}$ converge to $\mvM^{D_{\vac}^*}(\lambda_{\vac}^*)$. This explicates the ``universality'' phenomenon we expect in this regime.

\vskip5pt

\noindent{\upshape(ii)} Infinite third moment regime: This regime seems more non-trivial. Personal communication from Remco van der Hofstad suggests that the limits here differ substantially from critical percolation. Scaling limits for critical percolation in such heavy tailed random graphs were derived in \cite{bhamidi-hofstad-sen}, where scaling limits of maximal components in Aldous's multiplicative coalescent were established in terms of tilted inhomogeneous continuum random trees. One ramification of these  results (\cite[Theorem 1.2]{bhamidi-hofstad-sen}) is the continuum scaling limits of the maximal components in the critical regime of the so-called Norros-Reittu model where the driving weight sequence is assumed to have heavy tails with exponent $\tau\in (3,4)$. These were extended to critical percolation for the configuration model with heavy tailed degree sequence in \cite{SB-SD-vdH-SS}. For a full description of this random graph model as well as the corresponding limits we refer the interested reader to \cite{bhamidi-hofstad-sen}. Calculations by Remco van der Hofstad suggest (see \cite[Section 4.8.4]{hofstadStFlour} for a detailed discussion) that the maximal components in the critical regime for the VSRW model still scale like $n^{-1/3}$ and lie in the \erdos universality regime \cite{BBG-12}. In that sense, random walk percolation differs from bond percolation in this regime.

\vskip5pt
	
	Let us now say a few words on how one can go about proving the above conjecture (at least in the finite third moment setting).  As will become evident from the proofs, the result follows owing to the following three ingredients \begin{inparaenumi}
		\item Theorem \ref{thm:graphs-given-degree-scaling};
		\item a result of Cooper and Frieze \cite{cooper-frieze} which expresses the annealed measure for the vacant set problem in terms of the random graphs with prescribed degree sequence;
		\item refined bounds on the degree sequence of the vacant in the critical scaling window derived in \cite{cerny-teixeira}.
	\end{inparaenumi}
	Parts (i) and (ii) continue to hold for the vacant set problem for random walks on general graphs with prescribed degree sequence. Thus to extend our results to the vacant set problem for general graphs, all one needs is an extension of the refined bounds in (iii) to random walks on general graphs.

%

\medskip

\noindent\item {\bf Proof techniques:} The techniques used in this paper differ from the standard techniques used to show convergence of such random discrete objects to limiting random tree like metric spaces.
One standard technique (used in \cite{BBG-12,SBSSXW14}) is to construct an exploration process of the discrete object of interest that converges to the exploration process of a continuum random tree (see \cite{legall-survey,evans-book} for beautiful treatments), and encode the ``surplus'' edges as a random point process falling under the exploration, and show that this point process converges to a Poisson point process in the limit.
In this work, we use a different technique that requires less work.
We first prove convergence of the object of interest in the Gromov-weak topology, essentially showing that for each fixed $k\geq 2$, the distance matrix constructed from $k$ randomly sampled vertices converges in distribution to the distance matrix constructed from $k$ points appropriately sampled from the limiting structure.
This result, coupled with a global lower mass bound implies via general theory \cite{athreya-lohr-winter} that convergence occurs in the stronger Gromov-Hausdorff-Prokhorov sense.
In the context of critical random graphs, this technique was first used in \cite{bhamidi-hofstad-sen} to analyze the so-called rank-one critical inhomogeneous random graph.
\end{inparaenumaa}


\section{Definitions and limit objects }
\label{sec:not}
This section collects all the basic definitions and existing results in the literature that are needed in the proof for easy reference.

\subsection{Notation and conventions}\label{sec:notation}

\chnr{For any set $A$, we write $|A|$ or $\# A$ for its cardinality and $\ind\set{A}$ for the associated indicator function.
For any graph $H$, we write $V(H)$ and $E(H)$ for the set of vertices and the set of edges of $H$ respectively.
We write $|H|$ for the number of vertices in $H$, i.e., $|H|=|V(H)|$.}
We use the standard Landau notation of $o(\cdot)$, $O(\cdot)$, $\Theta(\cdot)$, and the corresponding \emph{order in probability} notation $o_P(\cdot)$ and $O_P(\cdot)$.
 We use $\probc$, $\weakc$, and \chr{$\convas$} to denote convergence in probability, weak convergence and almost-sure convergence. 
 
 \chnr{Throughout this paper, $C, C', c, c'$ will denote positive universal constants, and their
 	values may change from line to line. Special constants will be indexed as $c_1,c_2$ etc.}

\subsection{Gromov-Hausdorff-Prokhorov metric}
\label{sec:gh-mc}

We mainly follow \cite{EJP2116,AddBroGolMie13,metric-geometry-book}. 
Let us recall the Gromov-Hausdorff distance $d_{\GH}$ between metric spaces.  Fix two metric spaces $X_1 = (X_1,d_1)$ and $X_2 = (X_2, d_2)$. For a subset $C\subseteq X_1 \times X_2$, the distortion of $C$ is defined as
\begin{equation}
	\label{eqn:def-distortion}
	\dis(C):= \sup \set{|d_1(x_1,y_1) - d_2(x_2, y_2)|: (x_1,x_2) , (y_1,y_2) \in C}.
\end{equation}
A correspondence $\cR$ between $X_1$ and $X_2$ is a measurable subset of $X_1 \times X_2$ such that for every $x_1 \in X_1$, there exists at least one $x_2 \in X_2$ such that $(x_1,x_2) \in \cR$ and vice-versa. The Gromov-Hausdorff distance between the two metric spaces  $(X_1,d_1)$ and $(X_2, d_2)$ is defined as
\begin{equation}
\label{eqn:dgh}
	d_{\GH}(X_1, X_2) = \frac{1}{2}\inf \set{\dis(\cR): \cR \mbox{ is a correspondence between } X_1 \mbox{ and } X_2}.
\end{equation}

Suppose $(X_1, d_1)$ and $(X_2, d_2)$ are two metric spaces and $p_1\in X_1$, and $p_2\in X_2$. Then the {\it pointed Gromov-Hausdorff distance} between $\mvX_1:=(X_1, d_1, p_1)$ and $\mvX_2:=(X_2, d_2, p_2)$ is given by
\begin{align}
\label{eqn:dgh-pointed}
	d_{\GH}^{\point}(\mvX_1, \mvX_2) = \frac{1}{2}\inf \set{\dis(\cR): \cR \mbox{ is a correspondence between }X_1 \mbox{ and } X_2\mbox{ and }(p_1, p_2)\in C}.
\end{align}

We will use the Gromov-Hausdorff-Prokhorov distance that also keeps track of associated measures on the corresponding metric spaces. A metric measure space $(X, d , \mu)$ is a metric space $(X,d)$ with an associated finite measure $\mu$ on the Borel sigma algebra on $X$.  Given two metric measure spaces $(X_1, d_1, \mu_1)$ and $(X_2,d_2, \mu_2)$ and a measure $\pi$ on the product space $X_1\times X_2$, the discrepancy of $\pi$ with respect to $\mu_1$ and $\mu_2$ is defined as
\begin{equation}
	\label{eqn:def-discrepancy}
	D(\pi;\mu_1, \mu_2):= ||\mu_1-\pi_1|| + ||\mu_2-\pi_2||
\end{equation}
where $\pi_1, \pi_2$ are the marginals of $\pi$ and $||\cdot||$ denotes the total variation of signed measures. 
Define the function $d_{\GHP}$ as
\begin{equation}
\label{eqn:dghp}
	d_{\GHP}(X_1, X_2):= \inf\bigg\{ \max\bigg(\frac{1}{2} \dis(\cR),~D(\pi;\mu_1,\mu_2),~\pi(\cR^c)\bigg) \bigg\},
\end{equation}
where the infimum is taken over all correspondences $\cR$ and measures $\pi$ on $X_1 \times X_2$. 

Similar to \eqref{eqn:dgh-pointed}, we can define a ``{\it pointed Gromov-Hausdorff-Prokhorov distance}'' $d_{\GHP}^{\point}$ between two metric measure spaces $X_1$ and $X_2$ having two distinguished points $p_1$ and $p_2$ respectively by taking the infimum in \eqref{eqn:dghp} over all correspondences $\cR$ and measures $\pi$ on $X_1 \times X_2$ such that $(p_1, p_2)\in \cR$.

The function $d_{\GHP}$ is a pseudometric that defines an equivalence relation $\sim$: for two metric measure spaces $X$ and $Y$, $X \sim Y \Leftrightarrow d_{\GHP}(X,Y) = 0$. Let $\bar\sS$ be the space of equivalence classes of compact metric measure spaces and $\bar d_{\GHP}$ be the induced metric. Then by \cite{EJP2116}, $(\bar \sS, \bar d_{\GHP})$ is a complete separable metric space. 
Sometimes we will be interested in not just one metric space but an
infinite sequence of metric spaces. 
Then the relevant space will be $\bar\sS^{\bN}$ equipped with the product topology inherited from $\bar d_{\GHP}$.

To ease notation, we will continue to use $d_{\GHP}$ instead of $\bar d_{\GHP}$. 
Similarly we will write $X = (X, d, \mu)$ to denote both the metric space and the corresponding equivalence class.

\subsection{Gromov-weak topology}\label{sec:gromov-weak}
Here we mainly follow \cite{Winter-gromov-weak}. Introduce an equivalence relation on the space of complete and separable metric spaces that are equipped with a probability measure on the associated Borel $\sigma$-algebra by declaring two such spaces $(X_1, d_1, \mu_1)$ and $(X_2, d_2, \mu_2)$ to be equivalent when there exists an isometry $\psi:\mathrm{support}(\mu_1)\to\mathrm{support}(\mu_2)$ such that $\mu_2=\psi_{\ast}\mu_1:=\mu_1\circ\psi^{-1}$, i.e., the push-forward of $\mu_1$ under $\psi$ is $\mu_2$. Write $\sS_{*}$ for the associated space of equivalence classes. As before, we will often ease notation by not distinguishing between a metric space and its equivalence class.

Fix $l\geq 2$, and a complete separable metric space $(X, d)$. Then given a collection of points $\vx:=(x_1, x_2, \ldots, x_l)\in X^l$, let $\vD(\vx):= (d(x_i, x_j))_{i,j\in [l]}$ denote the symmetric matrix of pairwise distances between the collection of points. A function $\Phi\colon \sS_* \to \bR$ is called a polynomial of degree $l$ if there exists a bounded continuous function $\phi\colon \bR_+^{l^2}\to \bR$ such that
\begin{equation}\label{eqn:polynomial-func-def}
	\Phi((X,d,\mu)):= \int  \phi(\vD(\vx)) d\mu^{\otimes l}(\vx).
\end{equation}
Here $\mu^{\otimes l}$ is the $l$-fold product measure of $\mu$. Let $\mvPi$ denote the space of all polynomials on $\sS_*$.

\begin{defn}[Gromov-weak topology]
	\label{def:gromov-weak}
	A sequence $\{(X_n, d_n, \mu_n)\}_{n\geq 1} $ in $\sS_*$ is said to converge to $(X, d, \mu) \in \sS_*$ in the Gromov-weak topology if and only if $\Phi((X_n, d_n, \mu_n))\to \Phi((X, d, \mu))$ for all $\Phi\in \mvPi$.
\end{defn}
In \cite[Theorem 1]{Winter-gromov-weak} it is shown that $\sS_*$ is a Polish space under the Gromov-weak topology. It is also shown that, in fact, this topology can be completely metrized using the so-called Gromov-Prokhorov metric.

\chr{For any metric measure space $(X,d,\mu)$ and $\delta>0$, define
\[\kappa_{\delta}(X)=\kappa_{\delta}(X,d,\mu):=\inf_{x\in X}\bigg\{\mu\big\{y\ :\ d(y,x)\le\delta\big\}\bigg\}.\]
The following theorem gives a criterion for lifting Gromov-weak convergence to Gromov-Hausdorff-Prokhorov convergence.}
\begin{thm}[\cite{athreya-lohr-winter}, Theorem 6.1]\label{thm:gw-to-ghp}
\chr{Suppose $(X,d,\mu)$ and $(X_n,d_n,\mu_n)$, $n\geq 1$, are elements in $\sS_{*}$ such that 
\\
\noindent {\upshape (a)} $\mathrm{support}(\mu_n)=X_n$ for all $n\geq 1$ and $\mathrm{support}(\mu)=X$,
\\
\noindent {\upshape (b)} $(X_n,d_n,\mu_n)\to (X,d_n,\mu)$ with respect to Gromov-weak topology, and
\\
\noindent {\upshape (c)} $\liminf_{n\to\infty}\kappa_{\delta}(X_n)>0$ for all $\delta>0$.}

\chr{Then $(X,d)$ is compact and $(X_n,d_n,\mu_n)\to (X,d,\mu)$ with respect to Gromov-Hausdorff-Prokhorov topology.}
\end{thm}

\subsection{Spaces of trees with edge lengths, leaf weights, and root-to-leaf measures}
\label{sec:space-of-trees}
The rest of this section largely follows \cite{bhamidi-hofstad-sen}.  In the proof of the main results we need the following two spaces built on top of the space of discrete trees. The first space $\vT_{IJ}$ was formulated in \cite{aldous-pitman-edge-lengths,aldous-pitman-entrance} where it was used to study trees spanning a finite number of random points sampled from an inhomogeneous continuum random tree (ICRT). A more general space $\vT_{IJ}^*$ was used in the proofs in \cite{bhamidi-hofstad-sen}.
The index $I$ in $\vT_{IJ}$ and $\vT_{IJ}^*$ is needed for the purpose of keeping track of the number of marked ``hubs,'' i.e., vertices of high (or infinite) degrees in such trees (see \cite{aldous-pitman-edge-lengths,aldous-pitman-entrance,bhamidi-hofstad-sen} for a proper definition). For our purpose it will suffice to consider the case $I=0$. So we only define the space $\vT_{J}:=\vT_{0J}$ and $\vT_{J}^*:=\vT_{0J}^*$.
\medskip

\noindent{\bf The space $\vT_{J}$:} Fix $J\geq 1$. Let $\vT_{J}$ be the space of trees having the following properties:
\begin{enumeratea}
	\item There are exactly $J$ leaves labeled $1+, \ldots, J+$, and the tree is rooted at another labeled vertex $0+$.
	\item Every edge $e$ has a strictly positive edge length $l_e$.
\end{enumeratea}
A tree $\mvt\in \vT_{J}$ can be viewed as being composed of two parts:\\
(1) $\shape(\mvt)$ describing the shape of the tree (including the labels of leaves) but ignoring edge lengths. The set of all possible shapes $\vT_{J}^{\shape}$ is obviously finite for fixed $J$.\\
(2) The edge lengths $\vl(\mvt):= \big(l_e, e\in \mvt\big)$. Consider the product topology on $\vT_{J}$ consisting of the discrete topology on $\vT_{J}^{\shape}$ and the product topology on $\bR^d$.

\medskip

\noindent{\bf The space $\vT_{J}^*$:} We will need a slightly more general space. Along with the two attributes above in $\vT_{J}$, the trees in this space have the following two additional properties. Let $\cL(\mvt):= \big\{1+, \ldots, J+\}$ denote the collection of leaves in $\mvt$. Then every leaf $i+\in \cL(\mvt) $ has the following attributes:

\begin{enumeratea}
	\item[(d)] {\bf Leaf weights:} A nonnegative number \chr{$LW(i+)$.}
	\item[(e)] {\bf Root-to-leaf measures:} A probability measure $\nu_{\mvt,i}$ on the path $[0+,i+]$ connecting the root and the leaf $i+$. 
	Here the path is viewed as a line segment pointed at $0+$ and has the usual Euclidean topology. 
\end{enumeratea}
In addition to the topology on $\vT_{J}$, the space $\vT_{J}^*$ with these additional two attributes inherits the product topology on $\bR^{J}$ owing to leaf weights and $(d_{\GHP}^{\point})^J$ owing to the root-to-leaf measures.

Additionally, we include a special element $\partial$ in $\vT_{J}^*$. This will be useful in the proofs as we will view any rooted tree that does not have exactly $J$ distinct leaves as $\partial$, which will allow us to work entirely in the space $\vT_{J}^*$.
\subsection{Scaling limits of component sizes at criticality}
\label{sec:erdos-scaling-limit}
The starting point for establishing the metric space scaling limit is understanding the behavior of the component sizes.
We first set up some notation. Fix parameters $\alpha, \eta, \beta > 0$, and write $\mvmu = (\alpha, \eta, \beta)\in \bR_+^3$. Let $(B(s),~s\geq 0)$ be a standard Brownian motion. For $\lambda \in \bR$, define
\begin{equation}
\label{eqn:bm-lamb-kapp-def}
W^{\mvmu,\lambda}(s):= \frac{\sqrt{\eta}}{\alpha} B(s)+\lambda s - \frac{\eta s^2}{2\alpha^3},\qquad s\geq 0.
\end{equation}
Write $\overline{W}^{\mvmu,\lambda}$  for the process reflected at zero:
\begin{equation}
\label{eqn:reflected-pro-def}
	\overline{W}^{\mvmu,\lambda}(s) := W^{\mvmu,\lambda}(s) -\min_{0\le u\le s} W^{\mvmu,\lambda}(u),   \qquad s\geq 0.
\end{equation}
 Consider the metric space,
\begin{equation}
\label{eqn:ldown-def}
	\ldown:= \bigg\{\vx= (x_i:i\geq 1): x_1\geq  x_2 \geq \ldots \geq 0, \sum_{i=1}^\infty x_i^2 < \infty\bigg\},
\end{equation}
equipped with the natural metric inherited from $l^2$.
It was shown by Aldous in \cite{aldous-crit} that the excursions of $\overline{W}^{\mvmu,\lambda}$ from zero can be arranged in decreasing order of lengths as
\begin{equation}
\label{eqn:mvxi-def}
	\mvxi^{\mvmu}(\lambda) = \big(|\gamma_{\sss(i)}^{\mvmu}(\lambda)|, i \geq 1\big),
\end{equation}
where $|\gamma_{\sss(i)}^{\mvmu}(\lambda)|$ is the length of the $i$-th largest excursion $\gamma_{\sss(i)}^{\mvmu}(\lambda)$, and further $\mvxi^{\mvmu}(\lambda) \in \ldown$. Let $\cP_\beta$ be a rate $\beta$ Poisson process $\bR_+^2$ independent of $W^{\mvmu,\lambda}(\cdot)$.
For each $i\geq 1$, write $N_{\sss(i)}^{\mvmu}(\lambda)$ for the number of points of $\cP_\beta$ that fall under the excursion $\gamma_{\sss(i)}^{\mvmu}(\lambda)$.


Aldous in \cite{aldous-crit} studied the maximal components of the \erdos random graph in the critical regime and  proved a remarkable result that says that the sizes of the maximal components scaled by $n^{-2/3}$ and the number of surplus edges in the maximal components of $\ERRG(n^{-1}+\lambda n^{-4/3})$ converge jointly in distribution to
$\big(\big(|\gamma_{\sss(i)}^{\mvmu_{\er}}(\lambda)|, N_{\sss(i)}^{\mvmu_{\er}}(\lambda)\big),~ i\geq 1 \big)$,
where $\mvmu_{\er} = (1,1,1)$.
This result has since been generalized to a number of other random graph models. In the context of graphs with given degree sequence,  Nachmias and Peres \cite{nachmias-peres-random-regular} studied critical percolation on random regular graphs;  Riordan \cite{riordan2012phase} analyzed the configuration model with bounded degrees; Joseph \cite{joseph2014component} considered i.i.d. degrees.
A stronger result under finite third moment assumptions was obtained in \cite{dhara-hofstad-leeuwaarden-sen}. We will state a weaker version of this result next.

\begin{thm}[\cite{dhara-hofstad-leeuwaarden-sen}]\label{thm:cm-component-sizes}
Consider a degree sequence $\vd=\vd^{\sss(n)}$ satisfying Assumption \ref{ass:cm-deg} with the limiting random variable $D$ and define $\sigma_r:=\E[D^r]$, $r=1,2,3$.
Write $\chr{\cC_{\sss(i)}^n}(\lambda)$ for the $i$-th largest connected component of $\ \CM_n(\vd)$ (or $\cG_{n,\vd}$). Let
\[
N_{\sss(i)}^n(\lambda) := 
\big|E\big(\cC_{\sss(i)}^n(\lambda)\big)\big| - \big|\chr{\cC_{\sss(i)}^n}(\lambda)\big| +1
\]
denote the number of surplus edges in $\chr{\cC_{\sss(i)}^n}(\lambda)$. Then as $n \to \infty$,
\[\bigg(\big(n^{-2/3}|\cC_{\sss(i)}^n(\lambda)|, N_{\sss(i)}^n(\lambda)\big),~ i\geq 1\bigg)\weakc \mvZ^D(\lambda):= \big(\big(|\gamma_{\sss(i)}^{\mvmu_{D}}(\lambda)|, N_{\sss(i)}^{\mvmu_{D}}(\lambda)\big),~ i\geq 1 \big)
\]
with respect to product topology. Here $\mvmu_{D} = (\alpha_D, \eta_D, \beta_D)$ is given by
\begin{equation*}
\alpha_D = \sigma_1,\ \
\eta_D = \sigma_3\sigma_1-\sigma_2^2,\ \text{ and }\
\beta_D = 1/\sigma_1.
\end{equation*}
\end{thm}
This result, in a stronger form, can be found in \cite[Theorem 2 and Remark 5]{dhara-hofstad-leeuwaarden-sen}.
We will use this result in the next section to describe the limiting metric measure spaces arising in Section \ref{sec:res}.


\subsection{The limiting metric measure spaces}

A compact metric space $(X,d)$ is called a \emph{real tree} \cite{legall-survey,evans-book} if between every two points there is a unique geodesic such that this path is also the only non self-intersecting path between the two points. Functions encoding excursions from zero can be used to construct such metric spaces via a simple procedure. We describe this construction next.

For $0 < a < b <\infty$, an \emph{excursion} on $[a,b]$ is a continuous function $h \in C([a,b], \bR)$ with $h(a)=0=h(b)$ and $h(t) > 0$ for $t \in (a,b)$. The length of such an excursion is $b-a$. For $l \in(0,\infty)$, let $\cE_l$ be the space of all excursions on the interval $[0,l]$. Given an excursion $h \in \cE_l$, one can construct a real tree as follows. Define the pseudo-metric $d_h$ on $[0,l]$:
\begin{equation}
\label{eqn:d-pseudo}
	d_h(s,t):= h(s) + h(t) - 2 \inf_{u \in [s,t]}h(u), \; \mbox{ for } s,t  \in [0,l].
\end{equation}
Define the equivalence relation $s \sim t \Leftrightarrow d_h(s,t) = 0$. Let $[0,l]/\sim$ denote the corresponding quotient space and consider the metric space $\cT_h:= ([0,l]/\sim, \bar d_h)$, where $\bar d_h$ is the metric on the equivalence classes induced by $d_h$. Then $\cT_h$ is a real tree (\cite{legall-survey,evans-book}).
Let $q_h:[0,l] \to \cT_h$ be the canonical projection and write $\mu_{\cT_h}$ for the push-forward of the Lebesgue measure on $[0,l]$ onto $\cT_h$ via $q_h$. Further, we assume that $\cT_h$ is rooted at $\rho := q_h(0)$.  Equipped with $\mu_{\cT_h}$, $\cT_h$ is now a rooted compact metric measure space. Note that by construction, for any $x\in \cT_h$, the function $h$ is constant on $q_h^{-1}(x)$. Thus for each $x\in [0,l]$, we write $\hght(x) = h(q_h^{-1}(x))$ for the height of this vertex.

The Brownian continuum random tree defined below is a fundamental object in the literature of random real trees.

\begin{defn}[Aldous's continuum random tree \cite{Aldo91a}]\label{def:aldous-crt}
	Let $\ve$ be a standard Brownian excursion on $[0,1]$. Construct the random compact real tree $\cT_{2\ve}$ as in \eqref{eqn:d-pseudo} with $h=2\ve$. The associated measure $\mu_{\cT_{2\ve}}$ is supported on the collection of leaves of $\cT_{2\ve}$ almost surely.
\end{defn}
Write $\nu$ for the law of a standard Brownian excursion on the space of excursions on $[0,1]$ namely $\cE_1$.
For $k\geq 0$, let $\widetilde \ve_{(k)}$ be a random excursion with distribution $\tilde \nu_{k}$ given via the following Radon-Nikodym density with respect to $\nu$:
\begin{equation}\label{eqn:tilde-nu-k-def}
\frac{d \tilde \nu_{k} }{d\nu}(h)
= \frac{\left[\int_0^{1}h(u) du\right]^k}{\E\left[\left(\int_0^{1} \ve(u)du\right)^k\right]}
= \frac{\big(\int_{\cT_h}\hght(x)\mu_{\cT_h}(dx)\big)^k}{\E\big[\int_{\cT_{\ve}}\hght(x)\mu_{\cT_{\ve}}(dx)\big]^k}, \qquad h\in \cE_{1}.
	\end{equation}

\begin{constr}[The space $M^{\sss(k)}$]\label{constr:M-k}
Fix $k\geq 0$.
\begin{enumeratea}
\item Let $\widetilde \ve_{(k)}$ be as above, and write $\cT^\star=\cT_{2\widetilde\ve_{(k)}}$. Let $\mu_{\cT^\star}$ denote the associated measure.
\item Conditional on $\cT^\star$, sample $k$ leaves $\set{x_i: 1\leq i\leq k}$ in an i.i.d. fashion from $\cT^\star$ with density proportional to $\hght(x) \mu_{\cT^\star}(dx)$.
\item Conditional on the two steps above, for each of the sampled leaves $x_i$, sample a point $y_i$ uniformly at random on the line $[\rho, x_i]$. Identify $x_i$ and $y_i$, i.e., introduce the equivalence relation $x_i\sim y_i$, $1\leq i\leq k$, and form the quotient space $\cT_{\star}/\sim$.
\end{enumeratea}
Set $M^{\sss(k)}$ to be the resultant (compact) random metric measure space.
\end{constr}

Next recall the definition of $\mvZ^D(\lambda)$ from Theorem \ref{thm:cm-component-sizes}.
\begin{constr}[The sequence $\mvM^D(\lambda)$]\label{constr:M-D}
\hfill
\begin{enumeratea}
\item Sample $\mvZ^D(\lambda)=\big(\big(|\gamma_{\sss(i)}^{\mvmu_{D}}(\lambda)|, N_{\sss(i)}^{\mvmu_{D}}(\lambda)\big),\ i\geq 1\big)$.
For simplicity, write
\[\xi_i=|\gamma_{\sss(i)}^{\mvmu_{D}}(\lambda)|,\ \text{ and }\ N_i=N_{\sss(i)}^{\mvmu_{D}}(\lambda).\]
\item Conditional on $\mvZ^D(\lambda)$, construct the spaces $S_i$ independently for $i\geq 1$, where $S_i\equald M^{(N_i)}$.
\end{enumeratea}
Set
\[\mvM^D(\lambda)=\big(M_1^D(\lambda), M_2^D(\lambda),\ldots\big),\ \text{ where }\ M_i^D(\lambda)=\frac{\alpha_D\sqrt{\xi_i}}{\sqrt{\eta_D}}\cdot S_i,\quad\ i\geq 1.\]
\end{constr}

Note that the sequence $\mvM^D(\lambda)$ of limiting spaces depends only on the first three moments of the random variable $D$ (which is also true for $\mvZ^D(\lambda)$--the scaling limit of the component sizes and the number of surplus edges).

Finally, let $a_0$ be as in theorem \ref{thm:vacant-set-scaling}. Define
\begin{align}\label{eqn:lambda-vac-def}
\lambda_{\vac}=\frac{a_0(r-2)^2}{r(r-1)},\ \text{ and }\ p_{\vac}=\exp\bigg(-\frac{r\ln(r-1)}{(r-2)}\bigg),
\end{align}
and let $D_{\vac}$ be the mixture random variable
\begin{align}\label{eqn:D-vac-def}
D_{\vac}=(1-p_{\vac})\cdot\delta_0+p_{\vac}\cdot\mathrm{Binomial}\bigg(r,\ \frac{1}{r-1}\bigg).
\end{align}
\begin{constr}[The sequence $\mvM^{\vac}(a_0)$]\label{constr:M-vac}
Set
\[\mvM^{\vac}(a_0):=\mvM^{D_{\vac}}\big(\lambda_{\vac}\big).\]
\end{constr}

\begin{rem}\label{rem:erdos-renyi}
The \erdos scaling limit identified in \cite{BBG-12, BBG-limit-prop-11} can be recovered by taking the limiting random variable to be $D_{\er}\sim\mathrm{Poisson}(1)$, i.e., the scaling limit of $\ERRG(n^{-1}+\lambda n^{-4/3})$ (after rescaling the graph distance by $n^{-1/3}$) is given by
\[\mvM_{\er}(\lambda):=\mvM^{D_{\er}}(\lambda).\]
(Note that in this case, $\alpha_{D_{\er}}=\eta_{D_{\er}}=\beta_{D_{\er}}=1$.) The result for $\ERRG(n^{-1}+\lambda n^{-4/3})$ can be obtained from Theorem \ref{thm:graphs-given-degree-scaling} by observing the following two facts:
\begin{enumeratei}
\item The (random) degree sequence of $\ERRG(n^{-1}+\lambda n^{-4/3})$ satisfies Assumption \ref{ass:cm-deg} with limiting random variable $D_{\er}$.
\item Conditional on the event where the degree sequence equals $\vd$, $\ERRG(n^{-1}+\lambda n^{-4/3})$ is uniformly distributed over $\bG_{n,\vd}$.
\end{enumeratei}
\end{rem}

\section{Proof of Theorem \ref{lem:alternate-construction}}\label{sec:proof-lem-alternate-construction}
 \chr{Since Theorem \ref{lem:alternate-construction} does not require any of the ingredients required for the remaining theorems, we start by giving a quick proof of this result. } Fix a graph $G\in\bG_{\widetilde{\vd}}^{\con}$. Root the graph at the only neighbor of $1$ (recall that $\tilde d_1=1$), and remove the vertex $1$ and the edge incident to it. Suppose $H$ is the resulting rooted, labeled graph. We can construct a labeled plane tree from $H$ in the following way:
\begin{enumeratei}
\item Call the root $u_1$. Set the status of all its neighbors as ``discovered,'' and set the status of $u_1$ as ``explored.'' Shuffle all its neighbors uniformly and go to the ``leftmost'' neighbor and call it $u_2$.
\item When we are at $u_r$ ($r\ge 2$), search for all its neighbors (other than $u_{r-1}$) in the graph at that time. If none of these neighbors have been discovered previously, then shuffle them uniformly, set their status as ``discovered,'' set the status of $u_r$ as ``explored,'' and go to the leftmost neighbor and call it $u_{r+1}$.

If some of these neighbors have been previously discovered, then these edges create surplus. Suppose we have found $\chr{\ell_{0,r}}$ many surplus edges before exploring $u_r$, and at $u_r$ we found $\chr{\ell_{1,r}}$ many new surplus edges $e_1,\ldots,\chr{e_{\ell_{1,r}}}$.
Assume that $e_j=(u_r, y_j)$ and $y_1\prec_{\DF}\ldots\prec_{\DF}y_{\chr{\ell_{1,r}}}$. 
For $j=1,\ldots,\chr{\ell_{1,r}}$, delete the edge $e_j$, and create two leaves labeled $\widetilde m+2\chr{\ell_{0,r}}+2j-1$ and $\widetilde m+2\chr{\ell_{0,r}}+2j$, where $u_r=\overleftarrow{\widetilde m+2\chr{\ell_{0,r}}+2j-1}$ (i.e., $u_r$ is the parent of the leaf labeled $\widetilde m+2\chr{\ell_{0,r}}+2j-1$), and similarly $y_j=\overleftarrow{\widetilde m+2\chr{\ell_{0,r}}+2j}$.
Shuffle the neighbors of $u_r$ uniformly (including the newly created leaves), set their status as ``discovered,'' set the status of $u_r$ as ``explored,'' and move to the leftmost neighbor of $u_r$ and call it $u_{r+1}$. (Note that we do {\bf not} set the status of $\widetilde m+2\chr{\ell_{0,r}}+2j$, $j=1,\ldots,\chr{\ell_{1,r}}$, as discovered at this point.)

    If $u_r$ has no neighbors other than $u_{r-1}$, then go to the next (in the depth-first order) discovered but unexplored vertex and call it $u_{r+1}$.
\end{enumeratei}
Let $\mvt_{\ast}^{\lab}$ be the resulting labeled plane tree, and set
$\mvx_\ast=\{(\widetilde m+1,\widetilde m+2),\ldots,(\widetilde m+2k-1,\widetilde m+2k)\}$.
Note that $(\widetilde m+2j-1, \widetilde m+2j)$ is an admissible pair for $1\leq j\leq k$. 
Note also that the child sequence of $\mvt_\ast^{\lab}$ is the sequence $\mvxi$ defined in \eqref{eqn:children}. 
Thus $(\mvt_\ast^\lab,\mvx_\ast)$ is a random, labeled  element of $\bT_{\mvs}^{(k)}$. 
Let $\DF(H)$ be the set of all possible realizations of $(\mvt_\ast^\lab,\mvx_\ast)$.

\chr{In the above exploration process, when we are at the vertex labeled $j\in\{2,\ldots,m\}$, the shuffling can be done in exactly $(\tilde d_j-1)!$ many ways.
Since the vertices labeled $\widetilde m+1,\ldots,\widetilde m+2k$ all have degree one, when we are at any of these vertices, we simply proceed to the next step.
Further, each one of the $\prod_{j=2}^{\widetilde m}(\tilde d_j-1)!$ combinations of shufflings results in a unique realization of $(\mvt_{\ast}^{\lab},\mvx_\ast)$.
It thus follows that
\begin{align}\label{eqn:997}
	\big|\DF(H)\big|=\prod_{j=2}^{\widetilde m}(\tilde d_j-1)!.
\end{align}
Note also that for every $(\mvt^\lab,\mvx)\in\DF(H)$,
\begin{align}\label{eqn:998}
\bP\big((\widetilde\cT_{\mvs}^\lab, \widetilde\mvX)=(\mvt^\lab,\mvx)\big)=\frac{1}{|\bT_{\mvs}^{\sss(k)}|}
\times\frac{1}{(s_0-2k)!\times\prod_{i\geq 1} s_i!},
\end{align}
where $(\widetilde\cT_{\mvs}^\lab, \widetilde\mvX)$ are as defined before Theorem \ref{lem:alternate-construction}.
}

\chr{Now, clearly $\cI(\mvt^\lab,\mvx)=H$ for every $(\mvt^\lab,\mvx)\in\DF(H)$.
Conversely, if $(\mvt^\lab,\mvx)$ is a labeled element of $\bT_{\mvs}^{(k)}$ satisfying 
(a) the vertices of $\mvt^\lab$ are labeled $2,\ldots,\widetilde m+2k$,
(b) $\mvx=\{(\widetilde m+1, \widetilde m+2), \ldots, (\widetilde m+2k-1, \widetilde m+2k)\}$ and 
$(\widetilde m+1, \widetilde m+2)<<\ldots<<(\widetilde m+2k-1, \widetilde m+2k)$, and finally
(c) $\cI(\mvt^\lab,\mvx)=H$,
then $(\mvt^\lab,\mvx)\in\DF(H)$.
Indeed, in each step of the above exploration procedure, if the shuffling agrees with the planar order in $\mvt^\lab$, then for this combination of shufflings the resulting realization of $(\mvt_\ast^\lab,\mvx_\ast)$ will be $(\mvt^\lab,\mvx)$.
}
Thus
\begin{align}\label{eqn:34}
\bP\big(\cI(\widetilde\cT_{\mvs}^\lab, \widetilde\mvX)=H\big)
&=\sum_{(\mvt^\lab,\mvx)\in\DF(H)}\bP\big((\widetilde\cT_{\mvs}^\lab, \widetilde\mvX)=(\mvt^\lab,\mvx)\big)\notag\\
&=\frac{\prod_{j=2}^{\widetilde m}(\tilde d_j-1)!}{|\bT_{\mvs}^{\sss(k)}|}\times\frac{1}{(s_0-2k)!\times\prod_{i\geq 1} s_i!},
\end{align}
where the last step uses \eqref{eqn:997} and \eqref{eqn:998}.
Since this probability is constant and the map from $G$ to $H$ is a bijection, we get the desired result.\qed

\section{Properties of plane trees}\label{sec:proof-plane-tree}
We start by describing the setting and assumptions. 
Assume that for each $m\geq 1$, $\mvs^{\sss(m)}=(s_i^{\sss(m)},~i\ge 0)$ is a tenable \chr{CFD} for a tree on $m$ vertices (thus $\sum_{i\geq 0} s_i^{\sss(m)}=m$). 
\chnr{When there is no scope of confusion, we will simply write $\mvs$ and $s_i$ instead of $\mvs^{\sss (m)}$ and $s_i^{\sss(m)}$.}
Analogous to Assumption \ref{ass:degree}, we make the following assumption on $\{\mvs^{\sss(m)}\}_{m\geq 1}$:

\begin{ass}\label{ass:ecd}
There exists a p.m.f. $(p_0, p_1,\ldots)$ with
\[ p_0>0,\quad  \sum_{i\ge 1}i p_i=1,\quad\text{and}\quad\sum_{i\ge 0}i^2 p_i<\infty\]
such that
\[\frac{s_i}{m}\to p_i\ \text{ for }\ i\ge 0,\ \text{ and }\
\frac{1}{m}\sum_{i\ge 0}i^2 s_i\to\sum_{i\ge 0}i^2 p_i.\]
In particular, $\Delta_m:=\max\set{i: s_i\neq 0}=o(\sqrt{m})$.

We will write $\sigma^2 = \sum_i i^2 p_i - 1$ for the variance associated with the p.m.f. $(p_0, p_1,\ldots)$.
\end{ass}
Then the following was shown in \cite{broutin-marckert}.
\begin{theorem}[{\cite[Theorem 1]{broutin-marckert}}]
\label{thm:broutin-marckert}
Let $\cT_{\mvs}$ be a uniform element of $\bT_{\mvs}$ endowed with the {\bf uniform probability measure on $m$ vertices} and viewed as a metric measure space. Under Assumption \ref{ass:ecd}, as $m\to\infty$,
\[\frac{\sigma}{\sqrt{m}}\cT_{\mvs}\weakc\cT_{2\ve}\]
with respect to the GHP topology (see Definition \ref{def:aldous-crt}).
\end{theorem}

\begin{rem}
In \cite[Theorem 1]{broutin-marckert}, the convergence is stated to hold in the Gromov-Hausdorff sense. However, it is easy to see that the proof in fact implies convergence in the GHP sense.
\end{rem}

The following technical lemma collects all the ingredients necessary for proving our main results.
\chnr{Its proof can be read independently of the rest. So if the reader wishes, they can go through the statement of Lemma \ref{lem:plane-trees} and move on to the next section to read the proofs of the main theorems, and then come back to the proof of this lemma at leisure.}
\begin{lem}\label{lem:plane-trees}
Suppose Assumption \ref{ass:ecd} is satisfied by $\mvs=\mvs^{\sss (m)}$.
Let $\cT_{\mvs}$ be a uniform plane tree with \chr{CFD} $\mvs$. 
Then the following assertions hold.

\vskip5pt

\begin{inparaenumii}
\noindent\item\label{lem:leaves-counting-measure}
For each $k\ge 1$, we can construct independent random vectors $(U_m^{(1)}, V_m^{(1)}),\ldots,(U_m^{(k)}, V_m^{(k)})$ such that $U_m^{(j)}$, $j=1,\ldots,k$, have uniform distribution on the $s_0$ leaves and $V_m^{(j)}$, $j=1,\ldots,k$, have uniform distribution on the $m$ vertices, and
\[m^{-1/2}d_{\cT_{\mvs}}\big(U_m^{(j)}, V_m^{(j)}\big)\weakc 0\ \text{ for }j=1,\ldots,k.\]
In particular,
\[d_{\GHP}\big(\frac{1}{\sqrt{m}}\cT_{\mvs}, \frac{1}{\sqrt{m}}\cT_{\mvs}^{\cL}\big)\weakc 0,\]
where $\cT_{\mvs}^{\cL}$ denotes the metric measure space obtained when the underlying tree is endowed with the uniform probability measure on the set of leaves $\cL(\cT_{\mvs})$. (Recall that the measure on the space $\cT_{\mvs}$ is the uniform probability measure on {\bf all} vertices.)

\medskip

\noindent\item\label{lem:uniform-integrability}
Recall the definition of an admissible pair of leaves (Definition \ref{def:admissible}) and of the function $f_{\cT_{\mvs}}$ (Equation \eqref{eqn:ftu-atu-def}). Let $U_m$ be uniformly distributed over $\cL(\cT_{\mvs})$. Then for every $k\ge 1$,
\[\sup_m\ \E\left(\frac{|\mvA(\cT_{\mvs})|}{s_0\sqrt{m}}\right)^k
\le \sup_m\ \E\left(\frac{f_{\cT_{\mvs}}(U_m)}{\sqrt{m}}\right)^k
<\infty.\]

\medskip

\noindent\item\label{lem:independent-sampling}
For every $k\ge 1$,
\[\frac{1}{m^{3k/2}}\bigg(\big|\mvA(\cT_{\mvs})\big|^k-\big|\mvA_k^{\ord}(\cT_{\mvs})\big|\bigg)\weakc 0.\]

\medskip

\noindent\item\label{lem:finite-dim-convergence}
Let $k\ge 0$ and $\ell\ge 1$.
Suppose $u_1,\ldots,u_k\in\cL(\cT_{\mvs})$, and $v_1,\ldots,v_{\ell}$ are vertices of $\cT_{\mvs}$. 
Write $\mvu=(u_1,\ldots,u_k)$ and $\mvv=(v_1,\ldots,v_{\ell})$.
Let $\cT_{\mvs}(\mvu,\mvv)$ be the element of $\ \vT_{k+\ell}^*$ as in Section \ref{sec:space-of-trees} defined as follows:
If the subtree of $\cT_{\mvs}$ spanned by the root $\rho$ and $u_1,\ldots, u_k, v_1,\ldots,v_\ell$ does not have $(k+\ell)$ distinct leaves, then set $\cT_{\mvs}(\mvu, \mvv)=\partial$.
Otherwise set $\cT_{\mvs}(\mvu, \mvv)$ to be the subtree of $\cT_{\mvs}$ spanned by the root $\rho$ and  $u_1,\ldots, u_k, v_1,\ldots,v_\ell$, 
and for $1\leq i\leq k$,
attach the leaf value $m^{-1/2}f_{\cT_{\mvs}}(u_i)$ to $u_i$, and 
endow $[\rho, u_i]$ with a probability measure by assigning mass $p_x^{\sss(i)}$ to each $x\in[\rho, u_i)$, where
\[p_x^{\sss(i)}:=\frac{1}{f_{\cT_{\mvs}}(u_i)}\cdot
\#\bigg\{v\in\mvA\big(\cT_{\mvs}, u_i\big)\ \ \bigg| \ \
\overleftarrow{\overleftarrow{v}}=x\bigg\}
.\]
(The leaf values and root-to-leaf measures attached to $v_j$, $1\le j\le\ell$, are irrelevant in our proof and can be taken to be zero and $\delta_{\{\rho\}}$ respectively.)

Consider independent random variables $U_m^{(i)}$, $i=1,\ldots, k$, and $V_m^{(j)}$, $j=1,\ldots,\ell$, where $U_m^{(i)}$, $i=1,\ldots,k$ have uniform distribution on the $s_0$ leaves of $\cT_{\mvs}$, and $V_m^{(j)}$, $j=1,\ldots,\ell$ have uniform distribution on the $m$ vertices. 
Let 
$\mvU=(U_m^{(i)}, i=1,\ldots,k)$ and $\mvV=(V_m^{(j)}, j=1,\ldots,\ell)$.
Then
\[\frac{1}{\sqrt{m}}\cT_{\mvs}(\mvU, \mvV)\weakc \frac{1}{\sigma}\cT_{k,\ell},\]
where $\cT_{k,\ell}$ is the random element of $\ \vT_{k+\ell}^*$ constructed as follows: The shape of $\cT_{k,\ell}$ is that of the subtree of $\cT_{2\ve}$ spanned by $(k+\ell)$ points $x_1,\ldots,x_{k+\ell}$ sampled independently according to the mass measure $\mu_{\cT_{2\ve}}$. The leaf weight attached to $x_i$ is $p_0\sigma\cdot\hght(x_i)/2$, $i=1,\ldots,k$, and the measure on $[\rho,x_i]$ is the normalized line measure.

\medskip

\noindent\item\label{lem:conditional-expectation}
The following joint convergence holds:
\begin{align}\label{eqn:nnn}
\bigg(\frac{1}{\sqrt{m}}\cT_{\mvs}, \frac{|\mvA(\cT_{\mvs})|}{s_0\sqrt{m}}\bigg)
\weakc\bigg(\frac{1}{\sigma}\cT_{2\ve},\frac{p_0\sigma}{2}\int_{\cT_{2\ve}}\hght(x)\ \mu_{\cT_{2\ve}}(dx)\bigg)
\end{align}
with respect to product topology induced by GHP topology on the first coordinate and Euclidean topology on the second coordinate.
\end{inparaenumii}
\end{lem}
\subsection{Proof of Lemma \ref{lem:plane-trees}\eqref{lem:leaves-counting-measure}}
Let $\bT_{\mvs}^{\lab}$ be the set of all labeled plane trees on $m$ vertices with the following property:
The vertices are labeled by $[m]$ such that vertices labeled $1,\ldots,s_0$ are leaves, vertices labeled $s_0+1,\ldots, s_0+s_1$ have one child,..., and vertices labeled $(m-s_{{\Delta}_m}+1),\ldots, m$ have $\Delta_m$ many children. As before, we denote by $\xi_j$ the number of children of the vertex labeled $j$.

We will now describe a way of generating a tree uniformly distributed over $\bT_{\mvs}^{\lab}$.
Let $\pi$ be a uniform permutation on $[m]$. Let
\begin{align}\label{eqn:F}
S(j):=\sum_{i=1}^j (\xi_{\pi(i)}-1),\ \ j=1,\ldots,m.
\end{align}
Extend the definition of $\pi$ periodically by letting \chr{$\pi(j)=\pi(j -m)$ for all integers $j$}.
Let $i_0$ denote the location of the first \chr{global} minima of $(S(j),~ 1\leq j\leq m)$ and consider the Vervaat transform w.r.t. this location:
\begin{align}\label{eqn:F-exc}
S^{\exec}(j) = \sum_{i=1}^{j} \big(\xi_{\pi(i_0 + i)}-1\big),\ \ 0\leq j\leq m.
\end{align}
\chr{Let $\cT_{\mvs}^{\lab-}$ be the plane tree whose {\L}ukasiewicz path (see, e.g., \cite{legall-survey} for definition) is $\big(S^{\exec}(j),~0\leq j\leq m\big)$. 
Let $v_1\prec_{\DF} v_2\prec_{\DF} \ldots\prec_{\DF} v_m$ be the vertices of $\cT_{\mvs}^{\lab-}$ arranged in depth-first order (thus $v_1$ is the root of $\cT_{\mvs}^{\lab-}$). 
Label $v_j$ as $\pi(i_0+j)$, $1\leq j\leq m$.
Denote the resulting labeled plane tree by $\cT_{\mvs}^{\lab}$.
Note that $\cT_{\mvs}^{\lab}$ is a random element of $\bT_{\mvs}^{\lab}$, and $\cT_{\mvs}^{\lab-}$ is obtained from $\cT_{\mvs}^{\lab}$ by removing the labels but retaining the plane embedding.}
\begin{lem}\label{lem:generating-tree-from-permutation}
Let $\cT_{\mvs}^{\lab}$ and $\cT_{\mvs}^{\lab-}$ be as above.
\begin{enumeratea}
\item $\cT_{\mvs}^{\lab}\sim\mathrm{Unif}(\bT_{\mvs}^{\lab})$.
\item $\cT_{\mvs}^{\lab-}\equald\cT_{\mvs}$, i.e., $\cT_{\mvs}^{\lab-}\sim\mathrm{Unif}(\bT_{\mvs})$.
\end{enumeratea}
\end{lem}
\noindent{\bf Proof:} Each of the $(m-1)!$ rotation classes of the $m!$ permutations \chr{on $[m]$ gives rise to a unique realization of $\cT_{\mvs}^{\lab}$ and vice versa.
Consequently, $|\bT_{\mvs}^{\lab}|=(m-1)!$, and}
\[\bP(\cT_{\mvs}^{\lab}=\mvt^\lab)=\frac{1}{(m-1)!}=\frac{1}{|\bT_{\mvs}^{\lab}|}\ \text{ for any }\ \mvt^\lab\in\bT_{\mvs}^{\lab}.\]
This implies
\[\bP\big(\cT_{\mvs}^{\lab-}=\mvt\big)=\frac{\prod_{i\ge 0}s_i!}{(m-1)!}=\frac{1}{|\bT_{\mvs}|}\ \text{ for any }\ \mvt\in\bT_{\mvs}.\]
\qed

We now state a useful concentration inequality.
\begin{lem}\label{lem:concentration-uniform-permutation}
There exist universal constants $c_1, c_2>0$ such that for any $m\geq 1$ and probability vector $\mvq:=(q_1,\ldots,q_m)$,
\[\bP\bigg(\max_{j\in[m]}\bigg|\sum_{i=1}^j q_{\pi(i)}-\frac{j}{m}\bigg|\geq x\sigma(\mvq)\bigg)
\leq\exp\big(-c_1 x\log\log x\big),\ \text{ for }\ x\geq c_2,\]
where $\pi$ is a uniform permutation on $[m]$, and $\sigma(\mvq):=\sqrt{q_1^2+\ldots+q_m^2}$.
\chr{Consequently, 
\[\bP\bigg(\max_{1\leq j_1< j_2\leq m}\bigg|\sum_{i=j_1+1}^{j_2} q_{\pi(i)}-\frac{j_2-j_1}{m}\bigg|\geq 2x\sigma(\mvq)\bigg)
\leq 2\exp\big(-c_1 x\log\log x\big),\ \text{ for }\ x\geq c_2.\]
}
\end{lem}

\noindent{\bf Proof:}
The result is essentially contained in \cite[Lemma 4.9]{bhamidi-hofstad-sen}, and we only outline how to extract the result from its proof.
We can work with $\pi$ generated in the following way: let $X_1,\ldots,X_m$ be i.i.d. Unif$[0,1]$, and set $\pi(i)=j_i$, where $X_{j_1}<\ldots<X_{j_m}$. Write $X_{(i)}=X_{j_i}$. Then
\begin{align}\label{eqn:20}
\max_{j\in[m]}\bigg(\sum_{i=1}^j q_{\pi(i)}-\frac{j}{m}\bigg)
\leq\max_{j\in[m]}\bigg(\sum_{i=1}^j q_{\pi(i)}-X_{(j)}\bigg)+\max_{j\in[m]}\bigg|X_{(j)}-\frac{j}{m}\bigg|.
\end{align}
By the DKW inequality \cite{massart-DKW},
\begin{equation}
\label{eqn:dkw}
\pr\bigg(\max_{j\in [m]} \bigg|X_{\sss(j)}-\frac{j}{m}\bigg| \geq \sigma(\mvq) x\bigg)
\leq 2\exp\left(-2m\cdot\left(\sigma(\mvq)x\right)^2\right)
\leq 2\exp\left(-2x^2\right),
\end{equation}
where the last step uses the inequality: $m\sigma(\mvq)^2\geq(\sum_i q_i)^2=1$.
From \cite[Equations (4.14) and (4.15)]{bhamidi-hofstad-sen} we have
\begin{equation}\label{eqn:21}
\bP\bigg(\max_{j\in[m]}\bigg(\sum_{i=1}^j q_{\pi(i)}-X_{(j)}\bigg)\geq x\sigma(\mvq)\bigg)
\leq\exp\big(-c_3 x\log\log x\big),\ \text{ for }\ x\geq c_4.
\end{equation}
Combining \eqref{eqn:20}, \eqref{eqn:dkw}, and \eqref{eqn:21}, it follows that
\begin{equation}\label{eqn:22}
\bP\bigg(\max_{j\in[m]}\bigg(\sum_{i=1}^j q_{\pi(i)}-\frac{j}{m}\bigg)\geq x\sigma(\mvq)\bigg)
\leq\exp\big(-c_5 x\log\log x\big),\ \text{ for }\ x\geq c_6.
\end{equation}
Now note that
\begin{align*}
\max_{j\in[m]}\bigg(\frac{j}{m}-\sum_{i=1}^j q_{\pi(i)}\bigg)
\leq\max_{j\in[m]}\bigg(X_{(j)}-\sum_{i=1}^j q_{\pi(i)}\bigg)+\max_{j\in[m]}\bigg|X_{(j)}-\frac{j}{m}\bigg|.
\end{align*}
Further, the arguments used in the proof of \cite[Equation (4.15)]{bhamidi-hofstad-sen} can be used to prove a tail bound similar to \eqref{eqn:21} for  $\max_{j\in[m]}\big(X_{(j)}-\sum_{i=1}^j q_{\pi(i)}\big)$. Combining this observation with \eqref{eqn:22} yields the desired result.
\qed

\begin{lem}\label{lem:leaves-counting-measure-labeled-tree}
For each $m\geq 1$, let $\mvq=\mvq^{\sss(m)}=(q_1^{\sss(m)},\ldots, q_m^{\sss(m)})$ be a probability vector such that $q_{\max}:=\max_j q_j^{\sss (m)}\to 0$ as $m\to\infty$. Then for each $k\ge 1$, we can construct independent random vectors $(U_m^{(1)}, V_m^{(1)}),\ldots,(U_m^{(k)}, V_m^{(k)})$ on $\cT_{\mvs}^\lab$ such that $U_m^{(j)}$, $j=1,\ldots,k,$ are distributed according \chr{to} $\mvq$ and $V_m^{(j)}$, $j=1,\ldots,k,$ have uniform distribution on the $m$ vertices, and
\[m^{-1/2}d_{\cT_{\mvs}^\lab}\big(U_m^{(j)}, V_m^{(j)}\big)\weakc 0\ \text{ for }j=1,\ldots,k.\]
\end{lem}

\noindent{\bf Proof:}
Recall from around \eqref{eqn:F} \chnr{the construction of $\cT_{\mvs}^\lab$ using a uniform permutation $\pi$.} 
Set
\[ G(x)=\sum_{i=0}^{\lfloor mx\rfloor -1} q_{\pi(i+i_0+1)},\ \ x\in[0,1],\]
where as before $i_0$ is the location of the first global minima of $(S(j),~ 1\leq j\leq m)$.
Let $\sigma(\mvq)$ be as in Lemma \ref{lem:concentration-uniform-permutation}. Then $\sigma(\mvq)\geq q_{\max}\geq 1/m$\chr{, and consequently $y\sigma(\mvq)\geq 2/m$ if $y\geq 2$. Further,
\[
\big|G(x)-x\big|
\leq
\big|G(x)-\frac{\lfloor mx\rfloor}{m}\big|+\frac{1}{m}\ \ \text{ for all }\ \ x\in[0,1].
\]
Combining these observations with Lemma \ref{lem:concentration-uniform-permutation} yields, for $y\geq\max\{2, 4c_2\}$,
\begin{align}\label{eqn:9}
\bP\bigg(\sup_{x\in[0,1]}\bigg|G(x)-x\bigg|\geq y\sigma(\mvq)\bigg)
&\leq\bP\bigg(\sup_{x\in[0,1]}\bigg|G(x)-\frac{\lfloor mx\rfloor}{m}\bigg|\geq y\sigma(\mvq)/2\bigg)
\notag\\
&\leq 2\exp\big(-c_1 (y/4)\log\log(y/4)\big).
\end{align}
}
Since $\sigma(\mvq)\leq\sqrt{q_{\max}}\to 0$ by assumption, we get, in particular,
\begin{align}\label{eqn:10}
\sup_{u\in[0,1]}|G^{-1}(u)-u|\weakc 0,
\end{align}
\chnr{where $G^{-1}(u)=\inf\, \{x\in[0,1] : G(x)\geq u\}$.}

\chnr{Let $H_m$ be the height function of $\cT_{\mvs}^\lab$, i.e.,
\[
H_m(j)=d_{\cT_{\mvs}^\lab}\big(\pi(i_0+1), \pi(i_0+j+1)\big),\ \ j=0,\ldots,m-1.
\]
Then for any $1\leq j_1\leq j_2\leq m$,
\begin{align}\label{eqn:99}
\bigg|d_{\cT_{\mvs}^\lab}\big(\pi(i_0+j_1), \pi(i_0+j_2)\big)-
\bigg(H_m\big(j_1-1\big)+H_m(j_2-1)-2\cdot\min\displaystyle^\star H_m(k)\bigg)\bigg|\le 2,
\end{align}
where $\min^\star$ denotes minimum taken over $j_1-1\leq k\leq j_2-1$.} 
Now \cite[Proposition 5 and Lemma 7]{broutin-marckert} imply that  
\begin{align}\label{eqn:999}
\big(m^{-1/2}H_m(\lfloor mx\rfloor),\ x\in [0,1]\big)\weakc \big(2\ve(x)/\sigma,\ x\in[0,1]\big)
\end{align}
with respect to Skorohod $J_1$ topology on $D[0,1]$. 
Let $X\sim\mathrm{Unif}[0, 1]$. Let $V_m^{\sss(1)}$ (resp. $U_m^{\sss(1)}$) be the 
\chnr{vertex labeled $\pi\big(i_0+\lceil mX\rceil\big)$ (resp. the vertex labeled $\pi\big(i_0+\lceil mG^{-1}(X)\rceil\big)$) in $\cT_{\mvs}^\lab$.} 
Then $V_m^{\sss(1)}$ has uniform distribution on the $m$ vertices, and $U_m^{(1)}$ is distributed according to $\mvq$.
Now, using \eqref{eqn:10}, $|G^{-1}(X)-X|\weakc 0$. 
It thus follows from \eqref{eqn:99} and \eqref{eqn:999} that
\[m^{-1/2} d_{\cT_{\mvs}^\lab}\big(U_m^{(1)}, V_m^{(1)}\big)\weakc 0.\]
We can take $k$ independent copies of $X$ and repeat the same argument to complete the proof.\qed

\medskip

\noindent{\bf Completing the proof of Lemma \ref{lem:plane-trees}\eqref{lem:leaves-counting-measure}:}
By Lemma \ref{lem:generating-tree-from-permutation}, $\cT_{\mvs}\equald\cT_{\mvs}^{\lab-}$. Thus the claim follows by an application of Lemma \ref{lem:leaves-counting-measure-labeled-tree} with the choice
\begin{align*}
\hskip180pt&q_i=\frac{\ind\set{\xi_i=0}}{s_0}.\hskip160pt\blacksquare
\end{align*}

\subsection{Proof of Lemma \ref{lem:plane-trees}\eqref{lem:uniform-integrability}}
If $\mvf$ is a finite forest of plane trees, let the \chr{CFD} of $\mvf$ be
$\tilde\mvs(\mvf)=(\tilde s_i(\mvf), i\geq 0)$,
where $\tilde s_i(\mvf)$ is the number of vertices in $\mvf$ that have exactly $i$ children.
Note that for any sequence of integers $\tilde \mvs=(\tilde s_i, i\ge 0)$ satisfying
\[\tilde s_i\ge 0,\ \ \sum_i i \tilde s_i<\infty,\ \text{ and }\ \sum_i \tilde s_i-\sum_i i \tilde s_i\geq 1,\]
there exists a forest with \chr{CFD} $\tilde \mvs$. Such a forest has exactly 
$(\sum_i \tilde s_i-\sum_i i \tilde s_i)$ 
many trees and $\sum_i \tilde s_i$ many vertices.
Given such a sequence $\tilde \mvs$, let $\bF_{\tilde \mvs}$ denote the set of all plane forests with {\bf ranked roots} having \chr{CFD} $\tilde \mvs$. Thus each forest in $\bF_{\tilde \mvs}$ comes with an ordering of the roots so it makes sense to talk about the ``first'' tree of the forest, the ``second'' tree etc.

The following lemma gives a useful set of estimates.
\begin{lem}\label{lem:moments}
Let $\tilde\mvs$ be a tenable \chr{CFD} for a forest of plane trees, and let
\[\tilde z=\sum_i \tilde s_i-\sum_i i\tilde s_i,\ \ \tilde\Delta=\max_i\set{i:\tilde s_i\neq 0},\ \text{ and }\ \tilde m:=\sum_i\tilde s_i,\]
i.e., $\tilde z$ is the number of trees, $\tilde\Delta$ is the maximum number of children, and $\tilde m$ is the total number of vertices in any forest in $\bF_{\tilde\mvs}$. Sample a forest uniformly from $\bF_{\tilde\mvs}$, and let $X_j$ denote the number of children of the root of the $j$-th tree, $1\le j\le \tilde z$. Then for any $\sigma_1,\ldots,\sigma_r\ge 1$ and $1\le j_1<\ldots<j_r\le\tilde z$,
\begin{align}\label{eqn:23}
\E\bigg[X_{j_1}^{\sigma_1}\times\ldots\times X_{j_r}^{\sigma_r}\bigg]
\le r 2^r\bigg(\sum_i\frac{i^2\tilde s_i}{\tilde m}\bigg)^r
\bigg(1+\frac{\tilde\Delta}{\tilde z}\bigg)
(\tilde\Delta)^{\sigma_1+\ldots+\sigma_r-r},
\end{align}
whenever $r\le\tilde m/2$. 
As a consequence, the sum of the number of children of all roots in the randomly sampled forest satisfies the moment bound
\begin{align}\label{eqn:24}
\E\big[X_1+\ldots+X_{\tilde z}\big]^k\le K_{k}\bigg(1+\frac{\tilde\Delta}{\tilde z}\bigg)
\sum_{r=1}^k\tilde z^r\bigg(\sum_i\frac{i^2\tilde s_i}{\tilde m}\bigg)^r\tilde\Delta^{k-r}
\end{align}
whenever $k\le\tilde m/2$, and where $K_k$ is a constant depending only on $k$.
\end{lem}
\noindent{\bf Proof:} 
Using exchangeability, it is enough to consider $j_1=1,\ldots, j_r=r$ in \eqref{eqn:23}.
Recall (see, e.g., \cite[Equation 6.19]{pitman-book}) that
\begin{align}\label{eqn:pitman}
\big|\bF_{\tilde\mvs}\big|=\frac{\tilde z(\tilde m-1)!}{\prod_{i\ge 0}\tilde s_i!}.
\end{align}

\chnr{Let $\delta_0=(1,0,0\ldots)$, $\delta_1=(0, 1,0,\ldots)$, and similarly define $\delta_i$ for $i\geq 2$.}
\chr{Consider the case where $\tilde z\geq 2$, and either $i_1\neq i_2$ and $\tilde s_{i_1}\geq 1$ and $\tilde s_{i_2}\geq 1$, or $i_1=i_2=i$ and $\tilde s_i\geq 2$ . 
Then from any forest $F$ in $\bF_{\tilde\mvs}$ in which the roots of the first two trees have respectively $i_1$ and $i_2$ many children, we can delete the first two roots and the edges incident to them, declare the children of the first two roots in $F$ to be the roots of the newly created trees, and
then rank the roots of the resulting trees using the planar order to obtain a forest in $\bF_{\tilde\mvs-\delta_{i_1}-\delta_{i_2}}$.
The inverse map is straightforward: Given any forest $F$ in $\bF_{\tilde\mvs-\delta_{i_1}-\delta_{i_2}}$, create two new vertices $\rho_1$ and $\rho_2$, connect each of the first $i_1$ roots of $F$ by an edge to $\rho_1$, and connect each of the next $i_2$ roots of $F$ by an edge to $\rho_2$, and then declare $\rho_1$ and $\rho_2$ to be respectively the first and second roots of the resulting forest, and let the $i_1+i_2+j$-th root of $F$ be the $j+2$-th root of the resulting forest, $j\geq 1$.}
It thus follows that
\begin{align}\label{eqn:25}
\bP\big(X_1=i_1,X_2=i_2\big)=\frac{|\bF_{\tilde\mvs-\delta_{i_1}-\delta_{i_2}}|}{|\bF_{\tilde\mvs}|}.
\end{align}
When $i_1\neq i_2$,
\begin{align*}
\bP\big(X_1=i_1,X_2=i_2\big)
=\frac{(\tilde z-2+i_1+i_2)\tilde s_{i_1}\tilde s_{i_2}}{\tilde z(\tilde m-1)(\tilde m-2)}
\le 4\bigg(1+\frac{i_1+i_2}{\tilde z}\bigg)\frac{\tilde s_{i_1}\tilde s_{i_2}}{\tilde m^2},
\end{align*}
\chr{where the last inequality is true whenever $\tilde m\geq 4$.}
When $i_1=i_2=i$ \chr{and $\tilde m\geq 4$},
\begin{align*}
\bP\big(X_1=i,X_2=i\big)
=\frac{(\tilde z-2+2i)\tilde s_{i}(\tilde s_{i}-1)}{\tilde z(\tilde m-1)(\tilde m-2)}
\le 4\bigg(1+\frac{2i}{\tilde z}\bigg)\frac{\tilde s_{i}^2}{\tilde m^2}.
\end{align*}
In general, \chr{for $\tilde m\geq 2r$},
\begin{align*}
\bP\big(X_1=i_1,\ldots,X_r=i_r\big)
&\le 2^r\bigg(1+\frac{i_1+\ldots+i_r}{\tilde z}\bigg)\frac{\tilde s_{i_1}\times\ldots\times\tilde s_{i_r}}{\tilde m^r}\\
&\le r2^r\bigg(1+\frac{\tilde\Delta}{\tilde z}\bigg)\frac{\tilde s_{i_1}\times\ldots\times\tilde s_{i_r}}{\tilde m^r}.
\end{align*}
Note that
\[\sum_{i\ge 1}i^{\sigma_j}\frac{\tilde s_{i}}{\tilde m}
\le\tilde\Delta^{(\sigma_j-2)\vee 0}\sum_{i\ge 1}i^2\frac{\tilde s_{i}}{\tilde m}
\le\tilde\Delta^{(\sigma_j-1)}\sum_{i\ge 1}i^2\frac{\tilde s_{i}}{\tilde m}.\]
Combining the above, we \chr{see that for $\tilde m\geq 2r$},
\[\E\bigg[X_1^{\sigma_1}\times\ldots\times X_r^{\sigma_r}\bigg]
\le r 2^r\bigg(1+\frac{\tilde\Delta}{\tilde z}\bigg)\bigg(\sum_i\frac{i^2\tilde s_i}{\tilde m}\bigg)^r
\tilde\Delta^{\sum_{j=1}^r(\sigma_j-1)},
\]
which proves \eqref{eqn:23}. The bound in \eqref{eqn:24} follows by a direct expansion.\qed

\medskip

We now introduce some notation.
For any plane tree $\mvt$ and a vertex $u$, define
\[\mvB_1(\mvt,u)=\set{v\ :\ \overleftarrow{v}\in[\rho,u)}\chr{\setminus [\rho, u]},\ \text{ and }\
\mvB_2(\mvt,u)=\set{v\ :\ \chr{\overleftarrow{v}\in\mvB_1(\mvt, u)}},\]
where $\rho$ is the root of $\mvt$. Thus for any $u\in\cL(\mvt)$,
\begin{equation}
\label{eqn:ft-lessb2}
f_{\mvt}(u)=|\mvA(\mvt, u)|\le|\mvB_2(\mvt, u)|.
\end{equation}
For any plane tree $\mvt$ and a vertex $u$, define
\[\mvB_1^-(\mvt,u)=\set{v\ :\ \overleftarrow{v}\in[\rho,u),\ u\prec_{\DF} v},\ \text{ and }\
\mvB_1^+(\mvt,u)=\set{v\ :\ \overleftarrow{v}\in[\rho,u),\ v\prec_{\DF} u}\chr{\setminus [\rho, u)}.\]
So if our convention is to explore the children of a vertex from left to right in a depth-first search, then
$\mvB_1^{-}(\mvt,u)$ (resp. $\mvB_1^{+}(\mvt,u)$) is the collection of vertices that are at distance one from the path $[\rho,u)$ and lie on the right (resp. left) side of $[\rho,u)$.

For a plane tree $\mvt$ and $u\in\cL(\mvt)$, let $\anc^{\sss (1)}(\mvt,u)$ be the plane subtree of $\mvt$ whose vertex set is given by
\[\chr{V=[\rho, u]\cup\mvB_1(\mvt, u)},\]
and further, the vertex $u$ is marked in $\anc^{\sss (1)}(\mvt,u)$.

Now, from \eqref{eqn:At=sum-ftu} and \eqref{eqn:ft-lessb2} it is clear that
\begin{align}\label{eqn:12}
\E\left(\frac{|\mvA(\cT_{\mvs})|}{s_0\sqrt{m}}\right)^k
&=\E\bigg[\E\bigg(\frac{1}{\sqrt{m}}f_{\cT_{\mvs}}(U_m)\bigg|\cT_{\mvs}\bigg)\bigg]^k\notag\\
&\le\E\left(\frac{f_{\cT_{\mvs}}(U_m)^k}{m^{k/2}}\right)
\le\E\left(\frac{|\mvB_2(\cT_{\mvs}, U_m)|^k}{m^{k/2}}\right).
\end{align}

\chr{For any $i\geq 0$, and any plane tree $\mvt$ with root $\rho$ and a marked leaf $u$, write $k_i(\mvt)$ for the number of vertices in $[\rho, u]$ with $i$ many children. 
(Thus, $k_0(\mvt)=1$, as $u$ is the only leaf in $[\rho, u]$.)
Let $\mvk(\mvt)=\big(k_i(\mvt), i\geq 0\big)$.}
Let $\cF$ be the (random) plane forest with ranked roots obtained \chr{from $\cT_{\mvs}$} by deleting the vertices of $[\rho, U_m]$ and the edges incident to them, rooting the resulting trees at the vertices in $\mvB_1(\cT_{\mvs}, U_m)$, and ranking the roots using the depth-first order.
\chr{The CFD of $\cF$ is given by $\big(\mvs-\mvk\big(\anc^{\sss (1)}(\cT_{\mvs}, U_m)\big)\big)$.}
Clearly, there is a bijection between the realizations of $(\cT_{\mvs}, U_m)$ and those of $\big(\anc^{\sss (1)}(\cT_{\mvs}, U_m), \cF\big)$. 
Thus, for all possible realizations $\mvt$ of $\anc^{\sss (1)}(\cT_{\mvs}, U_m)$,
\begin{align}\label{eqn:forest}
\bP\big(\cF=\mvf,\ \anc^{\sss (1)}(\cT_{\mvs}, U_m)=\mvt\big)=\frac{1}{s_0|\bT_{\mvs}|}\ \text{ for all }\ \mvf\in\bF_{\tilde\mvs(\mvt)},
\end{align}
where 
$\tilde\mvs(\mvt):=\mvs-\mvk(\mvt)$, and further,
\[
\bP\big(\anc^{\sss (1)}(\cT_{\mvs}, U_m)=\mvt\big)=\frac{\big|\bF_{\tilde\mvs(\mvt)}\big|}{s_0|\bT_{\mvs}|}.
\]
Writing $\bP_{\mvt}=(\cdot\ |\ \anc^{\sss (1)}(\cT_{\mvs}, U_m)=\mvt)$, it follows that
\begin{equation}\label{eqn:30}
\bP_{\mvt}\big(\cF=\mvf\big)=\frac{1}{|\bF_{\tilde\mvs(\mvt)}|}\ \text{ for any }\ \mvf\in\bF_{\tilde\mvs(\mvt)}.
\end{equation}
Define the conditional expectation operator $\E_{\mvt}$ in an analogous fashion. Let
\[E:=\big\{\hght(U_m)< m/2\big\}.\]
\chr{By \cite[Theorem 1]{addario-berry}, there exists $m_0\geq 1$ such that for all $m\geq m_0$,
\begin{align}\label{eqn:13}
\bP(E^c)\le 7\exp(-\theta m),
\end{align}
where $\theta>0$ is a constant depending only on $p_1$ and $\sum_i i^2p_i$. (Here $(p_j, j\geq 0)$ is as in Assumption \ref{ass:ecd}.)} 
Let $\tilde z, \tilde m, \tilde\Delta$  be as in Lemma \ref{lem:moments} with $\tilde\mvs=\tilde\mvs(\mvt)$. Then $\tilde m\ge m/2$ on the event $E$, $\tilde s_i\le s_i$, and $\tilde\Delta\le\Delta$. 
\chr{In particular, on the event $E$,
\[
\sum_{i\geq 1}\frac{i^2\tilde s_i}{\tilde m}
\leq
2\sum_{i\geq 1}\frac{i^2 s_i}{m}=O(1),
\]
where the last step uses Assumption \ref{ass:ecd}.
}
Thus it follows from Lemma \ref{lem:moments} that
\begin{align*}
\ind_E\E_{\mvt}\bigg[|\mvB_2(\cT_{\mvs}, U_m)|^k\bigg]
\le K\sum_{r=1}^k\bigg(\tilde z^r\Delta^{k-r}+\tilde z^{r-1}\Delta^{k-r+1}\bigg)
\end{align*}
for some constant $K>0$.
Hence,
\begin{align}\label{eqn:14}
\E\bigg[\ind_E|\mvB_2(\cT_{\mvs}, U_m)|^k\bigg]
\le K\sum_{r=1}^k\bigg(\Delta^{k-r}\E\big[\big|\mvB_1(\cT_{\mvs}, U_m)\big|^r\big]
+\Delta^{k-r+1}\E\big[\big|\mvB_1(\cT_{\mvs}, U_m)\big|^{r-1}\big]\bigg).
\end{align}
Let $S^{\exec}$ be as in \eqref{eqn:F-exc}. Letting
\[\mvq:=(q_1,\ldots, q_m),\ \text{ where }\ q_i=\frac{\xi_i}{m-1},\ i=1,\ldots,m,\]
an application of Lemma \ref{lem:concentration-uniform-permutation} shows that
\[\bP\bigg(\max_{j\in[m]}\bigg|\sum_{i=1}^j q_{\pi(i+i_0)}-\frac{j}{m}\bigg|\geq 2x\sigma(\mvq)\bigg)
\leq 2\exp\big(-c_1 x\log\log x\big),\ \text{ for }\ x\geq c_2,\]
where $i_0$ has the same meaning as in \eqref{eqn:F-exc}. Since
\[\bigg|\sum_{i=1}^j \frac{\big(\xi_{\pi(i+i_0)}-1\big)}{m-1}-\bigg(\sum_{i=1}^j q_{\pi(i+i_0)}-\frac{j}{m}\bigg)\bigg|
=\frac{j}{m(m-1)}\leq\frac{1}{m-1},\]
we conclude that
\begin{align*}
\bP\bigg(\max_{j\in[m]}\big|S^{\exec}(j)\big|\ge y\sqrt{m}\bigg)\le\exp\big(-cy\log\log y\big),\ \ y\ge c'.
\end{align*}
This implies that the same tail bound holds for $|\mvB_1^-(\cT_{\mvs}, U_m)|$. 
Since $|\mvB_1^+(\cT_{\mvs}, U_m)|\equald |\mvB_1^-(\cT_{\mvs}, U_m)|$, and $|\mvB_1(\cT_{\mvs}, U_m)|=|\mvB_1^-(\cT_{\mvs}, U_m)|+|\mvB_1^+(\cT_{\mvs}, U_m)|$, we have
\begin{align}\label{eqn:15}
\bP\bigg(\big|\mvB_1(\cT_{\mvs}, U_m)\big|\ge y\sqrt{m}\bigg)\le\exp\big(-cy\log\log y\big),\ \ y\ge c'.
\end{align}
Combining \eqref{eqn:14} and \eqref{eqn:15}, we conclude that
\[\sup_m\ m^{-k/2}\E\bigg[\ind_E|\mvB_2(\cT_{\mvs}, U_m)|^k\bigg]<\infty.\]
From \eqref{eqn:13}, it follows that for $m\geq m_0$,
\[\E\bigg[\ind_E^c\cdot|\mvB_2(\cT_{\mvs}, U_m)|^k\bigg]\le 7m^k\exp(-\theta m).\]
We complete the proof of Lemma \ref{lem:plane-trees}\eqref{lem:uniform-integrability} by combining these observations with \eqref{eqn:12}. \qed

\subsection{Proof of Lemma \ref{lem:plane-trees}\eqref{lem:independent-sampling}}
\chr{It is enough to consider $k\geq 2$, since  $|\mvA(\cT_{\mvs})|=|\mvA_1^{\ord}(\cT_{\mvs})|$.}
Note that Lemma \ref{lem:plane-trees}\eqref{lem:uniform-integrability} shows that $|\mvA(\cT_{\mvs})|/m^{3/2}$ is tight. 
Now,
\begin{align}\label{eqn:16}
&\big|\mvA(\cT_{\mvs})\big|^k-\big|\mvA_k^{\ord}(\cT_{\mvs})\big|\notag\\
&\hskip10pt
=\#\bigg\{\big((u_1,v_1),\ldots,(u_k,v_k)\big)\in\mvA(\cT_{\mvs})^k\ \big|\ \set{u_i,v_i}\cap\set{u_j,v_j}\neq\emptyset\text{ for some } i\neq j\bigg\}\notag\\
&\hskip20pt
\le k^2\big|\mvA(\cT_{\mvs})\big|^{k-2}(R_{11}+2R_{12}+R_{22}),
\end{align}
where
\[R_{11}=\#\bigg\{\big((u_1,v_1),(u_2,v_2)\big)\in\mvA(\cT_{\mvs})^2\ \big|\ u_1=u_2\bigg\},\
R_{12}=\#\bigg\{\big((u_1,v_1),(u_2,v_2)\big)\in\mvA(\cT_{\mvs})^2\ \big|\ u_1=v_2\bigg\},\]
and
\[R_{22}=\#\bigg\{\big((u_1,v_1),(u_2,v_2)\big)\in\mvA(\cT_{\mvs})^2\ \big|\ v_1=v_2\bigg\}.\]
We have,
\begin{align*}
R_{11}
&=\sum_{u\in \cL(\cT_{\mvs})}\ \sum_{v_1\in \mvA(\cT_{\mvs},u)}
\#\bigg\{v_2\ \big|\ (u,v_2)\in\mvA(\cT_{\mvs})\bigg\}\\
&\hskip20pt\le m\sum_{u\in\cL(\cT_{\mvs})}f_{\cT_{\mvs}}(u)
=m\big|\mvA(\cT_{\mvs})\big|
=mO_{P}(m^{3/2})=o_P(m^3).
\end{align*}
A similar argument shows that $R_{12}+R_{22}=o_P(m^3)$. Combined with \eqref{eqn:16}, this gives
\[\big|\mvA(\cT_{\mvs})\big|^k-\big|\mvA_k^{\ord}(\cT_{\mvs})\big|
=O_P\big(m^{3(k-2)/2}\big)o_P(m^3)=o_P(m^{3k/2}),\]
as desired. \qed

\subsection{Proof of Lemma \ref{lem:plane-trees}\eqref{lem:finite-dim-convergence}}
We will make use of the following lemma.
\begin{lem}\label{lem:line-measure}
Let $U_m$ be uniformly distributed on the $s_0$ leaves. Then
\[\frac{1}{\sqrt{m}}\ \max_{k\le\hght(U_m)}
\Bigg|\#\bigg\{v\in\mvA\big(\cT_{\mvs}, U_m\big)\ \bigg| \
\hght\big(\overleftarrow{\overleftarrow{v}}\big)\le k\bigg\}-\frac{p_0\sigma^2 k}{2}\Bigg|\weakc 0.\]
\end{lem}

We first prove Lemma \ref{lem:plane-trees}\eqref{lem:finite-dim-convergence} assuming Lemma \ref{lem:line-measure}.

\medskip

\noindent{\bf Completing the proof of Lemma \ref{lem:plane-trees}\eqref{lem:finite-dim-convergence}:}
Theorem \ref{thm:broutin-marckert} shows that the shape of the subtree of $m^{-1/2}\cT_{\mvs}$ spanned by $(k+\ell)$ vertices sampled independently and uniformly from $[m]$ converges to the shape of $\sigma^{-1}\cT_{k,\ell}$. By Lemma \ref{lem:plane-trees}\eqref{lem:leaves-counting-measure}, the same conclusion holds if the first $k$ vertices are sampled independently and uniformly from the $s_0$ leaves, and the other $\ell$ vertices are sampled independently and uniformly from $[m]$. Convergence of the root-to-leaf measures and the leaf values is a consequence of Lemma \ref{lem:line-measure}.
\qed

The rest of this section is devoted to the proof of Lemma \ref{lem:line-measure}.
We start with the following lemma.
\begin{lem}\label{lem:concentration}
Let $\tilde\mvs=\tilde\mvs_{\kappa}$ be a sequence of \chr{CFD}s indexed by $\kappa$. We will suppress $\kappa$ in the notation most of the time. Let $\tilde m,\ \tilde\Delta,\ \tilde z$ and $X_1,\ldots, X_{\tilde z}$ be as in Lemma \ref{lem:moments}. Let $f:\bZ_{\ge 0}\to\bR_{\ge 0}$. Assume that
\begin{enumeratei}
\item $\tilde z\to\infty$;
\item there exists $a>0$ such that $\sum_{i\ge 0}f(i)\tilde s_i/\tilde m\to a$;
\item $\sup_{\kappa}\sum_{i\ge 0}f^2(i)\tilde s_i/\tilde m<\infty$;
\item $\tilde\Delta=o(\tilde z)$; and
\item $\max_{1\le i\le\tilde\Delta}f(i)=o(\tilde z)$.
\end{enumeratei}
Then
\[\frac{1}{\tilde z}\cdot\max_{1\le j\le\tilde z}\big|\sum_{i=1}^j f(X_i)-aj\big|\weakc 0.
\]
\end{lem}

\noindent{\bf Proof:}
An argument similar to the one used in \eqref{eqn:25} gives
\[
\bP\big(X_1=i\big)
=\frac{(\tilde z-1+i)\tilde s_{i}}{\tilde z(\tilde m-1)}.
\]
Hence
\begin{align}\label{eqn:26}
\E\bigg[\sum_{i\ge 1}^{\tilde z}\frac{f(X_i)}{\tilde z}\bigg]
=\E\big[f(X_1)\big]=\sum_{i=1}^{\tilde\Delta}
\bigg(\frac{(\tilde z-1+i)\tilde m}{\tilde z(\tilde m-1)}\bigg)\frac{f(i)\tilde s_{i}}{\tilde m}\to a.
\end{align}
Similarly, using \eqref{eqn:25}, a direct computation shows that
$\cov\big(f(X_1), f(X_2)\big)\to 0,$
which in turn implies
\begin{align}\label{eqn:27}
\var\bigg[\sum_{i\ge 1}^{\tilde z}\frac{f(X_i)}{\tilde z}\bigg]\to 0.
\end{align}
Combining \eqref{eqn:26} and \eqref{eqn:27}, we see that
\begin{align}\label{eqn:29}
\sum_{i\ge 1}^{\tilde z}f(X_i)/\tilde z\probc a.
\end{align}

Let $\hat S=(\hat S_0,\hat S_1,\ldots)$ denote the \chnr{frequency} distribution of $X_1,\ldots, X_{\tilde z}$.
Since
\begin{align*}
\bP\big(X_1=i_1,\ldots, X_{\tilde z}=i_{\tilde z}\big)=\frac{|\bF_{\mvs-\hat\mvs}|}{|\bF_{\mvs}|}
\end{align*}
for any $(i_1,\ldots,i_{\tilde z})$ with \chnr{frequency} distribution $\hat\mvs$,
\begin{align}\label{eqn:17}
\bP\big(X_1=i_1,\ldots, X_{\tilde z}=i_{\tilde z}\ \big|\ \hat S=\hat\mvs\big)=\frac{\prod_{i\geq 0}\hat s_i!}{\tilde z!}
\end{align}
Define $y_1,\ldots, y_{\tilde z}$ as follows:
\begin{align*}
y_1=\ldots=y_{\hat S_0}=0,\ \ y_{\hat S_0+1}=\ldots=y_{\hat S_0+\hat S_1}=1,\ldots.
\end{align*}
Then conditional on $\hat S$, the distribution \eqref{eqn:17} can be generated by uniformly permuting $y_1,\ldots,y_{\tilde z}$ and removing the $y$ labels. Set
\[\hat\mvq:=(\hat q_1,\ldots,\hat q_{\tilde z}),\ \text{ where }\
\hat q_i=\frac{f(y_i)}{\sum_{j=1}^{\tilde z}f(y_j)}.\]
From Lemma \ref{lem:concentration-uniform-permutation}, for a uniform permutation $\pi$ (independent of $\hat S$) on $\tilde z$ elements and $\eps>0$,
\begin{align}\label{eqn:28}
\bP\bigg(\max_{1\le j\le\tilde z}\bigg|\sum_{i=1}^j \hat q_{\pi(i)}-\frac{j}{\tilde z}\bigg|\ge\eps\ \bigg|\ \hat S\bigg)
\le \exp\bigg(-c_1\bigg(\frac{\eps}{\sigma(\hat\mvq)}\bigg)\log\log\bigg(\frac{\eps}{\sigma(\hat\mvq)}\bigg)\bigg)\
\text{ on }\ \big\{c_2\sigma(\hat\mvq)\leq\eps\big\}.
\end{align}
Since
\[\sigma(\hat\mvq)^2\le \hat q_{\max}
=\frac{\tilde z}{\sum_{1\le i\le\tilde z}f(X_i)}\times\max_{1\le i\le\tilde\Delta}\frac{f(i)}{\tilde z}
\probc 0,\]
we conclude from \eqref{eqn:28} that
\[
\bP\bigg(\max_{1\le j\le\tilde z}\bigg|
\frac{\sum_{i=1}^j f(X_i)}{\sum_{k=1}^{\tilde z}f(X_k)}-\frac{j}{\tilde z}\bigg|\ge\eps\bigg)\to 0,
\]
which combined with \eqref{eqn:29} yields the claim.
\qed

\medskip


\noindent{\bf Proof of Lemma \ref{lem:line-measure}:}
It follows from \cite[Proposition 5]{broutin-marckert} that
\begin{align}\label{eqn:31}
\frac{1}{\sqrt{m}}\ \max_{k\le\hght(U_m)}\bigg|\#\bigg\{v\ \big|\
 v\in\mvB_1^-(\cT_{\mvs}, U_m),\ \hght\big(\overleftarrow{v}\big)\le k\bigg\}-\frac{\sigma^2 k}{2}\bigg|\probc 0,
\end{align}
and a similar statement is true with $\mvB_1^+$ by symmetry.
Now recall from \eqref{eqn:30} that conditional on $\anc^{\sss (1)}(\cT_{\mvs},U_m)$, the forest $\cF$ (defined around \eqref{eqn:forest}) is uniformly distributed over $\bF_{\tilde\mvs}$, where $\tilde\mvs$ is the \chr{CFD} of the vertices in $[\rho, U_m]^c$.
We make the following observations about the forest $\cF$:
\begin{enumeratei}
\item The \chr{CFD} $\tilde\mvs$ of $\cF$ satisfies
\[
s_i-\hght(U_m)\leq\tilde s_i\leq s_i,
\]
where $\hght(U_m)=\Theta_P(\sqrt{m})$ by Theorem \ref{thm:broutin-marckert}.
\item Similarly, the number of vertices $\tilde m$ (say) in $\cF$ satisfies
\[
m-\hght(U_m)\leq\tilde m\leq m.
\]
\item \eqref{eqn:31} and its analogue for $\mvB_1^+$ combined with the fact $\hght(U_m)=\Theta_P(\sqrt{m})$ shows that the number of roots of $\cF$, namely $|\mvB_1(\cT_{\mvs}, U_m)|$ satisfies
\[|
\mvB_1(\cT_{\mvs}, U_m)|=\Theta_P(\sqrt{m}).
\]
\end{enumeratei}
Thus, when $\mvs$ satisfies Assumption \ref{ass:ecd}, the \chr{CFD} of $\cF$ and the function $f(i)=i$ satisfy the assumptions of Lemma \ref{lem:concentration} with $a=1$. This combined with \eqref{eqn:31} gives
\begin{align}\label{eqn:32}
\frac{1}{\sqrt{m}}\ \max_{k\le\hght(U_m)}\bigg|\#\bigg\{v\ \big|\
 \overleftarrow{v}\in\mvB_1^-(\cT_{\mvs}, U_m),\ \hght\big(\overleftarrow{\overleftarrow{v}}\big)\le k\bigg\}
 -\frac{\sigma^2 k}{2}\bigg|\probc 0,
\end{align}
and a similar statement is true with $\mvB_1^+$.

Let $\anc^{\sss(2)}(\mvt, u)$ be the plane subtree of $\mvt$ whose vertex set is given by
\[\chr{V^{\sss(2)}=[\rho, u]\cup\mvB_1(\mvt, u)\cup\mvB_2(\mvt, u)},\]
and further, the vertex $u$ is marked in $\anc^{\sss (2)}(\mvt,u)$.
Let $\cF^{\sss(2)}$ be the plane forest with ranked roots obtained by deleting the vertices of $\anc^{\sss (1)}(\cT_{\mvs},U_m)$ and the edges incident to them, rooting the resulting trees at the vertices of $\mvB_2(\cT_{\mvs}, U_m)$, and ranking them in the depth-first order.
Then conditional on $\anc^{\sss(2)}(\cT_{\mvs},U_m)$, $\cF^{\sss(2)}$ is again uniformly distributed over the set of plane forests with ranked roots with the remaining child sequence.
Further, reasoning similar to above shows that the \chr{CFD} of $\cF^{\sss(2)}$ and the function $f(i)=\ind\set{i=0}$ satisfies the assumptions of Lemma \ref{lem:concentration} with $a=p_0$.

Applying Lemma \ref{lem:concentration} to the forest $\cF^{\sss(2)}$ and the function $f(i)=\ind\set{i=0}$, and combining this with \eqref{eqn:32} yields the claim in Lemma \ref{lem:line-measure}.
\qed

\subsection{Proof of Lemma \ref{lem:plane-trees}\eqref{lem:conditional-expectation}}
First recall from Lemma \ref{lem:plane-trees} \eqref{lem:leaves-counting-measure} that $\cT_{\mvs}^{\cL}$ denoted the metric measure space obtained when the underlying tree is endowed with the uniform probability measure on $\cL(\cT_{\mvs})$.  Let $U_m^{\sss(i)}$, $1\leq i\leq k$, and $x_1,\ldots,x_k$ be as in Lemma \ref{lem:plane-trees}\eqref{lem:finite-dim-convergence}.
Lemma \ref{lem:plane-trees}\eqref{lem:leaves-counting-measure} together with Theorem \ref{thm:broutin-marckert} shows that for all $k\geq 1$,
\begin{align*}
&\left(\frac{1}{\sqrt{m}}\cT_{\mvs},\
\frac{1}{\sqrt{m}}\cT_{\mvs}^{\cL},\
\frac{1}{k\sqrt{m}}\bigg(\hght(U_m^{\sss(1)})+\ldots+\hght(U_m^{\sss(k)})\bigg)\right)\\
&\hskip60pt
\weakc\left(\frac{1}{\sigma}\cT_{2\ve},\
\frac{1}{\sigma}\cT_{2\ve},\
\frac{1}{k\sigma}\bigg(\hght(x_1)+\ldots+\hght(x_k)\bigg)
\right)
\end{align*}
with respect to product topology induced by GHP topology on the first two coordinates and Euclidean topology on $\bR$ on the third coordinate.
Thus by Lemma \ref{lem:line-measure},
\begin{align}\label{eqn:33}
\left(\frac{1}{\sqrt{m}}\cT_{\mvs},\
\frac{1}{k\sqrt{m}}\bigg(f_{\cT_{\mvs}}(U_m^{\sss(1)})+\ldots+f_{\cT_{\mvs}}(U_m^{\sss(k)})\bigg)\right)
\weakc\left(\frac{1}{\sigma}\cT_{2\ve},\
\frac{p_0\sigma}{2k}\bigg(\hght(x_1)+\ldots+\hght(x_k)\bigg)
\right).
\end{align}
Now for any $\eps>0$ and $k\geq 1$,
\begin{align*}
&\bP\bigg(\bigg|
\frac{1}{k\sqrt{m}}\bigg(f_{\cT_{\mvs}}(U_m^{\sss(1)})+\ldots+f_{\cT_{\mvs}}(U_m^{\sss(k)})\bigg)-\frac{|\mvA(\cT_{\mvs})|}{s_0\sqrt{m}}
\bigg|\geq\eps
\bigg)\\
&\hskip20pt
=\bP\bigg(\bigg|
\frac{1}{k\sqrt{m}}\bigg(f_{\cT_{\mvs}}(U_m^{\sss(1)})+\ldots+f_{\cT_{\mvs}}(U_m^{\sss(k)})\bigg)-
\frac{\E\big(f_{\cT_{\mvs}}(U_m^{\sss(1)})\big|\cT_{\mvs}\big)}{\sqrt{m}}
\bigg|\geq\eps
\bigg)\\
&\hskip40pt
\leq\frac{1}{\eps^2 k m}\E\bigg[\var\bigg(f_{\cT_{\mvs}}(U_m^{\sss(1)})\big|\cT_{\mvs}\bigg)\bigg]
\leq\frac{1}{\eps^2 k m}\E\bigg[f_{\cT_{\mvs}}(U_m^{\sss(1)})^2\bigg]\leq\frac{C}{\eps^2k},
\end{align*}
where the first equality holds because of \eqref{eqn:At=sum-ftu} and the last step uses Lemma \ref{lem:plane-trees}\eqref{lem:uniform-integrability}.
By a similar argument, we can show that
\begin{align*}
\bP\bigg(\bigg|
\frac{1}{k}\big(\hght(x_1)+\ldots+\hght(x_k)\big)-
\int_{\cT_{2\ve}}\hght(x)\mu_{\cT_{2\ve}}(dx)\bigg|
\geq\eps
\bigg)\leq\frac{C}{\eps^2k}.
\end{align*}
These observations combined with \eqref{eqn:33} yield
\[\bigg(\frac{1}{\sqrt{m}}\cT_{\mvs}, \frac{|\mvA(\cT_{\mvs})|}{s_0\sqrt{m}}\bigg)
\weakc\bigg(\frac{1}{\sigma}\cT_{2\ve},\frac{p_0\sigma}{2}\int_{\cT_{2\ve}}\hght(x)\ \mu_{\cT_{2\ve}}(dx)\bigg),\]
which is the desired result. \qed

\section{Asymptotics for connected graphs with given degree sequence}
\label{sec:proof-conn}
\chr{The aim of this section is to prove Theorems \ref{prop:condition-on-connectivity} and \ref{thm:number-of-connected-graphs}.
Recall Theorem \ref{lem:alternate-construction} that described an algorithm for generating the uniform measure on the space of simple connected graphs with a prescribed degree sequence with some fixed number $k$ of surplus edges. Using this Theorem together with the technical \chr{Lemma} \ref{lem:plane-trees}, we will complete the proofs of the above two theorems.}

\subsection{Proof of Theorem \ref{prop:condition-on-connectivity}}
It suffices to work with $k\geq 1$.
Let $\mvs=\mvs^{\sss(m)}$ be as in \eqref{eqn:35}.
Note that when $\{\widetilde\vd^{\sss(\widetilde m)}\}_{\widetilde m\geq 1}$ satisfies Assumption \ref{ass:degree}, $\{\mvs^{\sss(m)}\}_{m\geq 1}$ satisfies Assumption \ref{ass:ecd} with limiting p.m.f. $(p_i, i\geq 0)$, where
\[p_i:=\widetilde p_{i+1},\ \ i=0,1,\ldots.\]
In view of Theorem \ref{lem:alternate-construction}, it is enough to prove the result for $\widetilde m^{-1/2}\cI(\widetilde\cT_{\mvs}, \widetilde\mvX)$.

\begin{figure}[htbp]
	\centering
	\pgfmathsetmacro{\nodebasesize}{1} 
	
	\begin{tikzpicture}[scale=.25,  iron/.style={circle, minimum size=6mm, inner sep=0pt, ball color=red!20}, wat/.style= {circle, inner sep=3pt, ball color=green!20}
	]

	\node (1) [iron] at (0,0) {$\rho$};
	\node (11) [iron] at (-4,5) {$\overleftarrow{\overleftarrow{v}}$};
	\node (12) [iron] at (4,5) {};
	\node (111) [iron] at (-8,10) {$\overleftarrow{u}$};
	\node (112) [iron] at (1,10) {$\overleftarrow{v}$};
	\node (1121) [wat] at (5,15) {${v}$};

	\draw[blue, very thick] (1)-- (11) -- (111);
	\draw[blue, very thick] (1) -- (12);
	\draw[blue, very thick] (11) -- (112) -- (1121);
	\draw[red, very thick]   (11) to[out=240,in=190] (111);

	\end{tikzpicture}
	\caption{An example of the operation $\cQ$ applied on the tree $\mvt$ and admissible pair $(u,v)$ in Figure \ref{fig:admiss}. }
	\label{fig:Q-t-x}
\end{figure}

For a plane tree $\mvt$ and $\mvx=\big((u_1,v_1),\ldots,(u_k,v_k)\big)\in\mvA(\mvt)^k$, let $\cQ(\mvt, \mvx)$ be the space obtained by identifying $u_j$ and $\overleftarrow{\overleftarrow{v}}_j$ \chnr{(or equivalently, deleting $u_j$ and adding an edge between $\overleftarrow{u_j}$ and $\overleftarrow{\overleftarrow{v}}_j$)} for $1\leq j\leq k$; see Figure \ref{fig:Q-t-x} for an illustration.
We endow the space $\cQ(\mvt, \mvx)$ with the graph distance and the push-forward of the uniform probability measure on $\mvt$.

\chr{Identifying the vertices of $\cI(\widetilde\cT_{\mvs}, \widetilde\mvX)$ (resp. $\cQ(\widetilde\cT_{\mvs}, \widetilde\mvX)$) with vertices of $\widetilde\cT_{\mvs}$, for every $z\in\cQ(\widetilde\cT_{\mvs}, \widetilde\mvX)\setminus\{v_j : 1\leq j\leq k\}$, we can find a corresponding vertex in $\cI(\widetilde\cT_{\mvs}, \widetilde\mvX)$, which we also denote by $z$.
Let $\cR\subseteq\cQ(\widetilde\cT_{\mvs}, \widetilde\mvX)\times\cI(\widetilde\cT_{\mvs}, \widetilde\mvX)$ be the correspondence
\[
\cR=\big\{(z,z) : z\notin\{v_1,\ldots,v_k\}\big\}\cup\big\{(v_j,\overleftarrow{v_j}) : 1\leq j\leq k\big\}.
\]
Consider a geodesic path $\cP=(z=w_1,\ldots,w_r=z')$ in $\cQ(\widetilde\cT_{\mvs}, \widetilde\mvX)$ between $z,z'\in\cQ(\widetilde\cT_{\mvs}, \widetilde\mvX)\setminus\{v_1,\ldots,v_k\}$.
By replacing every consecutive occurrence of $\overleftarrow{u_j}, \overleftarrow{\overleftarrow{v_j}}$ (resp. $\overleftarrow{\overleftarrow{v_j}}, \overleftarrow{u_j}$) in $\cP$ by $\overleftarrow{u_j}, \overleftarrow{v_j}, \overleftarrow{\overleftarrow{v_j}}$ (resp. $\overleftarrow{\overleftarrow{v_j}}, \overleftarrow{v_j}, \overleftarrow{u_j}$), we get a path between $z$ and $z'$ in $\cI(\widetilde\cT_{\mvs}, \widetilde\mvX)$.
From this it follows that
\[
d_{\cI}(w,w')\leq d_{\cQ}(z,z')+k\ \ \text{ for all }\ \ (z,w), (z',w')\in \cR,
\]
where $d_\cI$ and $d_\cQ$ respectively denote the metrics on $\cI(\widetilde\cT_{\mvs}, \widetilde\mvX)$ and $\cQ(\widetilde\cT_{\mvs}, \widetilde\mvX)$.
By a similar argument we can show that 
\[
d_{\cQ}(z,z')\leq d_{\cI}(w,w')+k+2\ \ \text{ for all }\ \ (z,w), (z',w')\in \cR.
\]
Hence
\begin{align}\label{eqn:hhh}
\dis(\cR)\leq k+2.
\end{align}
Let $\pi$ be the probability measure on $\cQ(\widetilde\cT_{\mvs}, \widetilde\mvX)\times\cI(\widetilde\cT_{\mvs}, \widetilde\mvX)$ given by 
\[
\pi(\{z,z\})=\frac{1}{\widetilde m -1}\ \ \text{ for }\ \ z\in\cQ(\widetilde\cT_{\mvs}, \widetilde\mvX)\setminus\{v_1,\ldots,v_k\}.
\]
Clearly, $\pi(\cR^c)=0$.
Note also that the projection of $\pi$ onto $\cI(\widetilde\cT_{\mvs}, \widetilde\mvX)$ is the uniform probability measure on the vertices of $\cI(\widetilde\cT_{\mvs}, \widetilde\mvX)$.
Now, in $\cQ(\widetilde\cT_{\mvs}, \widetilde\mvX)$, the measure assigned to $\overleftarrow{\overleftarrow{v_j}}$ is at most $k/(\widetilde m-1+2k)$, $1\leq j\leq k$. 
Thus the discrepancy of $\pi$ is bounded by
\[
\frac{1}{2}
\Big[\Big(\frac{1}{\widetilde m-1}-\frac{1}{\widetilde m-1+2k}\Big)\times(\widetilde m-1)
+\frac{k^2}{\widetilde m-1+2k}
+\frac{k}{\widetilde m-1+2k}\Big]
\leq 
\frac{3k+k^2}{2\widetilde m}.
\]
Combining the last display with \eqref{eqn:hhh}, we get
\[d_{\GHP}\bigg(\frac{1}{\sqrt{\widetilde m}}\cI(\widetilde\cT_{\mvs},
\widetilde\mvX),\frac{1}{\sqrt{\widetilde m}}\cQ(\widetilde\cT_{\mvs}, \widetilde\mvX)\bigg)\leq \frac{k+2}{2\sqrt{\widetilde m}}+\frac{3k+k^2}{2\widetilde m}.
\]
}Thus it is enough to prove the result for the space $m^{-1/2}\cQ(\widetilde\cT_{\mvs}, \widetilde\mvX)$, where $m=\widetilde m-1+2k$ as defined in \eqref{eqn:def-m}.

Recall the Gromov-weak topology from Section \ref{sec:gromov-weak}. We will first prove convergence of $m^{-1/2}\cQ(\widetilde\cT_{\mvs}, \widetilde\mvX)$ in the Gromov-weak topology by making use of a technique from \cite{bhamidi-hofstad-sen} and then strengthen it to convergence in the GHP sense. Let $\Phi$ and $\phi$ be as in \eqref{eqn:polynomial-func-def}. Then
\begin{align}\label{eqn:1}
\E\left(\Phi\bigg(\frac{1}{\sqrt{m}}\cQ\big(\widetilde\cT_{\mvs}, \widetilde\mvX\big)\bigg)\right)
&=\frac{\sum_{(\mvt,\mvx)\in\bT_{\mvs}^{\sss(k)}}\Phi\bigg(\frac{1}{\sqrt{m}}\cQ\big(\mvt,\mvx\big)\bigg)}{|\bT_{\mvs}^{\sss(k)}|}\notag\\
&=\frac{\sum_{\mvt\in\bT_{\mvs}}\sum_{\mvx\in\mvA_k(\mvt)}
\Phi\bigg(\frac{1}{\sqrt{m}}\cQ\big(\mvt,\mvx\big)\bigg)\big/\big(|\bT_{\mvs}|\cdot s_0^k m^{k/2}\big)}{\sum_{\mvt\in\bT_{\mvs}}|\mvA_k(\mvt)|\big/\big(|\bT_{\mvs}|\cdot s_0^k m^{k/2}\big)}\notag\\
&=\frac{\sum_{\mvt\in\bT_{\mvs}}\sum_{\mvx\in\mvA_k^{\ord}(\mvt)}
\Phi\bigg(\frac{1}{\sqrt{m}}\cQ\big(\mvt,\mvx\big)\bigg)\big/\big(|\bT_{\mvs}|\cdot s_0^k m^{k/2}\big)}{\sum_{\mvt\in\bT_{\mvs}}|\mvA_k^{\ord}(\mvt)|\big/\big(|\bT_{\mvs}|\cdot s_0^k m^{k/2}\big)}.
\end{align}

\chr{Recall that $\cT_{\mvs}$ is a uniform plane tree with CFD $\mvs$.}
Then
\begin{align}\label{eqn:4}
\frac{\sum_{\mvt\in\bT_{\mvs}}\sum_{\mvx\in\mvA_k^{\ord}(\mvt)}
\Phi\bigg(\frac{1}{\sqrt{m}}\cQ\big(\mvt,\mvx\big)\bigg)}{|\bT_{\mvs}|\cdot s_0^k m^{k/2}}
&=\E\bigg[\sum_{\mvx\in\mvA_k^{\ord}(\cT_{\mvs})}
\Phi\bigg(\frac{1}{\sqrt{m}}\cQ\big(\cT_{\mvs},\mvx\big)\bigg)\frac{1}{s_0^km^{k/2}}\bigg]\\
&=\E\bigg[\sum_{\mvx\in\mvA(\cT_{\mvs})^k}
\Phi\bigg(\frac{1}{\sqrt{m}}\cQ\big(\cT_{\mvs},\mvx\big)\bigg)\frac{1}{s_0^km^{k/2}}\bigg]+o(1),\notag
\end{align}
where the second equality follows from Lemma \ref{lem:plane-trees} \eqref{lem:uniform-integrability} and \eqref{lem:independent-sampling}. Writing
\[\mvx=\big((u_1,y_1),\ldots,(u_k,y_k)\big), \
\sum\displaystyle_1=\sum_{\substack{u_1\in\cL(\cT_{\mvs})\\ \vdots\\ u_k\in\cL(\cT_{\mvs})}},
\ \text{ and }\
\sum\displaystyle_2=\sum_{\substack{y_1\in\mvA(\cT_{\mvs},u_1)\\ \vdots\\ y_k\in\mvA(\cT_{\mvs},u_k)}},
\]
we note that
\begin{align}\label{eqn:5}
\frac{1}{s_0^km^{k/2}}\sum_{\mvx\in\mvA(\cT_{\mvs})^k}
\Phi\bigg(\frac{1}{\sqrt{m}}\cQ\big(\cT_{\mvs},\mvx\big)\bigg)
&=\frac{1}{s_0^km^{k/2}}
\sum\displaystyle_1\sum\displaystyle_2\ \Phi\bigg(\frac{1}{\sqrt{m}}\cQ\big(\cT_{\mvs},\mvx\big)\bigg)\\
&=
\sum\displaystyle_1\bigg(\frac{1}{s_0^k}\bigg)\prod_{i=1}^k \bigg(\frac{f_{\cT_{\mvs}}(u_i)}{\sqrt{m}}\bigg)
\sum\displaystyle_2  \prod_{i=1}^k \bigg(\frac{1}{f_{\cT_{\mvs}}(u_i)}\bigg)\Phi\bigg(\frac{1}{\sqrt{m}}\cQ\big(\cT_{\mvs},\mvx\big)\bigg).\notag
\end{align}

We now recall some constructs from \cite{bhamidi-hofstad-sen}.
Recall the space $\vT_{J}^*$ from Section \ref{sec:space-of-trees}. 
Let $\mvt$ be an element in $\vT_{k+\ell}^*$
\chr{with root $0+$ and leaves $1+,\ldots, (k+\ell)+$}.
Recall \chr{from Section \ref{sec:space-of-trees}} that for each $i$, there is a probability measure $\nu_{\mvt,i }(\cdot)$ on the path  $[0+, i+]$ for $1\leq i\leq k+\ell$. For $1\leq i\leq k$, sample $y_i$ according to the distribution  $\nu_{\mvt,i}(\cdot)$ independently for different $i$ and identify $i+$ with $y_i$. Let $\mvt'$ denote the (random) space thus obtained, and let $d_{\mvt'}$ denote the induced metric on $\mvt'$. Define the function $g^{\sss(k)}_\phi:\vT_{k+\ell}^*\to \bR$ by
\begin{align}
\label{eqn:gphi-def}
g^{\sss(k)}_{\phi}(\mvt):=
\left\{
\begin{array}{l}
\E\left[\phi\left(d_{\mvt'}(i+, j+): k+1\leq i<j\leq k+\ell\right)\right],\text{ if }\mvt\neq\partial,\\
0,\text{ if }\mvt=\partial.
\end{array}
\right.
\end{align}
In words, we look at the expectation of $\phi$ applied to the pairwise distances between the last $\ell$ leaves after sampling $y_i$ on the path $[0+, i+]$ for $1\leq i\leq k$ and identifying $i+$ with $y_i$. Note that here the expectation is only taken over the choices of $y_i$.

Write $d_{\cQ}$ for the induced metric on the space $m^{-1/2}\cQ\big(\cT_{\mvs},\mvx\big)$,
and let $\sum_3=\sum_{v_1,\ldots,v_{\ell}\in [m]}$. Then
\[\Phi\bigg(\frac{1}{\sqrt{m}}\cQ\big(\cT_{\mvs},\mvx\big)\bigg)
=\sum\displaystyle_3\frac{1}{m^{\ell}}\phi\bigg(d_{\cQ}(v_i, v_j): 1\leq i<j\leq\ell\bigg).\]
Write $\mvu=(u_1,\ldots, u_k)$, $\mvv=(v_1,\ldots, v_{\ell})$,
and let $\mvU=(U_m^{\sss(i)}, 1\leq i\leq k),\ \mvV=(V_m^{\sss(j)}, 1\leq j\leq \ell)$, and $\cT_{\mvs}(\mvu,\mvv)$ be as in Lemma \ref{lem:plane-trees}\eqref{lem:finite-dim-convergence}. 
Then we immediately see that
\[
\bigg|\sum\displaystyle_2\prod_{i=1}^k
\bigg(\frac{1}{f_{\cT_{\mvs}}(u_i)}\bigg)\Phi\bigg(\frac{1}{\sqrt{m}}\cQ\big(\cT_{\mvs},\mvx\big)\bigg)
-\sum\displaystyle_3\frac{1}{m^{\ell}}
g_{\phi}^{\sss(k)}\bigg(\frac{1}{\sqrt{m}}\cT_{\mvs}(\mvu,\mvv)\bigg)\bigg|
\leq \|\phi\|_{\infty}\pr\big(\cT_{\mvs}(\mvu,\mvV)=\partial\big|\ \cT_{\mvs}\big).
\]
\chr{Now $\pr\big(\cT_{\mvs}(\mvU,\mvV)=\partial\big)\to 0$ as a consequence of Lemma \ref{lem:plane-trees}\eqref{lem:finite-dim-convergence}.
Further, by Lemma \ref{lem:plane-trees}\eqref{lem:uniform-integrability}, 
$\big\{\prod_{i=1}^k \big(f_{\cT_{\mvs}}(U_m^{\sss(i)})/\sqrt{m}\big)\big\}_{m\geq 1}$ is uniformly integrable.
Thus, combining the last display with \eqref{eqn:4} and \eqref{eqn:5} yields}
\begin{align}\label{eqn:6}
&\frac{\sum_{\mvt\in\bT_{\mvs}}\sum_{\mvx\in\mvA_k^{\ord}(\mvt)}
\Phi\bigg(\frac{1}{\sqrt{m}}\cQ\big(\mvt,\mvx\big)\bigg)}{|\bT_{\mvs}|\cdot s_0^k m^{k/2}}
=\E\bigg[\prod_{i=1}^k \bigg(\frac{f_{\cT_{\mvs}}(U_m^{\sss(i)})}{\sqrt{m}}\bigg)
g_{\phi}^{\sss(k)}\bigg(\frac{1}{\sqrt{m}}\cT_{\mvs}(\mvU,\mvV)\bigg)\bigg]+o(1).
\end{align}
Since the functional $g_{\phi}^{\sss(k)}$ is continuous on the space $\vT_{k+\ell}^*$ \cite[Proposition 4.25]{bhamidi-hofstad-sen}, combining \eqref{eqn:6} with Lemma \ref{lem:plane-trees}\eqref{lem:finite-dim-convergence} and using uniform integrability
(Lemma \ref{lem:plane-trees}\eqref{lem:uniform-integrability}),
we get
\begin{align}\label{eqn:7}
&\frac{\sum_{\mvt\in\bT_{\mvs}}\sum_{\mvx\in\mvA_k^{\ord}(\mvt)}
\Phi\bigg(\frac{1}{\sqrt{m}}\cQ\big(\mvt,\mvx\big)\bigg)}{|\bT_{\mvs}|\cdot s_0^k m^{k/2}}
\\
&\hskip40pt\to\left(\frac{p_0\sigma}{2}\right)^k
\E\bigg[\int_{x_1\in\cT_{2\ve}}\ldots\int_{x_{k+\ell}\in\cT_{2\ve}}\mu_{\cT_{2\ve}}^{\otimes k+\ell}(dx_1\ldots dx_{k+\ell})\prod_{i=1}^k\hght(x_i)\cdot
g_{\phi}^{\sss(k)}\bigg(\frac{1}{\sigma}\cT_{k,\ell}\bigg)\bigg].\notag
\end{align}
Taking $\Phi\equiv 1$ in \eqref{eqn:7}, we get
\begin{align}\label{eqn:3}
\frac{\sum_{\mvt\in\bT_{\mvs}}|\mvA_k^{\ord}(\mvt)|}{|\bT_{\mvs}|\cdot s_0^k m^{k/2}}
\to \left(\frac{p_0\sigma}{2}\right)^k\E\bigg[\bigg(\int_{\cT_{2\ve}}\height(x)\ \mu_{\cT_{2\ve}}(dx)\bigg)^k\bigg].
\end{align}

Combining the above, we conclude that
\begin{align*}
\E\left[\Phi\bigg(\frac{1}{\sqrt{m}}\cQ\big(\widetilde\cT_{\mvs}, \widetilde\mvX\big)\bigg)\right]
&\to\frac{\E\bigg[\int_{x_1}\ldots\int_{x_{k+\ell}}\mu_{\cT_{2\ve}}^{\otimes k+\ell}(dx_1\ldots dx_{k+\ell})\prod_{i=1}^k\hght(x_i)\cdot
g_{\phi}^{\sss(k)}\bigg(\frac{1}{\sigma}\cT_{k,\ell}\bigg)\bigg]}{\E\bigg[\bigg(\int_{\cT_{2\ve}}\hght(x)\ \mu_{\cT_{2\ve}}(dx)\bigg)^k\bigg]}\\
&=\E\bigg[\Phi\bigg(\frac{1}{\sigma}M^{\sss(k)}\bigg)\bigg].
\end{align*}
This shows that
\begin{align}\label{eqn:19}
\frac{1}{\sqrt{m}}\cQ\big(\widetilde\cT_{\mvs}, \widetilde\mvX\big)
\weakc\frac{1}{\sigma}M^{\sss(k)}
\end{align}
with respect to Gromov-weak topology.

\chr{We will now improve this convergence to GHP convergence using Theorem \ref{thm:gw-to-ghp}.}
From the definition of $(\widetilde\cT_{\mvs}, \widetilde\mvX)$ (given right below \eqref{eqn:35}), it is clear that
\begin{align*}
\bP\big(\widetilde\cT_{\mvs}=\mvt\big)
=\frac{|\mvA_k(\mvt)|}{|\bT_{\mvs}^{\sss(k)}|}=\frac{|\mvA_k^{\ord}(\mvt)|}{\sum_{\mvt'\in\bT_{\mvs}}|\mvA_k^{\ord}(\mvt')|}
\end{align*}
for any $\mvt\in\bT_{\mvs}$. Hence for any bounded continuous (w.r.t. GHP topology) $h$,
\[\E\bigg[h\bigg(\frac{1}{\sqrt{m}}\widetilde\cT_{\mvs}\bigg)\bigg]
=\frac{\E\bigg[h\bigg(\frac{1}{\sqrt{m}}\cT_{\mvs}\bigg)\cdot \big|\mvA_k^{\ord}(\cT_{\mvs})\big|s_0^{-k} m^{-k/2}\bigg]}{\E\big[|\mvA_k^{\ord}(\cT_{\mvs})|s_0^{-k}m^{-k/2}\big]}.\]
Using Lemma \ref{lem:plane-trees}\eqref{lem:independent-sampling}, and Lemma \ref{lem:plane-trees}\eqref{lem:conditional-expectation}
together with uniform integrability (\chr{Lemma} \ref{lem:plane-trees}\eqref{lem:uniform-integrability}), we conclude that
\[\E\bigg[h\bigg(\frac{1}{\sqrt{m}}\widetilde\cT_{\mvs}\bigg)\bigg]
\to\E\bigg[h\bigg(\frac{1}{\sigma}\cT_{2\widetilde{\ve}_{(k)}}\bigg)\bigg],
\]
where $\widetilde{\ve}_{(k)}$ is as defined before \eqref{eqn:tilde-nu-k-def}. Hence $m^{-1/2}\widetilde\cT_{\mvs}\weakc\sigma^{-1}\cT_{2\widetilde{\ve}_{(k)}}$ in the GHP sense, and in particular, for each $\delta>0$, $1/\kappa_{\delta}\big(m^{-1/2}\widetilde\cT_{\mvs}\big)$, $m\geq 1$, is a tight sequence of random variables.
This immediately implies that $1/\kappa_{\delta}\big(m^{-1/2}\cQ(\widetilde\cT_{\mvs}, \widetilde\mvX)\big)$, $m\geq 1$, is also a tight sequence of random variables for each $\delta>0$. Combining this with \eqref{eqn:19} and Theorem \ref{thm:gw-to-ghp}, we see that
\begin{align*}
\frac{1}{\sqrt{m}}\cQ\big(\widetilde\cT_{\mvs}, \widetilde\mvX\big)
\weakc\frac{1}{\sigma}M^{\sss(k)}
\end{align*}
with respect to GHP topology. This concludes the proof of Theorem \ref{prop:condition-on-connectivity}. \qed

\subsection{Proof of Theorem \ref{thm:number-of-connected-graphs}}\label{sec:proof-cor-number-of-connected-graphs}
Recall the relation between $m$ and $\widetilde m$ from \eqref{eqn:def-m}.
Now it follows from \eqref{eqn:34} that
\[\big|\bG_{\widetilde\vd}^{\con}\big|=
\frac{\big|\bT_{\mvs}^{\sss(k)}\big|\times(s_0-2k)!\times\prod_{i\geq 1} s_i!}{\prod_{j=1}^{\widetilde m}(\widetilde d_j-1)!}.
\]
Further,
\[
\big|\bT_{\mvs}^{\sss(k)}\big|=\sum_{\mvt\in\bT_{\mvs}}|\mvA_k(\mvt)|
=\sum_{\mvt\in\bT_{\mvs}}|\mvA_k^{\ord}(\mvt)|/k!.
\]
We thus have
\begin{align*}
\big|\bG_{\widetilde\vd}^{\con}\big|
&=\frac{s_0^k m^{k/2}\times\big|\bT_{\mvs}\big|\times (s_0-2k)!\times\prod_{i\geq 1} s_i!}{\prod_{j=1}^{\widetilde m}(\widetilde d_j-1)!\times k!}
\times\sum_{\mvt\in\bT_{\mvs}}\frac{|\mvA_k^{\ord}(\mvt)|}{\big|\bT_{\mvs}\big|s_0^k m^{k/2}}\\
&\sim \frac{s_0^k m^{k/2}\times(\widetilde m+2k-2)!\times (s_0-2k)!}{s_0!\times\prod_{j=1}^{\widetilde m}(\widetilde d_j-1)!\times k!}\times
\left(\frac{p_0\sigma}{2}\right)^k\E\bigg[\bigg(\int_0^1 2\ve(x)dx\bigg)^k\bigg],
\end{align*}
where the last step uses \eqref{eqn:3} and the expression for $|\bT_{\mvs}|$ from \eqref{eqn:pitman}.
Using the relations $m/\widetilde m\sim 1$ and $s_0!/(s_0-2k)!\sim s_0^k(mp_0)^{k}$, a simple rearrangement of terms completes the proof.
\qed

\section{Proof of Theorems \ref{thm:graphs-given-degree-scaling} and \ref{thm:vacant-set-scaling}}\label{sec:proof-main-deg-vac-thm}

We start with the distribution of the configuration model.
\begin{lem}[\cite{Hofs13}, Proposition 7.7]\label{lem:CM-distribution}
Let $G$ be a multigraph on vertex set $[n]$ in which there are $x_{ij}$ many edges between $i$ and $j$, $1\leq i<j\leq n$, and vertex $i$ has $x_{ii}$ many loops. Let
$\overline d_i=x_{ii}+\sum_{j=1}^n x_{ij}$
be the total degree of $i$ (note that a loop contributes two to the degree). Let
\[\overline\vd=(\overline d_1,\ldots,\overline d_n),\ \text{ and }\ \ell_n=\sum_{i=1}^n\overline d_i.\]
Then
\[\overline\pr_{n,\overline\vd}\big(G\big)=
\frac{1}{(\ell_n-1)!!}\times
\frac{\prod_{i\in[n]}\overline d_i!}{\prod_{i\in[n]}2^{x_{ii}}\prod_{1\leq i\leq j\leq n}x_{ij}!}.\]
\end{lem}
We now state two fundamental results about the configuration model and uniform simple graphs with prescribed degree.

\begin{enumeratea}
	\item From Lemma \ref{lem:CM-distribution} (see also \cite{Boll-book,mcdiarmid1998concentration}), it follows that conditional on being simple, the configuration model has the same distribution as $\cG_{n,\vd}$, i.e.,
	\begin{equation}
	\label{eqn:cm-conditional-uniform}
		\bP\big(\CM_n(\vd) \in \cdot\ \big|\ \CM_n(\vd) \in \bG_{n, \vd}\big)= \bP\big(\cG_{n,\vd}\in \cdot\big).
	\end{equation}
    \item By \cite[Theorem 1.1]{janson2009probability}, under Assumption \ref{ass:cm-deg}, there exists a constant $c > 0$ such that the probability that $\CM_n(\vd)$ is simple satisfies
	\begin{equation}
	\label{eqn:cm-simple-0}
		\bP\big(\CM_n(\vd) \in \bG_{n, \vd}\big) \to c, \qquad \text{ as } n\to\infty.
	\end{equation}
\end{enumeratea}

This connection between $\CM_n(\vd)$ and $\cG_{n,\vd}$ is a very useful tool as it enables one to prove certain results about the uniform simple graph with given degrees by first obtaining a similar result for the configuration model, and then using \eqref{eqn:cm-conditional-uniform} and \eqref{eqn:cm-simple-0} to deduce the same for the simple graph.

For any nonnegative random variable $X$ with $\chr{0<\E X<\infty}$, define the corresponding size biased random variable $X^\circ$ via
\[\bP\big(X^\circ\leq x\big)=\frac{\E\big[X\ind_{X\leq x}\big]}{\E\big[X\big]},\ \ x\in[0,\infty).\]

\begin{prop}\label{prop:CM-facts}
Assume that $\vd$ satisfies Assumption \ref{ass:cm-deg} with limiting random variable $D$, and let $D^\circ$ denote the corresponding size-biased random variable. Let
\[p_i^\circ:=\bP\big(D^\circ=i\big)=\frac{i \bP\big(D=i\big)}{\E\big[D\big]},\ \ i=1, 2, \ldots.\]
\begin{enumeratei}
\item
Let $\cC_{\sss(k)}$ be the $k$-th largest component of $\CM_n(\vd)$. Then the following hold for each $k\geq 1$:
\begin{align}
&\hskip28pt\frac{1}{|\cC_{\sss(k)}|}\sum_{v\in \cC_{\sss(k)}} d_v^2\probc \sum_{i\ge 1}i^2 p_i^{\circ}<\infty;\label{eqn:cm-sum-square} \\
&\hskip52pt\pr\big(\cC_{\sss(k)} \text{ is simple }\big)\to 1;\label{eqn:cm-simple}\\
&\frac{1}{|\cC_{\sss(k)}|}\#\set{v\in \cC_{\sss(k)}: d_v=i}\probc  p^{\circ}_i\ \text{ for }\ i\ge 1.\label{eqn:cm-degree}
\end{align}
\item Further \eqref{eqn:cm-sum-square} 
and \eqref{eqn:cm-degree} continue to hold if we replace $\CM_n(\vd)$ by $\cG_{n,\vd}$.
\end{enumeratei}
\end{prop}
\noindent {\bf Proof:}
Given a sequence $a_1,\ldots,a_{\ell}$ of positive real numbers, the (random) size-biased permutation $\pi(1),\ldots, \pi(\ell)$ can be obtained as follows:
\[\bP\big(\pi(1)=i\big)=\frac{a_i}{\sum_{j=1}^{\ell}a_j},\ \text{ and }\
\bP\big(\pi(k)=i\ \big|\ \cS_{k-1}\big)=\frac{a_i \chr{\ind_{\{i\notin\cS_{k-1}\}}}}{\sum_{j\notin\cS_{k-1}}a_j},\ \ k=2,\ldots,\ell,\]
where $\cS_k=\{\pi(1),\ldots,\pi(k)\}$.
It is a standard fact that the random graph $\CM_n(\vd)$ can be explored in a depth-first way so that the vertices appear as a size-biased permutation, where vertex $i$ has size $d_i$; see \cite[Section 5.1]{dhara-hofstad-leeuwaarden-sen} or \cite{riordan2012phase}.
It further follows from \cite[Lemma 15]{dhara-hofstad-leeuwaarden-sen} and Theorem \ref{thm:cm-component-sizes} that for every $\eps>0$,
there exists $T_{\eps}>0$ such that
\begin{align}\label{eqn:36}
\limsup_n\ \bP\big(\cC_{(k)}\text{ is explored by time }T_{\eps}n^{2/3}\big)\geq 1-\eps.
\end{align}
\cite[Lemma 5]{dhara-hofstad-leeuwaarden-sen} shows that for every $T>0$,
\begin{align}\label{eqn:37}
\sup_{0\leq u\leq T}\bigg|\frac{1}{n^{2/3}}\sum_{i=1}^{\lfloor un^{2/3}\rfloor}d_{\pi(i)}^2-\frac{\sigma_3 u}{\sigma_1}\bigg|\probc 0,
\end{align}
where $\sigma_r=\E[D^r]$, $r=1,2,3$. Combining \eqref{eqn:36} and \eqref{eqn:37}, we get
\[\frac{1}{n^{2/3}}\bigg(\sum_{v\in \cC_{\sss(k)}} d_v^2-\frac{\sigma_3}{\sigma_1}\big|\cC_{(k)}\big|\bigg)\probc 0.\]
Since $\sigma_3/\sigma_1=\E\big[{D^\circ}^2\big]=\sum_{i\geq 1}i^2 p_i^\circ$, the last display together with Theorem \ref{thm:cm-component-sizes} yields \eqref{eqn:cm-sum-square}.


\eqref{eqn:cm-simple} follows from \cite[Section 5.3]{dhara-hofstad-leeuwaarden-sen} or by following verbatim the proof of this exact result but under slightly different moment assumptions in \cite[Equation 7.6]{dhara-hofstad-leeuwaarden-sen-2}. \eqref{eqn:cm-degree} follows from \cite[Equation 6.4]{dhara-hofstad-leeuwaarden-sen}.

Finally, part (ii) follows from (i) by an application of \eqref{eqn:cm-conditional-uniform} and \eqref{eqn:cm-simple-0}.
\qed

\vskip12pt

\noindent{\bf Proof of Theorem \ref{thm:graphs-given-degree-scaling}(i).}\ \
Note that $\bP\big(D=1\big)>0$ under Assumption \ref{ass:cm-deg}. Hence $p_1^\circ>0$.
Further, under Assumption \ref{ass:cm-deg},
\[\sum_{i\geq 1}i p_i^\circ=\E\big[D^\circ\big]=\E D^2/\E D=2.\]
Hence, by Proposition \ref{prop:CM-facts} (ii), for every $k\geq 1$, $(d_v,\ v\in\cC_{(k)})$ satisfies Assumption \ref{ass:degree} (after a possible relabeling of vertices) with limiting p.m.f. $(p_1^\circ, p_2^\circ,\ldots)$.

Let $\cP$ denote the partition of $\cG_{n,\vd}$ into different components. Then conditional on $\cP$, each component is uniformly distributed over the set of simple, connected graphs with the degrees prescribed by the partition $\cP$. Further, different components are conditionally independent.
We thus conclude using Theorem \ref{thm:cm-component-sizes} and Theorem \ref{prop:condition-on-connectivity} that for every $k\geq 1$,
\[n^{-2/3}\big(|\cC_{(1)}|,\ldots,|\cC_{(k)}|\big)
\weakc\big(|\gamma_{\sss(1)}^{\mvmu_{D}}(\lambda)|,\ldots,|\gamma_{\sss(k)}^{\mvmu_{D}}(\lambda)|\big)\]
jointly with
\[\bigg(\frac{1}{\sqrt{|\cC_{(1)}|}}\cC_{(1)},\ldots, \frac{1}{\sqrt{|\cC_{(k)}|}}\cC_{(k)}\bigg)
\weakc\frac{\alpha_D}{\sqrt{\eta_D}}\big(S_1,\ldots,S_k\big)\]
in the GHP sense, where $S_i$ are as in Construction \ref{constr:M-D}, and $\eta_D$, and $\alpha_D$ are as in Theorem \ref{thm:cm-component-sizes}. (Here we have used the fact $\sum_{i\geq 1}i^2 p_i^\circ-4=\eta_D/\alpha_D^2$.)
Combining the two yields the result.
\qed

\vskip12pt

\noindent{\bf Proof of Theorem \ref{thm:graphs-given-degree-scaling}(ii).}\ \
By Proposition \ref{prop:CM-facts} (i), for every $k\geq 1$, $(d_v,\ v\in\cC_{(k)})$ satisfies Assumption \ref{ass:degree} (after a possible relabeling of vertices) with limiting p.m.f. $(p_1^\circ, p_2^\circ,\ldots)$.
As before, let $\cP$ denote the partition of $\CM_n(\vd)$ into different components. For each $k\geq 1$, define the event
\[E_k:=\big\{\cC_{(j)}\ \text{ is simple for all }\ 1\leq j\leq k\big\}.\]
Then note that by Lemma \ref{lem:CM-distribution}, conditional on the event $E_k\cap\{\cP=P\}$, $\cC_{(j)}$, $j\geq 1$ are independent, and for each $i\leq k$, $\cC_{(i)}$ is uniformly distributed over the set of simple, connected graphs with the degrees prescribed by the partition $P$. Since $\bP(E_k^c)\to 0$ by \eqref{eqn:cm-simple}, the result follows by imitating the argument used in the proof of Theorem \ref{thm:graphs-given-degree-scaling}(i).\qed

\vskip12pt

\noindent{\bf Proof of Theorem \ref{thm:vacant-set-scaling}(ii).}\ \
For every $u\geq 0$, let $\cV^u$ denote the vacant set left by a random walk on $\CM_n(\vd_r^{\sss(n)})$ run up to time $nu$.
Let $\cE^u$ be the set of all edges of $\CM_n(\vd_r^{\sss(n)})$ both of whose endpoints are in $\cV^u$, i.e.,
\[\cE^u:=\big\{\{v_1,v_2\}\in\CM_n(\vd_r^{\sss(n)}) :\ v_1, v_2\in\cV^u\big\}.\]
Define the vacant graph $\mvV^u$ by $\mvV^u:=([n],\cE^u)$, and let $\vD^u:=\big(\vD^u(j),\ j\in[n]\big)$ be the degree sequence of $\mvV^u$.
Then, by \cite[Proposition 3.1]{cerny-teixeira}, for any collection $A$ of multigraphs on $[n]$,
\begin{align}\label{eqn:38}
\overline\vP_{n,r}\big(\mvV^u\in A\big)
=\sum_{\vd}\overline\vP_{n,r}\big(\vD^u=\vd\big)\times\overline\pr_{n,\vd}\big(A\big).
\end{align}
In words, the vacant graph $\mvV^u$ can be generated in two steps: (i) sample the degree sequence $\vD^u$ under the annealed measure $\overline\vP_{n,r}$, and then (ii) construct a configuration model with this degree sequence.

Let $\mvs^{u}=\big(s_0^{u},\ldots, s_r^{u}\big)$ denote the (random) frequency distribution corresponding to $\vD^{u}$, i.e., $s_i^{u}=\#\{j\in [n] : \vD^u(j)=i\}$, $0\leq i\leq r$.
Then, by \cite[Equation 6.1]{cerny-teixeira},
for every $\eps>0$,
\begin{align*}
\overline\vP_{n,r}\bigg(\bigg|\frac{1}{n}s_i^{u_{\star}}-\pr\big(D_{\vac}=i\big)\bigg|\geq\eps\bigg)
\to 0,\ \ 0\leq i\leq r,
\end{align*}
where $u_{\star}$ is as in \eqref{eqn:ustar-def}, and $D_{\vac}$ is as in \eqref{eqn:D-vac-def}. The simple observation
\[\big|s_i^{u_{\star}}-s_i^{u_n}\big|\leq (|a_0|+1)n^{2/3}(r+1)\]
for large $n$ when $u_n$ satisfies \eqref{eqn:40} leads to
\begin{align}\label{eqn:39}
\overline\vP_{n,r}\bigg(\bigg|\frac{1}{n}s_i^{u_n}-\pr\big(D_{\vac}=i\big)\bigg|\geq\eps\bigg)
\to 0,\ \ 0\leq i\leq r.
\end{align}
Further, by \cite[Equation 6.4]{cerny-teixeira},
\begin{align}\label{eqn:41}
\overline\vP_{n,r}\bigg(n^{1/3}
\bigg|\frac{\sum_{i=0}^r (i^2-2i)s_i^{u_n}}{\sum_{i=0}^r i s_i^{u_n}}
-\lambda_{\vac}\bigg|\geq\eps\bigg)\to 0,
\end{align}
where $\lambda_{\vac}$ is as in \eqref{eqn:lambda-vac-def}. Combining \eqref{eqn:39} and \eqref{eqn:41}, we see that
the degree sequence $\vD^{u_n}$ satisfies Assumption \ref{ass:cm-deg} with limiting random variable $D_{\vac}$ and $\lambda=\lambda_{\vac}$.
\chr{(That Assumption \ref{ass:cm-deg}(ii) is satisfied follows directly from the fact that $\vD^{u_n}(j)\leq r$ for all $j\in [n]$.)}
In view of \eqref{eqn:38}, an application of Theorem \ref{thm:graphs-given-degree-scaling}(ii) completes the proof.\qed

\vskip12pt

\noindent{\bf Proof of Theorem \ref{thm:vacant-set-scaling}(i).}\ \
Let $\mvV^u,\ \vD^u$, and $\mvs^u$ be as in the proof of Theorem \ref{thm:vacant-set-scaling}(ii), but with $\cG_{n,r}$ as the underlying graph (instead of $\CM_n(\vd_r^{\sss(n)})$). By \cite[Lemma 7]{cooper-frieze}, the analogue of \eqref{eqn:38} is true in this case, i.e.,
\[\vP_{n,r}\big(\mvV^u\in A\big)
=\sum_{\vd}\vP_{n,r}\big(\vD^u=\vd\big)\times\pr_{n,\vd}\big(A\big),\]
for any collection $A$ of simple graphs on $[n]$. Further, using \eqref{eqn:cm-conditional-uniform} and \eqref{eqn:cm-simple-0},
we conclude that  \eqref{eqn:39} and \eqref{eqn:41} continue to hold in this case. We complete the proof by an application of Theorem \ref{thm:graphs-given-degree-scaling}(i).\qed

\section*{Acknowledgements}
The authors thank Christina Goldschmidt and Remco van der Hofstad for many helpful discussions, and Bal\'{a}zs R\'{a}th for useful comments on a preliminary version of the paper. We thank two anonymous referees for detailed comments on the initial submission which lead to significant improvements in the exposition of the paper. 
SB has been partially supported by NSF-DMS grants 1105581, 1310002, 160683, 161307 and SES grant 1357622 and ARO grant W911NF-17-1-0010.
SS has been supported in part by EPSRC grant EP/J019496/1, a CRM-ISM fellowship, and the Netherlands Organization for Scientific Research (NWO) through the Gravitation Networks grant 024.002.003.

\bibliographystyle{plain}
\bibliography{vacant_bib}

\end{document}